\documentclass[11pt,reqno]{amsart}
\usepackage[T1]{fontenc}
\usepackage[utf8]{inputenc}
\usepackage{lmodern,amsmath,amssymb,amsthm,mathtools,mathrsfs}
\usepackage[margin=1.08in]{geometry}
\usepackage{microtype,enumitem,graphicx}
\usepackage[hidelinks,bookmarksnumbered]{hyperref}
\allowdisplaybreaks[2]
\numberwithin{equation}{section}
\newtheorem{theorem}{Theorem}[section]
\newtheorem{proposition}[theorem]{Proposition}
\newtheorem{lemma}[theorem]{Lemma}
\newtheorem{corollary}[theorem]{Corollary}
\theoremstyle{definition}
\newtheorem{definition}[theorem]{Definition}
\newtheorem{example}[theorem]{Example}
\theoremstyle{remark}
\newtheorem{remark}[theorem]{Remark}
\newcommand{\C}{\mathbb C}
\newcommand{\R}{\mathbb R}

\newcommand{\Sh}{\mathcal K}
\newcommand{\AllSh}{\mathcal S}
\newcommand{\Inc}{\mathcal I}
\newcommand{\Hull}{\mathsf H}
\newcommand{\RedQ}{\mathsf Q}
\newcommand{\RedN}{\mathsf N}
\newcommand{\dd}{\,\mathrm d}
\newcommand{\eps}{\varepsilon}
\newcommand{\Ree}{\operatorname{Re}}
\newcommand{\Imm}{\operatorname{Im}}
\newcommand{\supp}{\operatorname{supp}}
\newcommand{\conv}{\operatorname{conv}}
\newcommand{\dist}{\operatorname{dist}}
\newcommand{\divz}{\operatorname{div}_0}
\newcommand{\Area}{\operatorname{Area}}
\newcommand{\ind}{\mathbf 1}
\newcommand{\Haus}{\mathcal H}
\newcommand{\Rexp}{\R_{\exp}}
\newcommand{\centroid}{\mathfrak c}
\newcommand{\Blev}{\mathfrak b}
\newcommand{\FixSh}{\mathcal K^{\mathrm{fix}}}
\newcommand{\rootmu}{\mu^{\mathrm{root}}}
\newcommand{\rootC}{\mathscr C}
\newcommand{\mom}{\mathfrak m}
\newcommand{\LogDer}{\Lambda}
\newcommand{\prz}{\operatorname{pr}_z}
\DeclareMathOperator{\wind}{wind}
\DeclareMathOperator{\Law}{Law}
\DeclareMathOperator{\Unif}{Unif}
\DeclareMathOperator{\acc}{acc}
\hypersetup{pdftitle={The full shadow of a polynomial III. Semialgebraicity and reflection rigidity},
 pdfauthor={Christian Hagg and Boris Shapiro},
 pdfkeywords={Polynomial shadows, semialgebraicity, reflection rigidity, Cauchy transforms, moment reconstruction}}
\title[The full shadow of a polynomial III]{The full shadow of a polynomial III. Semialgebraicity and reflection rigidity}
\author{Christian H\"agg}
\address{Department of Mathematics, Stockholm University, SE--106 91 Stockholm, Sweden}
\email{hagg@math.su.se}
\author{Boris Shapiro}
\address{Department of Mathematics, Stockholm University, SE--106 91 Stockholm, Sweden}
\email{shapiro@math.su.se}
\date{}
\begin{document}
\begin{abstract}
For a polynomial $P$, the zeros of all derivatives of the powers $P^n$
converge, as $n\to\infty$, to a compact set, the limiting core of $P$, and
their normalized counting measures converge to a canonical probability
measure supported on it. We prove that in each degree these cores form a
semialgebraic family, although their definition is transcendental. The
limiting measure splits into a planar density, a measure carried by
reflection axes of the root configuration, and explicit atoms at the roots.
This rigidity classifies all cores of zero area: they are segments or finite
unions of segments with a common endpoint. A conservation law for the Cauchy
transform identifies the holomorphic moments of the limiting measure with
those of random chords between roots, and reconstructs a polynomial of
degree $d$ from the first $d$ moments. We also establish analytic regularity
of the planar density and stability of the area, including convergence when
root masses vanish.
\end{abstract}
\maketitle
\setcounter{tocdepth}{1}
\begingroup\small\tableofcontents\endgroup

\section{Introduction}\label{sec:intro}
By the Gauss--Lucas theorem, the zeros of every derivative of a polynomial
lie in the convex hull of its zeros. This paper is concerned with how these
zeros fill the hull when the degree and the order of differentiation are
both large. We fix a polynomial $P$ of degree $d$ and consider the zeros of
all derivatives of all its powers, that is, of the polynomials
\[
 (P^n)^{(m)},\qquad n\ge1,\quad 0\le m<dn.
\]
For $P(z)=z^2-1$ and $m=n$ these are, by Rodrigues' formula and up to
normalization, the Legendre polynomials; for general $P$ they are called
\emph{Rodrigues' descendants} of $P$ \cite{BHS}. The closure of the set of
all their zeros is the \emph{full shadow} of $P$. The asymptotic
distribution of zeros under repeated differentiation has attracted
considerable attention, usually at a fixed ratio between the order of
differentiation and the degree; see \cite{BHS,MFR} and the references given
there. Here all orders of differentiation are considered simultaneously.

The first two parts of this series \cite{HS1,HS2} established the basic
structure of the full shadow. As $n\to\infty$, the set of all zeros of the
descendants $(P^n)^{(m)}$, $0\le m<dn$, converges in Hausdorff distance to a
compact connected set $\Sh_P$, the \emph{limiting core} of $P$, and the
normalized counting measures of these zeros converge weakly to a probability
measure $\rho_P$ with support $\Sh_P$, the \emph{shadow measure}. The full shadow consists of the
core and of isolated zeros outside it, of which there may be infinitely
many. Both limits are described by logarithmic potential theory. If $m/n$
tends to a ratio $\alpha\in(0,d)$, the zeros of $(P^n)^{(m)}$ equidistribute
according to a measure whose potential at a point $z$ is, up to an affine
normalization, the supremum of the levels $\ell$ for which $z$ can be joined
to infinity within the superlevel set
$\{w:\log|P(w)|-\alpha\log|w-z|\ge\ell\}$; the shadow measure is an average
of these measures over $\alpha$. We refer to this extremal problem as the
\emph{logarithmic separation problem}. Its solution is transcendental in
nature: to decide whether a given point belongs to $\Sh_P$, one has to
compare the values of $\log|P(w)|-\alpha\log|w-z|$ at different critical
points $w$. A priori, the core is therefore only known to be definable in
the o-minimal structure generated by the real exponential function
\cite{HS2}. This leaves two questions open: how complicated can the core
be, and what does the shadow measure look like?

We show that, in spite of this description, the core is an algebraic object
and the shadow measure has a rigid structure. First, the cores of the monic
polynomials of a given degree $d$ form a semialgebraic family, and every
boundary $\partial\Sh_P$ lies on a real algebraic curve whose degree is
bounded in terms of $d$ alone. Second, $\rho_P$ is the sum of an absolutely
continuous planar part, a part with a density along finitely many arcs, and
an atom at each root. The middle part can only live on lines of mirror
symmetry of the root configuration, so it is absent for a generic polynomial
of degree at least three. As a consequence, we classify the polynomials
whose core has zero area: up to an affine change of variable they are the
polynomials $z^{r_0}H(z^q)$ with $H$ real-rooted, and their shadows are
segments or finite unions of segments with a common endpoint.

Third, we prove a conservation law for the Cauchy transform $C_{\rho_P}$ of
the shadow measure: outside the core, $C_{\rho_P}'=-(P'/(dP))^2$. Seen from
outside the convex hull of the roots, $\rho_P$ therefore cannot be
distinguished from the distribution of a uniformly random point on the chord
joining two roots chosen independently, with probabilities proportional to
their multiplicities. This gives all holomorphic moments of $\rho_P$ in
closed form and shows that $\rho_P$ determines $P$ up to passing to powers;
it is also a main tool in the proofs of the geometric results. Finally, we
prove that the planar density of $\rho_P$ is real-analytic and positive away
from finitely many arcs and points, and that the area of the shadow is
stable under perturbations of $P$, also when root masses tend to zero.

The remainder of this introduction gives the precise definitions and
statements and then discusses earlier work and the strategy of the proofs.
Figure~\ref{fig:regimes} shows three typical shadows.

\subsection{Shadows and shadow measures}
For a nonzero polynomial $G$, let $Z(G)$ be the set of its distinct zeros and
let $\divz G=\sum_{a\in Z(G)}\operatorname{ord}_a(G)\delta_a$ be its
zero-counting measure. Throughout the paper, $P$ denotes a monic polynomial
\begin{equation}\label{eq:root-data}
 P(z)=\prod_{j=1}^s(z-a_j)^{r_j},\qquad
 d=\sum_{j=1}^s r_j,\qquad p_j=\frac{r_j}{d},
\end{equation}
where the roots $a_j$ are distinct and the multiplicities $r_j$ are positive
integers; multiplying $P$ by a nonzero constant changes none of the zero
sets and measures considered below. We write
\[
 \rootmu_P=\sum_{j=1}^s p_j\delta_{a_j},\qquad
 \Hull_P=\conv Z(P),\qquad \centroid_P=\sum_{j=1}^s p_ja_j
\]
for the \emph{normalized root measure} of $P$, the convex hull of its roots,
and their centroid.
For a finite, compactly supported measure $\nu$ on $\C$, we normalize the
logarithmic potential and the Cauchy transform by
\begin{equation}\label{eq:potential-conventions}
 U_\nu(z)=\int\log|z-\xi|\dd\nu(\xi),\qquad
 C_\nu(z)=\int\frac{\dd\nu(\xi)}{z-\xi}.
\end{equation}
In particular, the \emph{root Cauchy transform} and the \emph{root moments}
of $P$ are
\[
 \rootC_P(z):=C_{\rootmu_P}(z)=\frac{P'(z)}{dP(z)},\qquad
 \mom_k(P):=\int \xi^k\dd\rootmu_P(\xi).
\]

Set $R_{m,n}=(P^n)^{(m)}$ for $n\ge1$ and $0\le m<dn$. These polynomials are
the \emph{descendants} of $P$, and we call $n$ their \emph{power level}. Define
\[
 \mathcal Z_{P,n}=\bigcup_{m=0}^{dn-1}Z(R_{m,n}),\qquad
 \AllSh_P=\overline{\bigcup_{n\ge1}\mathcal Z_{P,n}}.
\]
We call $\AllSh_P$ the \emph{full shadow} of $P$. If $m/n\to\alpha\in(0,d)$,
the probability measures $(dn-m)^{-1}\divz R_{m,n}$ converge weakly to a
probability measure $\mu_{\alpha,P}$ \cite[Theorem~3.3]{HS1}, which we call
the \emph{fixed-ratio measure} at ratio $\alpha$. Define
\begin{equation}\label{eq:core-measure}
 \FixSh_{\alpha,P}=\supp\mu_{\alpha,P},\qquad
 \rho_P=\frac2{d^2}\int_0^d(d-\alpha)\mu_{\alpha,P}\dd\alpha,
 \qquad \Sh_P=\supp\rho_P.
\end{equation}
The set $\FixSh_{\alpha,P}$ is the \emph{fixed-ratio support}, the measure
$\rho_P$ is the \emph{shadow measure}, and $\Sh_P$ is the \emph{limiting
core}, or simply the \emph{core}, of $P$. The weight $d-\alpha$ reflects
the number $dn-m$ of zeros of $R_{m,n}$: equivalently, $\rho_P$ is the weak
limit of the pooled measures
\begin{equation}\label{eq:pooled-powers}
 \bar\nu_{P,n}:=\frac{2}{dn(dn+1)}
       \sum_{m=0}^{dn-1}\divz R_{m,n},
\end{equation}
the normalized counting measures of all zeros of the descendants at power
level $n$. The zero sets $\mathcal Z_{P,n}$ themselves converge to $\Sh_P$ in
Hausdorff distance. The remaining points
\[
 T_P:=\AllSh_P\setminus\Sh_P
\]
of the full shadow are isolated zeros of descendants, each occurring at only
finitely many power levels, and they form a locally finite subset of
$\C\setminus\Sh_P$; we call them \emph{transients}. Thus the core and the
full shadow have the same interior and area
\cite[Theorems~5.1, 5.3, and~5.7]{HS1}. We refer to $\Sh_P$ and $\AllSh_P$ as
the two \emph{shadows} of $P$. The difference between them matters for the
set $\AllSh_P$ itself, but not for its planar geometry.

We use $\dd A$ for planar Lebesgue measure, $\Haus^1$ for arclength,
and $d_{\mathrm H}$ for Hausdorff distance. The symbol $\Longrightarrow$
denotes weak convergence of measures, convergence of compact sets is always
understood in Hausdorff distance, and $[a,b]$ denotes the closed segment
with endpoints $a,b\in\C$. The ratio $\alpha=m/n$ compares the order of
differentiation with the exponent $n$ of $P$. The ratio relative to the
total degree is $t=m/(dn)=\alpha/d$; it will be used for general atomic
measures in Sections~\ref{sec:weighted} and~\ref{sec:area-stability}.

\subsection{Algebraic geometry and reflection rigidity}
A subset of $\R^N$ is \emph{semialgebraic} if it is a finite Boolean
combination of sets defined by real polynomial equalities and inequalities.
The space of monic polynomials of degree $d$ is identified with $\R^{2d}$
through the real and imaginary parts of the coefficients. A
\emph{reflection axis} of $P$ is a line whose reflection preserves the
\emph{root multiset} of $P$, that is, its roots together with their
multiplicities.

\begin{theorem}[Geometry of the core]\label{thm:intro-tame}
For each degree $d$, the family $\{(P,z):z\in\Sh_P\}$, where $P$ ranges over
the monic polynomials of degree $d$, is semialgebraic.
There is an integer $D_d$ such that every boundary $\partial\Sh_P$ lies
in a real algebraic curve of degree at most $D_d$. The cores in each degree
have finitely many homeomorphism types; each core is homeomorphic to a
finite simplicial complex, has finitely many complementary components, and
has a boundary of finite length.

The shadow measure has a decomposition
\begin{equation}\label{eq:decomposition-intro}
 \rho_P=f_P\dd A+
 \vartheta_P\,\Haus^1\!\restriction\Gamma_P+
 \sum_{j=1}^s p_j^2\delta_{a_j},
\end{equation}
where the densities $f_P$ and $\vartheta_P$ are nonnegative and integrable
and $\Gamma_P$ is a finite union of bounded $C^1$ arcs of finite length. The
second term is supported on the reflection axes of $P$.

If $s\ge2$, put
\[
 \mathcal B_P^{\mathrm{pl}}=\overline{\operatorname{int}\Sh_P},\qquad
 \mathcal L_P=\overline{\Sh_P\setminus\mathcal B_P^{\mathrm{pl}}}.
\]
Then
\begin{equation}\label{eq:full-body-intro}
 \AllSh_P=\mathcal B_P^{\mathrm{pl}}\cup\mathcal L_P\cup T_P,
 \qquad \mathcal B_P^{\mathrm{pl}}=\supp(f_P\dd A).
\end{equation}
The set $\mathcal L_P$ is empty or a finite union of straight segments
on reflection axes, and its intersection with $\mathcal B_P^{\mathrm{pl}}$
is finite.
\end{theorem}

In \eqref{eq:decomposition-intro}, the first term is the \emph{planar part}
of $\rho_P$ and the second term is its \emph{curve part}; the mass of the
latter is the \emph{curve mass}. We call $\mathcal B_P^{\mathrm{pl}}$ the
\emph{planar body} and $\mathcal L_P$ the \emph{axial set} of $P$. If $P$
has no reflection axis, its core equals its planar body and the only
singular part of $\rho_P$ consists of the root atoms. More generally,
reflection rigidity gives a complete classification of the polynomials whose
shadows have zero area. A finite union of segments with a common endpoint is
called a \emph{radial star}.

\begin{theorem}[Zero-area classification]\label{thm:p3-dimension-summary}
Suppose that $P$ has at least two distinct roots. The core has zero area
if and only if, after a change $P(z)\mapsto cP(az+b)$ with $ac\ne0$, the
polynomial has the form
\[
 z^{r_0}H(z^q),\qquad q\ge1,\quad r_0\ge0,
\]
where $H$ is a nonconstant real-rooted polynomial with $H(0)\ne0$.
In this case the full shadow equals the core and is a segment or a radial
star. In all other cases the core has nonempty interior.
\end{theorem}

\begin{figure}[t]
\centering
\includegraphics{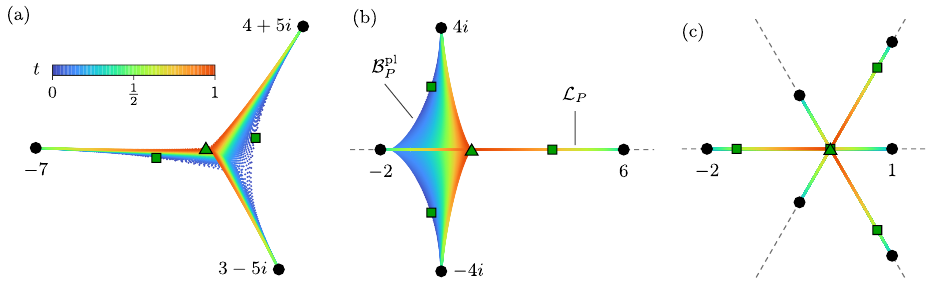}
\caption{Zeros of $(P^n)^{(m)}$, $0\le m<dn$, colored by $t=m/(dn)$.
Black disks: roots; green squares: critical points; green triangle:
$\centroid_P$; dashed: reflection axes.
(a)~$P=(z+7)(z-3+5i)(z-4-5i)$, $n=110$: no reflection axis, and
$\Sh_P=\overline{\operatorname{int}\Sh_P}$.
(b)~$P=(z^2+16)(z+2)(z-6)$, $n=120$: the planar body
$\mathcal B_P^{\mathrm{pl}}$ and the two segments of $\mathcal L_P$ on
the real axis.
(c)~$P=(z^3-1)(z^3+8)$, $n=60$: every root lies on an axis, and every
zero lies on a radial star.}
\label{fig:regimes}
\end{figure}

In the zero-area case the shadow measure is given by explicit sampling
formulas (Theorems~\ref{thm:collinear} and~\ref{thm:radial}); negative roots
of $H$ are allowed and produce a second family of segments, interlaced with
the first. A finite test on the coefficients detects the reflection axes and
decides the zero-area alternative (Proposition~\ref{prop:coefficient-test}).
Figure~\ref{fig:regimes} illustrates the three possibilities: no reflection
axis, an axis off which some roots lie, and every root on an axis.

\subsection{Conservation and reconstruction}
The analytic identity behind several of these results is the
\emph{conservation law}
\begin{equation}\label{eq:Cprime-intro}
 C_{\rho_P}'(z)=-\rootC_P(z)^2,\qquad z\in\C\setminus\Sh_P.
\end{equation}
It holds on every component of the complement of the core, including the
bounded ones. Near infinity it follows from an exact identity, valid for
every polynomial, between the logarithmic derivatives of its successive
derivatives (Proposition~\ref{prop:finite-chord}); the characterization of
the fixed-ratio measures by their potentials then gives
\eqref{eq:Cprime-intro} throughout the complement
(Theorem~\ref{thm:Cprime}).

The conservation law has a probabilistic interpretation. Let $X,Y$ be
independent random variables with law $\rootmu_P$, and let $\mathsf U$ be
uniform on $[0,1]$ and independent of them. The \emph{chord measure}
\[
 \sigma_P=\Law((1-\mathsf U)X+\mathsf UY)
\]
has the same logarithmic potential as $\rho_P$ outside $\Hull_P$
(Theorem~\ref{thm:chord}). Consequently, the \emph{holomorphic moments} of
$\rho_P$ are
\begin{equation}\label{eq:moments-intro}
 \int z^k\dd\rho_P(z)
   =\frac1{k+1}\sum_{j=0}^k\mom_j(P)\mom_{k-j}(P),\qquad k\ge0.
\end{equation}
When the roots are not collinear, this identifies the holomorphic moments of
$\rho_P$ and $\sigma_P$, rather than the measures themselves.

\begin{theorem}[Finite-moment reconstruction]\label{thm:p3-inverse-summary}
The first $d$ holomorphic moments of $\rho_P$ determine the monic
polynomial $P$ of degree $d$. Explicitly, if
$\mathfrak s_j=\int z^j\dd\rho_P(z)$, then $\mom_0=1$ and
\[
 \mom_j=\frac12\left((j+1)\mathfrak s_j-
           \sum_{i=1}^{j-1}\mom_i\mom_{j-i}\right),\qquad 1\le j\le d,
\]
recover the root moments, and Newton's identities recover $P$.
The Cauchy transform of $\rho_P$ near infinity also determines the
normalized root measure through $\rootC_P=\sqrt{-C_{\rho_P}'}$, with the
branch asymptotic to $z^{-1}$.
\end{theorem}

If the whole measure $\rho_P$ is known, the reconstruction is simpler and
does not require the degree: the atoms of $\rho_P$ have masses $p_j^2$, so
their positive square roots recover $\rootmu_P$. Thus two monic polynomials
have the same shadow measure precisely when they are powers of a common
monic polynomial whose root multiplicities have greatest common divisor
one. Such a polynomial will be called \emph{primitive}. Passing from a
polynomial to one of its powers can change the transient set, while it
leaves the core and the shadow measure unchanged.

\subsection{Earlier work and outline of the proofs}
B\o gvad, H\"agg, and Shapiro \cite{BHS} studied the descendants of $P$ at a
fixed differentiation ratio and proved limit theorems for the distribution
of their zeros when $P$ and $P'$ have simple zeros. Their shadow-boundary
conjecture \cite[Conjecture~7.2]{BHS}, stated in degree at least three under
these simplicity conditions and for roots in convex position that do not
form a regular polygon, is expressed through critical values of the rational
map $w\mapsto w-\alpha P(w)/P'(w)$. Our description of the boundary
(Theorem~\ref{thm:global-boundary}) includes these critical values together
with a second alternative, which arises from pairs of inverse images with
equal values of both $|P|$ and $|P'|$. We thus obtain algebraic containment
of the boundary and semialgebraicity of the entire coefficient family; the
description of the boundary by these critical values alone and the concavity
asserted in that conjecture are separate questions.

Mart\'inez-Finkelshtein and Rakhmanov \cite{MFR} established exterior
Cauchy-transform limits for iterated differentiation of general polynomial
arrays and related the evolution to a Burgers equation and fractional
free convolution. Here, averaging over all orders of differentiation yields
an exact identity in finite degree and a quantitative exterior estimate
(Theorem~\ref{thm:universal-exterior-chord}). We also establish the
conservation law on the bounded components of the complement of the core,
which is essential for the continuation argument of
Section~\ref{sec:algebraic}.

The paper is organized as follows. Section~\ref{sec:incidence-selection}
recalls the logarithmic separation problem, the fixed-ratio measures, and
the continuity results that we need from \cite{HS1,HS2}.
Section~\ref{sec:chord} proves the conservation law and the reconstruction
from moments. In Section~\ref{sec:balance}, local changes between the
harmonic branches of the potential give the decomposition of the shadow
measure and the alternatives for boundary points. In
Section~\ref{sec:algebraic}, the conservation law forces the core to contain
each edge of a finite algebraic graph entirely or not at all, and continuity
in the coefficients upgrades this to semialgebraicity of the family.
Reflection rigidity (Section~\ref{sec:rigidity}) and the explicit segment
and radial models (Section~\ref{sec:models}) yield the zero-area
classification. Section~\ref{sec:planar-stability} establishes regularity of
the planar density, and Sections~\ref{sec:weighted}
and~\ref{sec:area-stability} extend the arguments to finite atomic measures,
including vanishing masses and convergence of areas. Algebraic selection at
rational ratios and the potential-theoretic effect of transients are treated
in the appendices.

\section{Potentials, selected saddles, and support incidence}
\label{sec:incidence-selection}

\subsection{The logarithmic separation problem}
With the conventions \eqref{eq:potential-conventions}, we have
$2\partial_zU_\nu=C_\nu$, $\Delta U_\nu=2\pi\nu$, and
$\partial_{\bar z}C_\nu=\pi\nu$ in the sense of distributions, where
$\partial_z=(\partial_x-i\partial_y)/2$ and
$\partial_{\bar z}=(\partial_x+i\partial_y)/2$.
We write $C_{\alpha,P}:=C_{\mu_{\alpha,P}}$ and introduce the polynomials
$\RedQ_P,\RedN_P$ and the rational functions $\LogDer_P,h_P,F_{P,\alpha}$ by
\begin{equation}\label{eq:reduced-data}
 \begin{split}
 \RedQ_P(w)&=\prod_{j=1}^s(w-a_j),\qquad
 \RedN_P(w)=\sum_{j=1}^s r_j\prod_{k\ne j}(w-a_k),\\
 \LogDer_P(w)&=\frac{P'(w)}{P(w)}=\frac{\RedN_P(w)}{\RedQ_P(w)},\qquad
 h_P(w)=\frac{\RedQ_P(w)}{\RedN_P(w)},\qquad
 F_{P,\alpha}(w)=w-\alpha h_P(w).
 \end{split}
\end{equation}
The polynomials $\RedQ_P,\RedN_P$ are coprime. In particular,
$h_P=P/P'$ is understood as a reduced rational function.

Fix an \emph{observation point} $z\notin Z(P)$ and a ratio $0<\alpha<d$,
and consider the real \emph{landscape}
\[
 L_{P,z,\alpha}(w)=\log|P(w)|-\alpha\log|w-z|
\]
as a function of $w$. We call the $z$-plane the \emph{observation plane};
the adjective \emph{spatial} always refers to the variable $z$, as opposed
to the ratio. The landscape tends to $+\infty$ at $w=z$ and at $w=\infty$,
which we call the two \emph{marked ends}, and to $-\infty$ at the roots. For
a real level $\ell$, adjoin the marked ends to the superlevel set:
\[
 E_\ell(P,z,\alpha)
 =\{w\in\C\setminus(Z(P)\cup\{z\}):L_{P,z,\alpha}(w)\ge\ell\}
       \cup\{z,\infty\}\subset\widehat\C,
\]
where $\widehat\C$ is the Riemann sphere. The \emph{barrier} is
\begin{equation}\label{eq:barrier-definition}
 \Blev_P(z,\alpha)=
 \sup\{\ell:z\text{ and }\infty\text{ belong to the same connected
 component of }E_\ell(P,z,\alpha)\}.
\end{equation}
The \emph{logarithmic separation problem} is the determination of this
level. At the barrier level, let $\mathcal U_z$ and $\mathcal U_\infty$ be
the components of the strict superlevel set
$\{L_{P,z,\alpha}>\Blev_P(z,\alpha)\}\cup\{z,\infty\}$ containing
the respective marked ends.

\begin{figure}[t]
\centering
\includegraphics{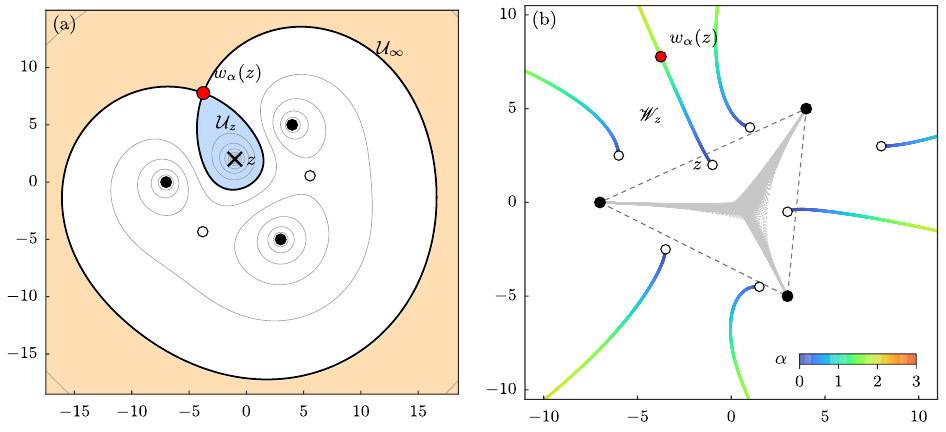}
\caption{(a)~The landscape $L_{P,z,\alpha}$ for the cubic of
Figure~\ref{fig:regimes}(a), $z=-1+2i$ (cross) and $\alpha=\tfrac32$.
Blue and orange: the components $\mathcal U_z$ and $\mathcal U_\infty$
of $\{L_{P,z,\alpha}>\Blev_P(z,\alpha)\}$, which meet at the connecting
saddle $w_\alpha(z)$ (red); black: the level set at height
$\Blev_P(z,\alpha)$; white disks: the other two saddles.
(b)~Selected saddle paths $\mathscr W_z$ of Theorem~\ref{thm:Cprime},
colored by $\alpha$, from seven points $z\notin\Sh_P$ (white disks);
gray: the zeros of Figure~\ref{fig:regimes}(a); dashed:
$\partial\Hull_P$.}
\label{fig:landscape}
\end{figure}

An \emph{ordinary saddle}, or simply a \emph{saddle}, is a critical
point of $L_{P,z,\alpha}$ away from the roots and marked ends; the value of
$L_{P,z,\alpha}$ at a saddle is its \emph{critical height}. Equivalently, a
saddle is a solution of
\begin{equation}\label{eq:reduced-saddle}
 (w-z)\RedN_P(w)-\alpha\RedQ_P(w)=0,
 \qquad F_{P,\alpha}(w)=z.
\end{equation}
The first of these, the \emph{reduced saddle equation}, is a polynomial
equation in $w$ of degree $s$ with leading coefficient $d-\alpha$.
A saddle is \emph{connecting}, or \emph{selected}, precisely when it
belongs to $\overline{\mathcal U_z}\cap\overline{\mathcal U_\infty}$;
see Figure~\ref{fig:landscape}(a).
It is \emph{simple} when $F_{P,\alpha}'(w)\ne0$; a solution with
$F_{P,\alpha}'(w)=0$ is a \emph{saddle collision}. With local branches of
the logarithms, the \emph{holomorphic phase} is
\[
 \Phi_{P,z,\alpha}(w)=\log P(w)-\alpha\log(w-z),\qquad
 \Ree\Phi_{P,z,\alpha}=L_{P,z,\alpha}.
\]
The \emph{order} of a saddle $w$ is the least $\kappa\ge2$ for which
$\partial_w^\kappa\Phi_{P,z,\alpha}(w)\ne0$; a simple saddle has order two.

The next proposition collects the facts about separating contours and
selection that are used below.
\begin{proposition}[Separating contours and local selection]
\label{prop:selection-facts}
\cite[Lemmas~2.1, 2.5--2.7, and~2.9]{HS1}
For $z\notin Z(P)$ and $0<\alpha<d$:
\begin{enumerate}[label=\textup{(\roman*)}]
\item The barrier equals
\[
 \inf_\gamma\max_{w\in\gamma}L_{P,z,\alpha}(w),
\]
where $\gamma$ ranges over the \emph{separating contours}, that is, the
Jordan curves avoiding $Z(P)$ and having winding number one about $z$.
There is a separating contour passing
through the connecting saddles and lying strictly below the barrier
elsewhere. This statement allows saddles of arbitrary finite order.
At a unique simple connecting saddle, the landscape has a nondegenerate
maximum along this contour and a strict minimum along a path joining the two
marked ends through that saddle.
\item If all connecting saddles at a parameter value are simple, their
local analytic continuations $w_1,\ldots,w_k$ satisfy
\begin{equation}\label{eq:local-barrier-maximum}
 \Blev_P(z,\alpha)=\max_{1\le j\le k}
                  L_{P,z,\alpha}(w_j(z,\alpha))
\end{equation}
in a neighborhood of that value. This holds jointly with the polynomial
coefficients, even if other saddles are multiple. The corresponding
contours and paths persist with uniform strict height gaps on compact
pieces away from their selected saddles.
\item On a \emph{root-free} disk in the observation plane (a disk
containing no root of $P$) on which
$\Blev_P(\cdot,\alpha)$ is harmonic, the connecting saddle is unique
and simple. It varies holomorphically in $z$, and the separating contour
through this single saddle and its strict height gap persist under small
changes of $z$, $\alpha$, and the coefficients.
\end{enumerate}
\end{proposition}

\subsection{Potentials and limits of zero sets}
The normalization of the potentials involves the entropy term
\[
 \mathfrak h_d(\alpha)
 =d\log d-\alpha\log\alpha-(d-\alpha)\log(d-\alpha).
\]
We record the results on potentials and on supports together, since the
arguments use both the weak limits of the zero-counting measures and the
Hausdorff limits of the zero sets.

\begin{proposition}[Fixed-ratio measures and selected inverses]
\label{prop:fixed-facts}
\cite[Theorem~3.3, Propositions~3.4--3.5, Lemma~3.6, and Theorem~3.9]{HS1}
The measure $\mu_{\alpha,P}$ is a probability measure supported in
$\Hull_P$, and its potential is
\begin{equation}\label{eq:barrier-potential}
 U_{\mu_{\alpha,P}}(z)
    =\frac{\Blev_P(z,\alpha)-\mathfrak h_d(\alpha)}{d-\alpha}
       \qquad(z\notin Z(P)),
\end{equation}
with the subharmonic extension at the roots. Its atomic part is
\begin{equation}\label{eq:fixed-atoms}
 \sum_{j=1}^s\frac{(r_j-\alpha)_+}{d-\alpha}\,\delta_{a_j},
 \qquad x_+:=\max\{x,0\}.
\end{equation}
The measures vary weakly continuously, and their supports vary continuously
in Hausdorff distance, jointly in $P$ and $\alpha\in(0,d)$. The zero sets
$Z(R_{m,n})$ converge to these supports whenever $m/n\to\alpha$;
the convergence is locally uniform in these parameters.

On $\C\setminus\FixSh_{\alpha,P}$ the function $C_{\alpha,P}$ has
no zeros, and
\begin{equation}\label{eq:selected-inverse}
 \begin{split}
 w_\alpha(z)&=z+\frac{\alpha}{(d-\alpha)C_{\alpha,P}(z)},\\
 F_{P,\alpha}(w_\alpha(z))&=z,\qquad
 \LogDer_P(w_\alpha(z))=(d-\alpha)C_{\alpha,P}(z).
 \end{split}
\end{equation}
At an observation point that is not a root, $w_\alpha(z)$ is the unique
selected simple saddle; we call $w_\alpha$ the \emph{selected inverse}
branch of $F_{P,\alpha}$, or simply the selected inverse.
It avoids the roots of $P$ and the poles of $h_P$, is holomorphic in $z$,
and is locally real-analytic in $\alpha$.
\end{proposition}

\begin{proposition}[Endpoints, pooling, and the full shadow]
\label{prop:global-facts}
\cite[Theorems~4.5, 5.1, 5.3, 5.6, and~5.7; Corollaries~5.4 and~5.12]{HS1}
As $\alpha\downarrow0$ and $\alpha\uparrow d$, respectively,
\begin{align*}
 \mu_{\alpha,P}&\Longrightarrow\rootmu_P,
 &\FixSh_{\alpha,P}&\longrightarrow Z(P)\cup Z(P'),\\
 \mu_{\alpha,P}&\Longrightarrow\delta_{\centroid_P},
 &\FixSh_{\alpha,P}&\longrightarrow\{\centroid_P\},
\end{align*}
where the support limits are in Hausdorff distance. The pooled measures
in \eqref{eq:pooled-powers} converge weakly to $\rho_P$, and
\begin{equation}\label{eq:core-union}
 \Sh_P=\overline{\bigcup_{0<\alpha<d}\FixSh_{\alpha,P}},
 \qquad \mathcal Z_{P,n}\longrightarrow\Sh_P.
\end{equation}
The latter convergence is uniform for $P$ in compact sets of monic
polynomials of a fixed degree, and $P\mapsto\Sh_P$ is continuous in
Hausdorff distance. The core is connected
and contains $Z(P)\cup Z(P')\cup\{\centroid_P\}$. If $s\ge2$, it has
no isolated points and
\[
 \AllSh_P=\Sh_P\,\dot\cup\,T_P,
 \quad \acc(\AllSh_P)=\Sh_P,
 \quad \operatorname{int}\AllSh_P=\operatorname{int}\Sh_P,
 \quad \partial\AllSh_P=\partial\Sh_P\cup T_P,
\]
where $\acc(E)$ denotes the set of limits of sequences of distinct points
of $E$, and $T_P$ is countable and locally finite in $\C\setminus\Sh_P$.
For $s=1$, both shadows reduce to the single root.
\end{proposition}

Integrating \eqref{eq:fixed-atoms} against the weight in
\eqref{eq:core-measure} gives the atomic part of the shadow measure:
\begin{equation}\label{eq:shadow-atoms}
 \rho_P^{\mathrm{at}}=\sum_{j=1}^s p_j^2\delta_{a_j},\qquad
 \rho_P(\{a_j\})=\frac2{d^2}\int_0^{r_j}(r_j-\alpha)\dd\alpha.
\end{equation}
The shadow measure and the core are unchanged when $P$ is replaced by a
positive integer power. They are covariant under affine changes of variables:
if $\phi(z)=az+b$, $a\ne0$, and $\widetilde P$ is the monic scalar
multiple of $P\circ\phi$, then
$\rho_{\widetilde P}=(\phi^{-1})_*\rho_P$ and
$\Sh_{\widetilde P}=\phi^{-1}(\Sh_P)$; the same covariance holds for
the full shadow \cite[Proposition~5.2]{HS1}.

\subsection{Definable sets and stratifications}\label{subsec:tame}
A subset of $\R^N$ is \emph{definable} here if it is given, with
real parameters, by a first-order formula in the real exponential field
$\Rexp=(\R;<,+,\cdot,\exp)$; a function is definable if its graph is.
This structure is \emph{o-minimal}:
every definable subset of the real line is a finite union of intervals
and points \cite{Wilkie}. A definable set $S$ admits finite $C^k$
\emph{stratifications}, compatible with any specified finite family of
definable subsets and with the ranks of specified definable maps. Such a
stratification is a partition of $S$ into connected $C^k$ submanifolds, its
\emph{strata}, which satisfies the \emph{frontier condition}: for each
stratum $C$, the set $(\overline C\setminus C)\cap S$ is a union of
lower-dimensional strata. The dimension of a definable set is the
largest dimension of a stratum. Projection and closure preserve
definability, and closure preserves dimension \cite{vdD}.

\begin{proposition}[Definable polynomial families]\label{prop:polynomial-tame}
\cite[Theorem~4.1]{HS2}
The barrier, the fixed-ratio supports, and the families that they form over
the monic polynomials of a fixed degree are definable. Their incidence sets
therefore admit finite stratifications compatible with prescribed definable
subsets and with the ranks of prescribed projections.
\end{proposition}
\begin{proof}
The definability statement of \cite{HS2} is joint in all parameters and
includes every root multiplicity pattern. There the ratio and the barrier
height are normalized by the degree; multiplying both by $d$ recovers the
ratio $\alpha$ and the barrier $\Blev_P$ for monic polynomials, as in
\eqref{eq:polynomial-weighted-bridge}. For the frontier condition, start with a
compatible smooth cell partition on which all prescribed map ranks are
constant. Put all cell closure differences and all lower-dimensional
cells into one closed definable exceptional set of smaller dimension.
The connected components of the remaining top-dimensional cells retain
the same chart maps. Recursively stratify the exceptional set, compatibly
also with the closures of these components. This gives the relative
frontier condition and retains all original constant ranks.
\end{proof}

A \emph{bounded definable curve}, that is, a bounded definable set of
dimension at most one, is a finite union
of points and $C^1$ arcs of finite length: subdividing into arcs on which
both coordinates are monotone bounds the length by the sum of their
coordinate variations. We refer to this as the \emph{monotone-arc argument}.
We shall also use semialgebraic triangulation,
quantifier elimination, and compatible topological trivializations
\cite{BCR,Hardt}. Semialgebraic triviality means that a finite partition
of parameter space into semialgebraic pieces makes the given family
topologically a product over each piece, simultaneously for the specified
subsets.

For an analytic branch of saddles $w=w(z,\alpha)$ with critical height
$\psi(z,\alpha)=\Ree\Phi_{P,z,\alpha}(w(z,\alpha))$, differentiation
at the critical point gives
\begin{equation}\label{eq:selected-spatial}
 2\partial_z\psi=\frac{\alpha}{w-z},\qquad
 \partial_\alpha\psi=-\log|w-z|.
\end{equation}
Distinct saddles thus have distinct spatial gradients.

\subsection{The compact incidence space}
\begin{figure}[t]
\centering
\includegraphics{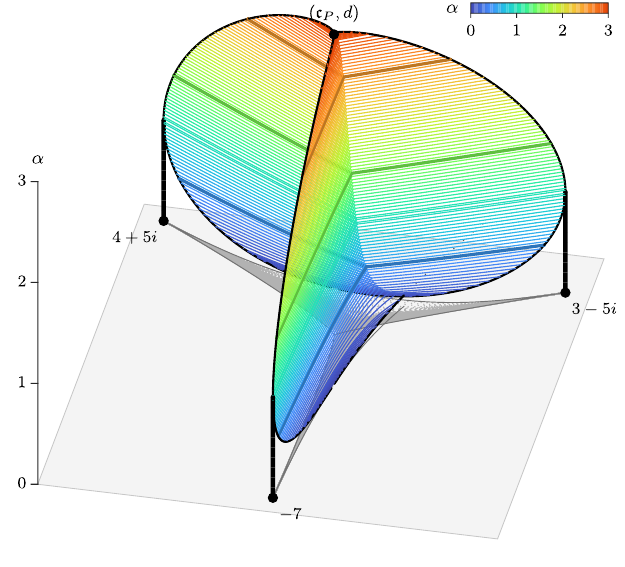}
\caption{The incidence space $\Inc_P$ for the cubic of
Figure~\ref{fig:regimes}(a): the supports $\FixSh_{\alpha,P}$, traced
numerically, at height $\alpha$ and colored by $\alpha$ (thicker:
$\alpha=\tfrac12,\tfrac32,\tfrac52$). Black: critical values of
$F_{P,\alpha}$ that are endpoints of $\FixSh_{\alpha,P}$, at height
$\alpha$; dotted: the other critical values; thick segments: the
segments $\{a_j\}\times[0,1]$ above the roots. The floor shows the zeros of
$(P^{110})^{(m)}$, approximating $\prz(\Inc_P)=\Sh_P$.}
\label{fig:incidence}
\end{figure}

For fixed $P$, define the \emph{support incidence space} by
\begin{equation}\label{eq:incidence}
 \Inc_P:=\overline{\{(z,\alpha):0<\alpha<d,\ z\in \FixSh_{\alpha,P}\}}
 \subset \Hull_P\times[0,d],
 \qquad \prz(z,\alpha):=z.
\end{equation}
The closure includes the \emph{endpoint fibers} over $\alpha=0$ and
$\alpha=d$, which are described in Proposition~\ref{prop:global-facts}.
The set $\Inc_P$ is compact and definable, and
$\prz(\Inc_P)=\Sh_P$. Away from roots and saddle collisions, the barrier is
locally one of finitely many harmonic critical-height branches. Its Riesz
measure, namely its distributional Laplacian divided by $2\pi$,
is therefore supported where at least two selected heights agree.
By \eqref{eq:selected-spatial}, distinct saddles have distinct spatial
gradients, so the sets on which two of their heights agree are locally
smooth arcs. Consequently the fiber of $\Inc_P$ over each ratio
$\alpha\in(0,d)$ has dimension at most one, and the dimension theory of
definable sets, together with fiberwise closure, gives
$\dim\Inc_P\le2$. Figure~\ref{fig:incidence} shows $\Inc_P$ for a cubic.

\begin{proposition}[Projection principle]\label{prop:projection}
Choose a finite $C^2$ stratification of $\Inc_P$ compatible with the
endpoint fibers and the vertical \emph{root lines}
$\{a_j\}\times[0,d]$, and refine it so that the restriction of $\prz$ to each
stratum has constant rank. Every
point of $\Inc_P$ that projects to $\partial\Sh_P$ lies in a stratum on which $\prz$ has
rank less than two. Hence the boundary is contained in the union of the images of the
rank-deficient strata.
\end{proposition}
\begin{proof}
If the restriction of $\prz$ to a stratum has rank two at a point, then it
maps a neighborhood of that point in the stratum onto a neighborhood of the
image point in the observation plane. This neighborhood lies in
$\prz(\Inc_P)=\Sh_P$, so the image point is interior.
\end{proof}

The rank-two strata account for planar contributions. Images of the
remaining strata have dimension at most one.

\begin{lemma}[Closedness of the support incidence away from the roots]\label{lem:incidence-closed}
The set
\[
 \{(z,\alpha):z\notin Z(P),\ 0<\alpha<d,
                  \ z\in\FixSh_{\alpha,P}\}
\]
is relatively closed in $(\C\setminus Z(P))\times(0,d)$.
\end{lemma}
\begin{proof}
At a point $(z,\alpha)$ of $(\C\setminus Z(P))\times(0,d)$ outside this set,
the fixed-ratio potential is harmonic on a root-free disk about $z$.
Proposition~\ref{prop:selection-facts}\textup{(iii)} gives one simple
connecting saddle, and by
Proposition~\ref{prop:selection-facts}\textup{(ii)} its critical height
continues to equal the barrier on a neighborhood of $(z,\alpha)$. This
height is harmonic in the spatial variable for every nearby ratio, so that
neighborhood contains no point of the set.
\end{proof}

\section{Conservation and moment reconstruction}\label{sec:chord}
An exact identity for an arbitrary polynomial shows that, outside the convex
hull of its roots, the pooled zero measure of its derivatives is close to a
chord measure. We then use the selected inverse to prove the conservation
law on every component of the complement of the core. The resulting
identities recover the root moments and constrain the closed curves that
avoid the core.

\subsection{An exact finite-degree identity}
\begin{proposition}[Telescoping logarithmic derivatives]\label{prop:finite-chord}
Let $G$ be a polynomial of degree $N\ge1$. Its normalized root measure and root Cauchy transform are
\[
 \rootmu_G:=\frac1N\divz G,
 \qquad
 \rootC_G:=C_{\rootmu_G}=\frac1N\frac{G'}G.
\]
For $0\le m<N$, define the logarithmic derivative of the $m$th derivative by
\[
 \mathfrak q_m:=\frac{G^{(m+1)}}{G^{(m)}},
\]
and adopt the terminal convention $\mathfrak q_N:=0$. Define the \emph{pooled measure}
of the zeros of the derivatives of $G$ and its Cauchy transform by
\[
 \nu_G^{\mathrm{pool}}:=
 \frac{2}{N(N+1)}\sum_{m=0}^{N-1}\divz G^{(m)},
 \qquad
 C_G^{\mathrm{pool}}:=C_{\nu_G^{\mathrm{pool}}}.
\]
Then the following identity between rational functions holds:
\begin{equation}\label{eq:finite-telescope}
 \bigl(C_G^{\mathrm{pool}}\bigr)'
 =-\frac{N}{N+1}\rootC_G^{\,2}
 -\frac1{N(N+1)}\sum_{m=0}^{N-1}
 \left(\frac{\mathfrak q_m'}{\mathfrak q_m}\right)^2.
\end{equation}
For $G=P^n$, the second term is $O(n^{-1})$ locally uniformly on every set sufficiently far
from $\Hull_P$.
\end{proposition}
\begin{proof}
Logarithmic differentiation of
$\mathfrak q_m=G^{(m+1)}/G^{(m)}$ gives the elementary identity
\[
 \mathfrak q_{m+1}=\mathfrak q_m+\frac{\mathfrak q_m'}{\mathfrak q_m}.
\]
It is an identity between rational functions and holds also when some
$G^{(m)}$ has multiple zeros. Squaring and subtracting gives
\[
 \mathfrak q_{m+1}^2-\mathfrak q_m^2
 =2\mathfrak q_m'+\left(\frac{\mathfrak q_m'}{\mathfrak q_m}\right)^2.
\]
Sum over $m=0,\dots,N-1$, use
$\mathfrak q_0=N\rootC_G$ and $\mathfrak q_N=0$, and divide by $N(N+1)$. Since
$C_G^{\mathrm{pool}}=2\sum_m\mathfrak q_m/(N(N+1))$, the result is \eqref{eq:finite-telescope}.

Suppose the roots of $G$ lie in $|z|\le R$. By the Gauss--Lucas theorem, so do the roots of
every derivative. For $|z|>3R+1$, the Cauchy transform of any probability measure supported
in that disk is uniformly separated from zero after multiplication by $z$, while its
derivative is bounded by a constant times $|z|^{-2}$. Thus
$|\mathfrak q_m'/\mathfrak q_m|$ is bounded independently of $m$ and $N$. There are $N$
summands and the prefactor is $O(N^{-2})$, which proves the final assertion.
\end{proof}

For $G=P^n$, one has
$\rootC_G=\rootC_P=P'/(dP)$ and
$\nu_G^{\mathrm{pool}}=\bar\nu_{P,n}$. Passing to the limit in
\eqref{eq:finite-telescope} already gives
$C_{\rho_P}'=-\rootC_P^2$ near infinity and an $O(n^{-1})$ error for each fixed moment of
$\bar\nu_{P,n}$.

\subsection{A uniform exterior estimate}
\label{subsec:universal-exterior}

The preceding identity gives a quantitative result for every polynomial.
For a polynomial $G$ of degree $N$, let $\sigma_G$ be its \emph{chord
measure}: the law of a point chosen uniformly on the segment joining two
independent random roots with law $\rootmu_G=N^{-1}\divz G$.
All transforms below are evaluated outside $\Hull_G=\conv Z(G)$.

\begin{theorem}[Uniform exterior chord approximation]
\label{thm:universal-exterior-chord}
For every nonconstant polynomial $G$ of degree $N$, every $z\notin \Hull_G$,
and $\delta=\dist(z,\Hull_G)$,
\begin{align}
 \left|\bigl(C_G^{\mathrm{pool}}-C_{\sigma_G}\bigr)'(z)\right|
        &\le\frac{2}{(N+1)\delta^2},
        \label{eq:exterior-chord-derivative}\\
 \left|C_G^{\mathrm{pool}}(z)-C_{\sigma_G}(z)\right|
        &\le\frac{2}{(N+1)\delta}.
        \label{eq:exterior-chord-transform}
\end{align}
If $\Gamma$ is a rectifiable Jordan curve surrounding $\Hull_G$, at
positive distance $\delta_\Gamma$ from it, and $f$ is holomorphic on a
neighborhood of $\Gamma$ and its interior, then
\begin{equation}\label{eq:exterior-holomorphic-test}
 \left|\int f\dd\nu_G^{\mathrm{pool}}-\int f\dd\sigma_G\right|
 \le\frac{\operatorname{length}(\Gamma)}
             {\pi(N+1)\delta_\Gamma}\max_\Gamma|f|.
\end{equation}
\end{theorem}
\begin{proof}
We first estimate each summand in
\eqref{eq:finite-telescope}. Choose a nearest point $b\in \Hull_G$ to $z$
and put $e=(z-b)/\delta$. Convexity gives
$\Ree(\overline e(z-a))\ge\delta$ for every $a\in \Hull_G$.
For any nonconstant polynomial $\widetilde G$ with roots $a_j\in \Hull_G$,
listed with multiplicity, write $q=\widetilde G'/\widetilde G$. Then
\[
 |q'(z)|\le\sum_j|z-a_j|^{-2},\qquad
 \Ree(e q(z))=
   \sum_j\frac{\Ree(\overline e(z-a_j))}{|z-a_j|^2}
      \ge\delta\sum_j|z-a_j|^{-2}.
\]
Consequently $q(z)\ne0$ and
\begin{equation}\label{eq:logder-exterior-universal}
                   |q'(z)/q(z)|\le1/\delta.
\end{equation}
By the Gauss--Lucas theorem, this applies to $\widetilde G=G^{(m)}$ for
every $m<N$, in which case $q=\mathfrak q_m$.
The chord integral for a fixed pair of roots, differentiated in $z$ and
averaged over the pairs, gives $C_{\sigma_G}'=-\rootC_G^2$. Therefore
\eqref{eq:finite-telescope} yields
\[
 \bigl(C_G^{\mathrm{pool}}-C_{\sigma_G}\bigr)'
 =\frac{\rootC_G^2}{N+1}
       -\frac1{N(N+1)}\sum_{m=0}^{N-1}
                       (\mathfrak q_m'/\mathfrak q_m)^2.
\]
Since $|\rootC_G|\le1/\delta$, the derivative estimate follows.

On the exterior normal ray $z+te$, $t\ge0$, the distance to $\Hull_G$
is at least $\delta+t$. Both probability Cauchy transforms tend to
zero at infinity, so integration of their derivative difference along
this ray gives
\[
 |C_G^{\mathrm{pool}}(z)-C_{\sigma_G}(z)|
 \le\frac2{N+1}\int_0^\infty(\delta+t)^{-2}\dd t
 =\frac2{(N+1)\delta}.
\]
Finally, the winding number $\iota$ of $\Gamma$ about the points of
$\Hull_G$ is constant and nonzero. Cauchy's integral formula and Fubini's
theorem express $\iota$ times the difference of the two integrals of $f$ as
$(2\pi i)^{-1}$ times the contour integral
of $f(C_G^{\mathrm{pool}}-C_{\sigma_G})$. Since $|\iota|\ge1$, the factor
$1/(2\pi)$ and \eqref{eq:exterior-chord-transform} give
\eqref{eq:exterior-holomorphic-test}.
\end{proof}

\begin{corollary}[Exterior universality for arbitrary bounded root arrays]
\label{cor:arbitrary-array-exterior}
Let $G_j$ be polynomials of degrees $N_j\to\infty$ whose roots all lie in
one compact convex set $K$, and suppose that
$N_j^{-1}\divz G_j\Longrightarrow\eta$, where $\eta$ is an arbitrary
probability measure supported in $K$. If $\sigma_\eta$ is the law of a
uniform point on the chord between two independent samples from $\eta$, then
\[
 C_{G_j}^{\mathrm{pool}}\longrightarrow C_{\sigma_\eta}
                 \quad\text{locally uniformly on }\C\setminus K.
\]
Every weak subsequential limit of $\nu_{G_j}^{\mathrm{pool}}$ has the
same logarithmic potential as $\sigma_\eta$ on $\C\setminus K$, and the
moments of the pooled measures converge:
\begin{equation}\label{eq:arbitrary-array-chord-moments}
 \lim_{j\to\infty}\int z^k\dd\nu_{G_j}^{\mathrm{pool}}
  =\frac1{k+1}\sum_{\ell=0}^k
          \left(\int z^\ell\dd\eta\right)
          \left(\int z^{k-\ell}\dd\eta\right)
                    \qquad(k\ge0).
\end{equation}
For $G_j=P^j$, \eqref{eq:exterior-chord-transform} gives the explicit
error $2/[(dj+1)\dist(z,\Hull_P)]$ throughout $\C\setminus\Hull_P$.
\end{corollary}
\begin{proof}
Products of weakly convergent probabilities on the compact set $K$
converge weakly; the continuous chord map then gives
$\sigma_{G_j}\Longrightarrow\sigma_\eta$. Their Cauchy transforms
converge locally uniformly off $K$, and the finite-degree estimate
proves the first assertion. The Gauss--Lucas theorem gives weak compactness
of the pooled measures. Every weak limit point has the stated Cauchy
transform outside $K$;
the difference of the logarithmic potentials has zero derivative along
each exterior normal ray and tends to zero at infinity, so it vanishes.
The expansion of the exterior Cauchy transform at infinity identifies
all complex moments. Compactness thus gives the same limit for every
subsequence of moments, and hence convergence of the whole sequence.
Finally $C_{\sigma_\eta}'=-C_\eta^2$; equating coefficients at infinity
gives the convolution in \eqref{eq:arbitrary-array-chord-moments}.
\end{proof}

\subsection{Conservation on complementary components}
\begin{theorem}[Cauchy-transform conservation law]\label{thm:Cprime}
For every $z\notin\Sh_P$,
\begin{equation}\label{eq:Cprime}
 C_{\rho_P}'(z)
 =-\frac1{d^2}\left(\frac{P'(z)}{P(z)}\right)^2
 =-\rootC_P(z)^2.
\end{equation}
More precisely, let $w_\alpha(z)$ be the selected inverse from
Proposition~\ref{prop:fixed-facts}. The \emph{selected saddle path} $\mathscr W_z$
is the path
\[
 \alpha\longmapsto w_\alpha(z),\qquad 0<\alpha<d,
\]
oriented by increasing $\alpha$. It tends to $z$ as $\alpha\downarrow0$ and to infinity as
$\alpha\uparrow d$, and
\begin{equation}\label{eq:path-primitive}
 C_{\rho_P}(z)=\frac1{d^2}\int_{\mathscr W_z}
 \left(\frac{P'(w)}{P(w)}\right)^2\dd w.
\end{equation}
\end{theorem}
\begin{proof}
Fix a disk compactly contained in $\C\setminus\Sh_P$. Every
$\mu_{\alpha,P}$ is supported away from that disk, so its Cauchy transform is holomorphic in
$z$ there and continuous in $\alpha$. Proposition~\ref{prop:fixed-facts} and the
implicit-function theorem make
$w_\alpha$ locally real-analytic in $\alpha$, because
$F_{P,\alpha}'(w_\alpha)\ne0$. The selected branch avoids the roots of $P$ and the poles of $h_P$.

As $\alpha\downarrow0$, the transforms $C_{\alpha,P}$ converge uniformly on the disk to
$\rootC_P$. Every zero of $\rootC_P$ belongs to $\Sh_P$ by
Proposition~\ref{prop:global-facts}; hence $w_\alpha(z)\to z$. As $\alpha\uparrow d$, the
same proposition gives
$C_{\alpha,P}(z)\to(z-\centroid_P)^{-1}$ uniformly on the disk, and therefore
$w_\alpha(z)\to\infty$.

Recall $\LogDer_P=P'/P$. Along the selected path, the saddle relation
$\alpha=(w-z)\LogDer_P(w)$ gives
\[
 \LogDer_P(w)\dd\alpha
 =\bigl(\LogDer_P(w)^2+(w-z)\LogDer_P(w)\LogDer_P'(w)\bigr)\dd w
 =\tfrac12\dd\bigl((w-z)\LogDer_P(w)^2\bigr)
  +\tfrac12\LogDer_P(w)^2\dd w.
\]
The boundary term vanishes at both endpoints: it contains the factor $w-z$ at the initial
endpoint and is $O(w^{-1})$ at infinity. The definition \eqref{eq:core-measure} of $\rho_P$
and \eqref{eq:selected-inverse} now yield
\[
 C_{\rho_P}(z)
 =\frac2{d^2}\int_0^d(d-\alpha)C_{\alpha,P}(z)\dd\alpha
 =\frac2{d^2}\int_0^d\LogDer_P(w_\alpha(z))\dd\alpha,
\]
and the preceding identity between differentials gives \eqref{eq:path-primitive}, with
$|(d-\alpha)C_{\alpha,P}(z)|\le(d-\alpha)/\dist(z,\Sh_P)$.

To prove \eqref{eq:Cprime}, set
\[
 \widehat C_{\alpha,P}(z)
 =(d-\alpha)C_{\alpha,P}(z)
 =\LogDer_P(w_\alpha(z)).
\]
Differentiating
$w_\alpha=z+\alpha/\widehat C_{\alpha,P}$ and
$\widehat C_{\alpha,P}=\LogDer_P(w_\alpha)$ gives
\[
 \widehat C_{\alpha,P}\,\partial_\alpha\widehat C_{\alpha,P}
 =\partial_z\widehat C_{\alpha,P},
 \qquad
 \partial_\alpha(\widehat C_{\alpha,P}^2)
 =2\partial_z\widehat C_{\alpha,P}.
\]
If the disk is at distance $\delta>0$ from $\Sh_P$, then
$|\partial_z\widehat C_{\alpha,P}|\le(d-\alpha)/\delta^2$. Differentiation may therefore be
interchanged with integration in $\alpha$. The endpoint limits are
$\widehat C_{\alpha,P}(z)\to\LogDer_P(z)$ as $\alpha\downarrow0$ and
$\widehat C_{\alpha,P}(z)\to0$ as $\alpha\uparrow d$. Consequently,
\[
 C_{\rho_P}'(z)
 =\frac2{d^2}\int_0^d\partial_z\widehat C_{\alpha,P}(z)\dd\alpha
 =\frac1{d^2}
  [\widehat C_{\alpha,P}(z)^2]_{\alpha=0}^{\alpha=d}
 =-\frac{\LogDer_P(z)^2}{d^2}.
\]
The unrescaled transforms $C_{\alpha,P}(z)$ extend continuously and without
zeros to $\alpha\in[0,d]$, so they are bounded away from zero. The first of
the two identities displayed above also gives $\partial_\alpha\widehat C_{\alpha,P}
=(\partial_z C_{\alpha,P})/C_{\alpha,P}$, which is bounded. In
\eqref{eq:path-primitive} the path integral means the improper parameter
integral of $\LogDer_P(w_\alpha)^2\partial_\alpha w_\alpha$. Differentiating
the explicit selected inverse identifies this integrand with
$\widehat C_{\alpha,P}-\alpha\partial_\alpha\widehat C_{\alpha,P}$;
thus the path integral is absolutely convergent.

The argument applies to every disk compactly contained in $\C\setminus\Sh_P$,
and hence to every component of this complement.
\end{proof}

Figure~\ref{fig:landscape}(b) shows selected saddle paths for several
points outside the core.

\subsection{Random chords and moments}
\begin{theorem}[The chord law]\label{thm:chord}
Let $\sigma_P=\Law((1-\mathsf U)X+\mathsf UY)$ be the chord measure from the introduction. Then
\[
 C_{\rho_P}=C_{\sigma_P},\qquad
 U_{\rho_P}=U_{\sigma_P}
 \qquad\text{on }\C\setminus \Hull_P.
\]
The moment identity \eqref{eq:moments-intro} holds for every $k\ge0$.
\end{theorem}
\begin{proof}
For $a,b\in Z(P)$ and $z\notin\Hull_P$, the elementary identity
\[
 \int_0^1\frac{\dd t}{(z-(1-t)a-tb)^2}
 =\frac{1}{(z-a)(z-b)}
\]
follows from the primitive $t/((z-a)(z-(1-t)a-tb))$, including $a=b$.
Differentiating the chord Cauchy integral and averaging over the independent
random roots gives $C_{\sigma_P}'=-\rootC_P^2$. By Theorem~\ref{thm:Cprime},
$C_{\rho_P}$ has the same derivative on $\C\setminus\Hull_P$, and both
transforms tend to zero at infinity. On the outward normal ray from any
point outside the compact convex hull, their difference has zero derivative
and zero limit, so it vanishes at the initial point. Their logarithmic
potentials also differ by a constant along each such ray; the common
probability-mass normalization at infinity fixes this constant to zero.

For the moments, expand the random chord directly:
\[
 \mathbb E[((1-\mathsf U)X+\mathsf UY)^k]
 =\sum_{j=0}^k\binom kj\mom_j(P)\mom_{k-j}(P)
 \int_0^1(1-t)^jt^{k-j}\dd t
 =\frac1{k+1}\sum_{j=0}^k\mom_j(P)\mom_{k-j}(P).
\]
Equivalently, comparing Taylor coefficients at infinity in
$C_{\sigma_P}'=-\rootC_P^2$ gives the same convolution. Equality of the
exterior transforms identifies these quantities with the moments of $\rho_P$.
\end{proof}

\begin{figure}[t]
\centering
\includegraphics{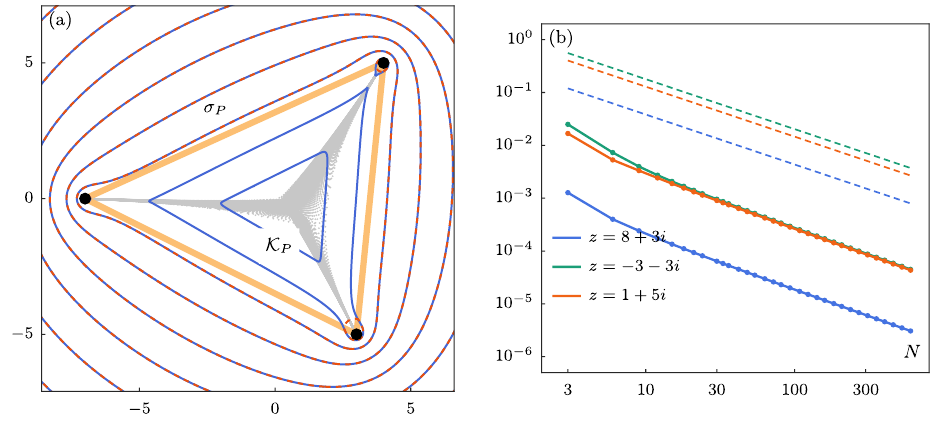}
\caption{(a)~Level curves of the logarithmic potentials of
$\bar\nu_{P,n}$ (solid; $n=110$) and of $\sigma_P$ (dashed) for the cubic
of Figure~\ref{fig:regimes}(a). They agree outside $\Hull_P$, where
$U_{\rho_P}=U_{\sigma_P}$, but not inside. Orange: the edges carrying
$\sigma_P$; gray: zeros.
(b)~$|C_G^{\mathrm{pool}}(z)-C_{\sigma_G}(z)|$ for $G=P^n$, $N=3n$, at
three points $z\notin\Hull_P$, with the bounds $2/((N+1)\delta)$ of
\eqref{eq:exterior-chord-transform} (dashed).}
\label{fig:chord}
\end{figure}

Figure~\ref{fig:chord}(a) compares the potentials of $\rho_P$ and
$\sigma_P$ for a cubic, and Figure~\ref{fig:chord}(b) illustrates the
rate in Theorem~\ref{thm:universal-exterior-chord}.

\subsection{What the shadow measure remembers}
\label{rev:inverse-measure}
Recall that a monic polynomial is \emph{primitive} if the greatest common
divisor of its root multiplicities is one.

\begin{theorem}[The fibers of the shadow measure]
\label{rev:measure-fibers}
For nonconstant monic polynomials $P,Q$, the following are equivalent:
\begin{enumerate}[label=\textup{(\roman*)}]
\item $\rho_P=\rho_Q$;
\item their normalized root measures are equal;
\item $P={P_*}^k$ and $Q={P_*}^\ell$ for a common primitive monic polynomial
$P_*$ and positive integers $k,\ell$.
\end{enumerate}
In particular, the monic polynomials with a given shadow measure form
exactly the fiber
$\{P_*,{P_*}^2,{P_*}^3,\ldots\}$ for a unique primitive $P_*$. If the degree $d$ is
known, the first $d$ holomorphic moments of $\rho_P$ determine $P$.

More explicitly, put $\mathfrak s_j=\int z^j\,\dd\rho_P(z)$, $\mom_0=1$, and
recursively define
\begin{equation}\label{rev:inverse-moments}
 \mom_j=\frac12\left((j+1)\mathfrak s_j-
            \sum_{i=1}^{j-1}\mom_i \mom_{j-i}\right),\qquad j\ge1.
\end{equation}
Then $\mom_j$ is the $j$th moment of the normalized root measure. Near
infinity the normalized root Cauchy transform is recovered by
\begin{equation}\label{rev:inverse-square-root}
 \rootC_P(z)=\sqrt{-C_{\rho_P}'(z)},
 \qquad \rootC_P(z)=z^{-1}+O(z^{-2}),
\end{equation}
where the asymptotic condition specifies the square-root branch.
\end{theorem}
\begin{proof}
The atomic formula \eqref{eq:shadow-atoms} gives
\begin{equation}\label{eq:atomic-reconstruction}
 \rootmu_P=\sum_{a:\,\rho_P(\{a\})>0}
                 \sqrt{\rho_P(\{a\})}\,\delta_a.
\end{equation}
Thus \textup{(i)} implies \textup{(ii)} directly. The exterior
reconstruction in \eqref{rev:inverse-square-root} uses less information:
Theorem~\ref{thm:Cprime} gives the square identity. The Cauchy transform
of the compactly supported probability measure $\rootmu_P$ has the
convergent expansion
$\rootC_P(z)=z^{-1}+O(z^{-2})$ at infinity. Any two square roots with this asymptotic
have nonzero sum near infinity, so factoring the difference of their
squares proves equality there. The rational continuation has the root
locations as poles and the normalized multiplicities as residues. Finally,
\eqref{rev:inverse-moments} follows by isolating the two terms containing
$\mom_j$ in the moment identity of Theorem~\ref{thm:chord}.

Write the distinct roots of $P$ as $a_j$, their multiplicities
as $r_j$, and $g=\gcd_j r_j$. Put $r_j^*=r_j/g$ and
\[
 {P_*}(z)=\prod_{j=1}^s(z-a_j)^{r_j^*}.
\]
Then $P={P_*}^g$ and the multiplicities of $P_*$ have gcd one. If $Q$
has the same normalized root measure, it has the same distinct roots.
Writing its multiplicities as $\tilde r_j$, its degree as $\tilde d$, and
$\tilde g=\gcd_j\tilde r_j$,
equality of root masses gives $\tilde d r_j=d\tilde r_j$. Taking gcds gives
$\tilde d g=d\tilde g$;
cancellation therefore gives $r_j/g=\tilde r_j/\tilde g$ for every $j$. Hence
$Q={P_*}^{\tilde g}$. This also proves uniqueness of the primitive representative,
and establishes \textup{(ii)} implies \textup{(iii)}.
Conversely, the pooled measures \eqref{eq:pooled-powers} of ${P_*}^k$ at
power level $n$ and of $P_*$ at power level $kn$ agree exactly, including their
normalizations and all their derivative orders, because
$({P_*}^k)^n={P_*}^{kn}$. Letting $n\to\infty$ proves
$\rho_{P_*^k}=\rho_{P_*}$, so \textup{(iii)} implies \textup{(i)}.

For the finite reconstruction at known degree, note that the power sums of
the roots of $P$, counted with multiplicity, are $d\mom_j$ for
$1\le j\le d$. Newton's identities recover the elementary symmetric
functions $e_j$ of the roots by $e_0=1$ and
\[
 j e_j=d\sum_{i=1}^j(-1)^{i-1}e_{j-i}\mom_i,\qquad 1\le j\le d.
\]
Then $P(z)=\sum_{j=0}^d(-1)^j e_jz^{d-j}$. The one-root case is
included: its shadow measure is the point mass at that root.
\end{proof}

\begin{corollary}[The discrete ambiguity in recovering the full shadow]
\label{rev:measure-versus-full}
If $P_*$ is primitive, all polynomials in the fiber of $\rho_{P_*}$, that
is, all powers ${P_*}^k$, have the same core, while
\begin{equation}\label{rev:power-fiber-shadows}
 \AllSh_{P_*^k}\subseteq\AllSh_{P_*},\qquad
 \bigcap_{k\ge1}\AllSh_{P_*^k}=\Sh_{P_*},\qquad
 \AllSh_{P_*^k}\longrightarrow\Sh_{P_*}\quad\text{in Hausdorff distance}.
\end{equation}
Thus the shadow measure determines both the common core and the largest
full shadow in its fiber. It need not determine
which of the full shadows $\AllSh_{P_*^k}$ is present when the exponent $k$
is unspecified.
\end{corollary}
\begin{proof}
Each descendant $(({P_*}^k)^n)^{(m)}=({P_*}^{kn})^{(m)}$ of ${P_*}^k$ is a
descendant of $P_*$, proving the
inclusion. The intersection and convergence follow from the removal of transients
under passage to powers \cite[Proposition~5.10]{HS1}.
The measure determines $P_*$ by Theorem~\ref{rev:measure-fibers}, hence
also determines the two sets so specified. Different members of the fiber
can indeed have different full shadows. The polynomial
\[
 Q_*(z)=z^4+\tfrac12z^3-\tfrac{15}{32}z^2+\tfrac12z+1
\]
has a transient by \cite[Theorem~6.4]{HS1}. Let $P_*$ be the primitive
polynomial of which $Q_*$ is a power, so that $Q_*={P_*}^r$ for some
positive integer $r$.
Choose a transient $x$ of $Q_*$. Power invariance of the core and the inclusion
$\AllSh_{P_*^r}\subseteq\AllSh_{P_*}$ make $x$ a transient of $P_*$.
By \cite[Proposition~5.10]{HS1}, $x$ is absent from
$\AllSh_{P_*^k}$ for all sufficiently large $k$, although
$\rho_{P_*^k}=\rho_{P_*}$.
\end{proof}

\subsection{Periods and winding transport}
In the next proposition, $D$ is a bounded open set whose boundary consists
of finitely many disjoint regular simple $C^1$ loops, and boundary integrals
are sums over these loops. Each loop $\gamma$ is positively oriented, that
is, it has $D$ on its left: in the local chart
$w\mapsto\gamma(t+\Ree w)+i\Imm w\,\gamma'(t)$, the set $D$ occupies $\Imm w>0$.
The winding numbers of the loops consequently add up to $1$ on $D$ and to
$0$ outside its closure.
\begin{proposition}[Moments inside a boundary disjoint from the core]\label{prop:component-moments}
Let $D$ be a bounded open set with smooth positively oriented boundary and
$\partial D\cap\Sh_P=\varnothing$. Then
\begin{align}
 \oint_{\partial D}\rootC_P(z)^2\dd z&=0,\label{eq:res-zero}\\
 \int_Dz^k\dd\rho_P(z)&=
 \frac{1}{2\pi i(k+1)}
 \oint_{\partial D}z^{k+1}\rootC_P(z)^2\dd z,
 \qquad k\ge0.\label{eq:component-moments}
\end{align}
In particular, $\rho_P(D)=0$ if $D$ contains no root, while
$\rho_P(D)=1$ if $D$ contains every root.
\end{proposition}
\begin{proof}
Integrating $C_{\rho_P}'=-\rootC_P^2$ around every boundary component gives
\eqref{eq:res-zero}. Fubini's theorem in the Cauchy integral, followed by integration by
parts, gives
\[
 2\pi i\int_Dz^k\dd\rho_P(z)
 =\oint_{\partial D}z^kC_{\rho_P}(z)\dd z
 =\frac1{k+1}\oint_{\partial D}z^{k+1}\rootC_P(z)^2\dd z.
\]
The core is connected by \cite[Theorem~5.6]{HS1} and contains all roots.
Since it avoids $\partial D$, it lies wholly in $D$ or wholly outside $D$.
A root fixes the side, proving both mass assertions.
\end{proof}

\begin{lemma}[Winding-number transport]\label{lem:winding-transport}
Let $\gamma$ be a closed piecewise smooth curve in $\C\setminus\Sh_P$. For each root $a_j$,
\begin{equation}\label{eq:winding-transport}
 \wind(\gamma,a_j)=\wind(\gamma,\centroid_P).
\end{equation}
\end{lemma}
\begin{proof}
By the endpoint limits in Proposition~\ref{prop:global-facts}, the measures
$\mu_{\alpha,P}$, $0<\alpha<d$, extend weakly continuously to $[0,d]$ by
$\mu_{0,P}:=\rootmu_P$ and $\mu_{d,P}:=\delta_{\centroid_P}$.
All these measures are carried by $\Sh_P$, so their Cauchy transforms
$C_{\alpha,P}$ are jointly
continuous on $[0,d]\times(\C\setminus\Sh_P)$ and are nonzero there.
At the lower endpoint this uses the inclusion of all roots and critical
points in the core; at the upper endpoint the transform is
$(z-\centroid_P)^{-1}$.

For a root $a_j$, define the rescaled displacement
\[
 H_\alpha(z)=\frac{d-\alpha}{d}(z-a_j)
                 +\frac{\alpha}{d\,C_{\alpha,P}(z)}.
\]
For $0<\alpha<d$ it equals $(d-\alpha)(w_\alpha(z)-a_j)/d$ and is nonzero
by Proposition~\ref{prop:fixed-facts}. Its endpoint values are
$H_0(z)=z-a_j$ and $H_d(z)=z-\centroid_P$, both nonzero outside the core.
Thus $H_\alpha\circ\gamma$ is a continuous homotopy of nonvanishing
closed curves between these two displacement curves. Their winding
numbers about zero are equal, proving \eqref{eq:winding-transport}.
\end{proof}

The continuity of the supports also gives paths to the centroid. Precisely, for
each $z\in\FixSh_{\alpha_0,P}$ there is a continuous definable path
$\gamma:[\alpha_0,d]\to\Sh_P$ with
\[
 \gamma(\alpha_0)=z,\qquad \gamma(d)=\centroid_P,\qquad
 \gamma(\alpha)\in\FixSh_{\alpha,P}\quad(\alpha_0\le\alpha<d)
\]
\cite[Lemma~4.7 and Corollary~4.8]{HS2}. The existence of these paths is compatible with the presence
of \emph{holes}, that is, bounded components of the complement of the core.
Lemma~\ref{lem:winding-transport} supplies an additional
constraint on each loop that avoids the core: the residue condition
\eqref{eq:res-zero} alone can
hold for a proper symmetric subset of the roots.

\section{Switching sheets and the shadow measure}
\label{sec:balance}

The nonatomic part of a fixed-ratio measure is carried by the interfaces
across which the selected saddle changes; an equality between two heights
that are not selected contributes no mass. We compute the measure on these
interfaces
and then project it to the observation plane. The projection separates
the planar and curve parts of the shadow measure and identifies the
possible regular boundary points.

\subsection{Selected switching charts}
\begin{definition}[Regular switching chart]\label{def:switch}
A \emph{regular switching chart} consists of an open set $V\subset\C\times(0,d)$ and two
distinct functions $w_1,w_2$ on $V$, holomorphic in $z$ and real-analytic in $\alpha$, with
the following two properties. First,
\[
 F_{P,\alpha}(w_j(z,\alpha))=z,\qquad
 F_{P,\alpha}'(w_j(z,\alpha))\ne0,\qquad
 P(w_j)P'(w_j)(w_j-z)\ne0,
\]
where $F_{P,\alpha}'$ denotes the derivative with respect to $w$; thus $w_1$ and $w_2$ are
distinct branches of simple ordinary saddles. Second, with compatible local branches of the
logarithms and
\[
 \Psi_j(z,\alpha):=\Phi_{P,z,\alpha}(w_j(z,\alpha)),
 \qquad
 \psi_j(z,\alpha):=\Ree\Psi_j(z,\alpha),
\]
the chart is \emph{switching}, that is,
\[
 \Blev_P(z,\alpha)=\max\{\psi_1(z,\alpha),\psi_2(z,\alpha)\}
 \quad\text{on }V,
\]
and both strict inequalities $\psi_1>\psi_2$ and $\psi_2>\psi_1$ occur in $V$. Define
\[
 g_{12}:=\psi_1-\psi_2,
 \qquad
 \Sigma_{12}:=\{g_{12}=0\}\cap V.
\]
The set $\Sigma_{12}$ is the local \emph{switching sheet}, and $g_{12}=0$ is the
\emph{switching equation}. Each branch must be selected on
its indicated side; equality of two unselected heights does not define a switching chart.
\end{definition}

The maximum is a local representation of the barrier, involving only the
two branches that exchange selection. By continuity, their heights are equal
on the interface; that the barrier is their maximum, rather than their
minimum, is forced by the positivity of its Riesz measure.

\begin{proposition}[Derivatives of a switching pair]\label{prop:dictionary}
On a regular switching chart,
\begin{align}
 2\partial_z g_{12}
 &=\alpha\left(\frac1{w_1-z}-\frac1{w_2-z}\right),\label{eq:g-z}\\
 \partial_\alpha g_{12}
 &=\log\frac{|w_2-z|}{|w_1-z|}
 =\log\frac{|h_P(w_2)|}{|h_P(w_1)|}.\label{eq:g-alpha}
\end{align}
In particular, $\nabla_zg_{12}\ne0$, where $\nabla_z$ is the real gradient in the two
real coordinates of $z$. Using the same local branches of the logarithms, define the
auxiliary phase
\[
 \Theta_\alpha(w):=(1-\alpha)\log P(w)+\alpha\log P'(w).
\]
At an ordinary saddle,
\begin{equation}\label{eq:theta}
 \Psi_j=\Theta_\alpha(w_j)-\alpha\log\alpha,
 \qquad
 \frac{\dd}{\dd\alpha}\Ree\Theta_\alpha(w_j(z,\alpha))
 =1-\log|h_P(w_j)|.
\end{equation}
The complex identity uses compatible logarithm branches; the real identities do not depend
on those choices. Moreover,
\begin{equation}\label{eq:g-second}
 \partial_\alpha^2g_{12}=
 \Ree\left(
 \frac{h_P'(w_2)}{1-\alpha h_P'(w_2)}-
 \frac{h_P'(w_1)}{1-\alpha h_P'(w_1)}
 \right).
\end{equation}
\end{proposition}
\begin{proof}
At a saddle the derivative of $\Phi_{P,z,\alpha}$ with respect to $w$ vanishes.
Differentiation in the two base variables therefore gives
\[
 \partial_z\Psi_j=\frac\alpha{w_j-z},\qquad
 \partial_\alpha\Psi_j=-\log(w_j-z).
\]
Taking real parts and subtracting proves \eqref{eq:g-z}--\eqref{eq:g-alpha}; the
distinctness of $w_1$ and $w_2$ gives $\nabla_zg_{12}\ne0$. Choose local branches of
$\log P$ and $\log P'$ and put
$\log(w_j-z)=\log\alpha+\log P(w_j)-\log P'(w_j)$. By the saddle equation
$w_j-z=\alpha h_P(w_j)$, this is a compatible branch of the logarithm of the
displacement $w_j-z$, and substitution gives
\eqref{eq:theta}. Finally, differentiating
$z=w_j-\alpha h_P(w_j)$ at fixed $z$ yields
\[
 \partial_\alpha w_j=\frac{h_P(w_j)}{1-\alpha h_P'(w_j)}.
\]
A second derivative of \eqref{eq:g-alpha} now gives \eqref{eq:g-second}.
\end{proof}

For a fixed ratio in the chart, write
\[
 V_\alpha:=\{z:(z,\alpha)\in V\},\qquad
 \mathscr T_{12,\alpha}:=\{z\in V_\alpha:g_{12}(z,\alpha)=0\}.
\]
By Proposition~\ref{prop:dictionary}, every nonempty
$\mathscr T_{12,\alpha}$ is a smooth real-analytic curve. We call it the local \emph{tie curve}:
along it the two critical heights are \emph{tied}, that is, equal.

A point of the switching sheet is \emph{stationary in the ratio} when $\partial_\alpha
g_{12}=0$. Two distinct ordinary saddles with equal values of both $|P|$ and $|P'|$ form a
\emph{balanced pair}. A \emph{balance fold} is a point of the switching sheet with
$\partial_\alpha g_{12}=0$ and $\partial_\alpha^2g_{12}\ne0$. Its projection is
\emph{exposed} if it lies locally on $\partial\Sh_P$; it is an \emph{internal density ridge}
if other selected strata fill both sides. The rank of the projection is
always the rank of its real differential
restricted to the relevant stratum. The \emph{ordinary saddle-collision
locus} consists of the parameter values at which two inverse branches of
$F_{P,\alpha}$ coalesce at a point where $P(w)P'(w)\ne0$.

\subsection{Projection of the zero measure}
Hausdorff measures use the Euclidean metric and are normalized to agree
with Lebesgue measure in the matching dimension.
We use the following change-of-variables convention. If $\Pi:M\to\R^\ell$
is a $C^1$ map from a smooth $k$-dimensional submanifold $M$, with
$k\ge\ell$, its \emph{$\ell$-dimensional Jacobian} $J_\ell\Pi$ is the
product of the $\ell$ singular values of its differential, including
zeros when the rank is less than $\ell$. The \emph{coarea formula}
reads
\[
 \int_M u\,J_\ell\Pi\dd\Haus^k
 =\int_{\R^\ell}\left(\int_{\Pi^{-1}(y)}u\dd\Haus^{k-\ell}\right)\dd y
\]
for nonnegative compactly supported Borel functions $u$ on $M$.
When $k=\ell$, the inner integral is a sum over preimages; this is
the \emph{area formula}. These formulas, used below only on smooth
charts, convert densities on a surface into densities on its image
or along its fibers; see \cite{Federer}.

Fix a regular switching chart and a nonnegative continuous cutoff $\chi$ compactly supported
in it. Define the corresponding weighted contribution $\rho_\chi$ by
\[
 \int \varphi(z)\dd\rho_\chi(z)
 :=\frac2{d^2}\int_0^d(d-\alpha)
       \int \varphi(z)\chi(z,\alpha)\dd\mu_{\alpha,P}(z)\dd\alpha
\]
for bounded continuous test functions $\varphi$. When several charts are used, a partition
of unity prevents double counting.

\begin{theorem}[Local density formula]\label{thm:density}
On a regular switching chart, the nonatomic part of the fixed-ratio measure carried by the chart is
\begin{equation}\label{eq:line-density}
 \mu_{\alpha,P}^{\mathrm{na}}\!\restriction V_\alpha
 =\frac{|\nabla_zg_{12}|}{2\pi(d-\alpha)}\,
 \Haus^1\!\restriction\mathscr T_{12,\alpha}.
\end{equation}
Suppose that every contributing preimage of $z$ in the chart satisfies $\partial_\alpha
g_{12}\ne0$; equivalently, the differential of the projection is surjective at every such
preimage. Then the planar density of the weighted average is
\begin{equation}\label{eq:density}
 \frac{\dd\rho_\chi}{\dd A}(z)=
 \frac1{\pi d^2}
 \sum_{\substack{\alpha:(z,\alpha)\in V,\ g_{12}(z,\alpha)=0}}
 \frac{\chi(z,\alpha)|\nabla_zg_{12}(z,\alpha)|^2}
 {|\partial_\alpha g_{12}(z,\alpha)|}.
\end{equation}
The sum is over solutions satisfying $\partial_\alpha g_{12}\ne0$. Thus
\eqref{eq:density} gives the rank-two contribution and makes no
assertion about a stationary rank-one sheet.
\end{theorem}
\begin{proof}
The distributional Laplacian of the maximum of two harmonic functions has, along a smooth
interface where their gradient difference is nonzero, arclength density equal to the
magnitude of that difference. Applied to the fixed-ratio potential, this jump formula gives
\eqref{eq:line-density}. Equivalently, integrate by parts on the two sides of the tie curve;
the outward normal derivatives differ by $|\nabla_zg_{12}|$.

The averaging factor $2(d-\alpha)/d^2$ cancels the denominator $d-\alpha$ in
\eqref{eq:line-density}. Before projection, the resulting measure on $\Sigma_{12}$ is
\begin{equation}\label{eq:incidence-density}
 \frac1{\pi d^2}
 \frac{\chi|\nabla_zg_{12}|^2}{|\nabla_{z,\alpha}g_{12}|}
 \Haus^2\!\restriction\Sigma_{12}.
\end{equation}
Here $\nabla_{z,\alpha}g_{12}$ is the full real gradient on
$\C\times\R\simeq\R^3$. The two-dimensional Jacobian of
$\prz(z,\alpha)=z$, restricted to $\Sigma_{12}$, is
$|\partial_\alpha g_{12}|/|\nabla_{z,\alpha}g_{12}|$. The area formula therefore gives
\eqref{eq:density}. For fixed $z$, the function $\alpha\mapsto g_{12}(z,\alpha)$ is
real-analytic, so its zeros with $\partial_\alpha g_{12}\ne0$ are locally finite; on a
component where it vanishes identically, $\partial_\alpha g_{12}=0$ and there are no such zeros.
Compact support inside the chart therefore leaves only finitely many nonzero terms.
\end{proof}

\begin{theorem}[Balance-fold singularity]\label{thm:fold}
Let $(z_*,\alpha_*)$ be a point of a regular switching chart satisfying
\[
 g_{12}=\partial_\alpha g_{12}=0,
 \qquad
 \kappa_*:=\partial_\alpha^2g_{12}(z_*,\alpha_*)\ne0.
\]
Then the projection of the switching sheet has a fold at $(z_*,\alpha_*)$: locally, points
on one side of the fold image have two nearby preimages on the sheet, while points on the
other side have none. Set
\[
 \mathfrak n_*:=|\nabla_zg_{12}(z_*,\alpha_*)|,
 \qquad
 \mathbf n:=\frac{\nabla_zg_{12}(z_*,\alpha_*)}{\mathfrak n_*},
\]
and let $\tau_\perp$ be signed Euclidean distance from $z_*$ in the direction $\mathbf n$. The
side with two preimages is the \emph{occupied side}. For a cutoff equal to one near the
fold, write $D_\chi^{\mathrm{fold}}$ for the sum in \eqref{eq:density} over these two nearby
preimages. Using this canonical representative, along the normal line on the occupied side,
\begin{equation}\label{eq:fold-density}
 D_\chi^{\mathrm{fold}}(z_*+\tau_\perp\mathbf n)
 =\frac{\sqrt2\,\mathfrak n_*^{3/2}}{\pi d^2\sqrt{|\kappa_*|}}
 |\tau_\perp|^{-1/2}\bigl(1+o(1)\bigr).
\end{equation}
On the other side this local contribution vanishes. The fold image is an exposed boundary
only if no other selected stratum fills the vacant side; otherwise it is an internal density
ridge.
\end{theorem}
\begin{proof}
Because $\nabla_zg_{12}\ne0$, choose local real coordinates in the observation plane and
solve the switching equation as $y=\mathfrak y(x,\alpha)$. At the base point,
$\partial_\alpha\mathfrak y=0$ and $\partial_\alpha^2\mathfrak y\ne0$. The
parameter-dependent Morse lemma gives a fold with two preimages on one side and none on the
other. Put $\widehat\alpha:=\alpha-\alpha_*$. Along the Euclidean normal at $z_*$, Taylor
expansion gives
\[
 g_{12}=\mathfrak n_*\tau_\perp+\tfrac12\kappa_*\widehat\alpha^2
 +O(\tau_\perp^2+|\tau_\perp\widehat\alpha|+|\widehat\alpha|^3).
\]
On the occupied side the two solutions satisfy
\[
 |\widehat\alpha|
 =\left(\frac{2\mathfrak n_*|\tau_\perp|}{|\kappa_*|}\right)^{1/2}(1+o(1)),
 \qquad
 |\partial_\alpha g_{12}|
 =\bigl(2\mathfrak n_*|\kappa_*||\tau_\perp|\bigr)^{1/2}(1+o(1)).
\]
Their two contributions in \eqref{eq:density}, together with
$|\nabla_zg_{12}|=\mathfrak n_*+o(1)$, give \eqref{eq:fold-density}.
Contributions of other parts of the incidence space are nonnegative. They
cannot cancel the
singularity, although they may fill the opposite side and make the
fold internal. The critical set of the projection on the sheet, where
$\partial_\alpha g_{12}=0$, is a graph $\alpha=c(x)$; it is a curve and
therefore has zero area in the parameters of the sheet. By Fubini's theorem
and the fixed-ratio law \eqref{eq:line-density}, it carries no mass of the
averaged measure. The density $D_\chi^{\mathrm{fold}}$ is therefore the
density of the finite measure $\rho_\chi$, hence integrable. The fold image,
a real-analytic curve, has zero area and carries no mass.
\end{proof}

\begin{remark}[Odd stationary order]
If instead
$\partial_\alpha g_{12}=\partial_\alpha^2g_{12}=0$ and
$\partial_\alpha^3g_{12}\ne0$, the image of the switching sheet contains a full neighborhood
of $z_*$. Indeed, in the representation $y=\mathfrak y(x,\alpha)$, the values at
$\alpha_*\pm\eps$ lie on opposite sides of $\mathfrak y(x_*,\alpha_*)$ for every
sufficiently small fixed $\eps$. This remains true for nearby $x$, and the
intermediate-value theorem fills a common intervening interval. Thus an odd-order stationary
point cannot by itself create an exposed boundary.
\end{remark}

\subsection{The global measure decomposition}
\begin{theorem}[Planar, curve, and atomic parts]\label{thm:decomposition}
The decomposition \eqref{eq:decomposition-intro} holds. Away from the root atoms, the
curve part is generated by two-dimensional selected incidence strata whose projection to the
$z$-plane has rank one. Such strata are locally stationary families of support arcs. By
contrast, an isolated fold curve inside a rank-two switching sheet creates an integrable
singularity of the planar density, not an additional arclength measure.
\end{theorem}
\begin{proof}
Consider the positive \emph{incidence measure}
\[
 \dd\Omega_P(z,\alpha)=\frac2{d^2}(d-\alpha)
       \dd\mu_{\alpha,P}(z)\dd\alpha,\qquad 0<\alpha<d.
\]
Its spatial projection is $\rho_P$. Its restriction to the root lines gives
exactly the atoms in
\eqref{eq:shadow-atoms}; write $\Omega_P^\circ$ for the remaining measure.

First, $\mu_{\alpha,P}\restriction(\C\setminus Z(P))\ll\Haus^1$.
Indeed, the first-variation formula, valid also at multiple connecting saddles
\cite[Theorem~A.1]{HS1}, bounds the one-sided spatial derivatives of the barrier
by the spatial differentials of the landscape at the connecting saddles
(compare \eqref{eq:selected-spatial}). Nearby connecting saddles stay in a compact region
away from the roots and the observation point, so these bounds are locally
uniform. The one-sided mean-value inequality and \eqref{eq:barrier-potential}
give local Lipschitz continuity of $U_{\mu_{\alpha,P}}$ off $Z(P)$, including
at collisions. If this potential is $L$-Lipschitz on a disk of radius $2r$,
the circle-mean formula gives
\[
 \mu_{\alpha,P}(\overline B(z,r))\log2
 \le U^{(2r)}(z)-U^{(r)}(z)\le3Lr,
\]
where $U^{(r)}$ is the potential with the kernel $\log|z-\xi|$ replaced by
$\log\max\{r,|z-\xi|\}$. Covering null sets by small disks proves the asserted
absolute continuity with respect to arclength. By Fubini's theorem, the spatial
marginal of $\Omega_P^\circ$ is therefore also absolutely continuous with
respect to $\Haus^1$.

The incidence space $\Inc_P$ is bounded and definable, of dimension at
most two. Apply Proposition~\ref{prop:polynomial-tame}, retaining constant ranks
of both the ratio and spatial projections. Strata of dimension at most one
have finite fibers outside finitely many ratios and hence zero
$\Omega_P^\circ$-mass. A connected stratum on which the differential of the
ratio vanishes has constant ratio and likewise contributes nothing.

On each remaining two-dimensional stratum, the inverse-function theorem
gives local coordinates $(x,\alpha)$ consisting of the ratio and one
spatial coordinate $x=\ell(z)$, where $\ell$ is a real linear coordinate
function on the observation plane. In these coordinates the stratum is the
graph of a spatial map $\gamma(x,\alpha)$ satisfying
$\ell(\gamma(x,\alpha))=x$. Its slices at fixed ratio are regular injective
arcs. Their arclength law and Fubini's theorem give
$\Omega_P^\circ\ll\dd x\dd\alpha$ in these
coordinates. If the spatial rank is two, the spatial map has smooth local
inverses, so its projected contribution is absolutely continuous with
respect to area. If the rank is one, $\partial_\alpha\gamma$ is a multiple
of $\partial_x\gamma$. Applying $\ell$ shows that this multiple is zero.
Thus $\gamma$ is independent of $\alpha$ on a smaller product neighborhood:
the original support arc is stationary.

The spatial images of the rank-one strata are bounded definable curves. By
the monotone-arc argument of Section~\ref{subsec:tame}, they decompose
into finitely many regular $C^1$ arcs of finite length and finitely many
points; the points
have zero nonatomic mass. Let $\Gamma_P$ be the union of these arcs. It has
zero area and finite $\Haus^1$ measure. Countably many local charts
suffice to patch the rank-two area laws. Consequently the spatial marginal
of $\Omega_P^\circ$ is absolutely continuous with respect to area off
$\Gamma_P$, and its restriction to $\Gamma_P$ is absolutely continuous with
respect to arclength. The
Radon--Nikodym theorem supplies the nonnegative integrable densities in
\eqref{eq:decomposition-intro}. Finally, Theorem~\ref{thm:fold} proves that
the critical set of an isolated fold has zero mass and that its projected
density is integrable; it adds no curve measure. Restoring the root atoms
completes the decomposition.
\end{proof}

\subsection{Regular boundary points}
\begin{theorem}[The regular boundary alternative]\label{thm:balance-boundary}
Suppose $(z_*,\alpha_*)$ lies on the switching sheet $\Sigma_{12}$ of a regular switching
chart and $z_*\in\partial\Sh_P$. Then
\begin{equation}\label{eq:stationary}
 g_{12}(z_*,\alpha_*)=\partial_\alpha g_{12}(z_*,\alpha_*)=0.
\end{equation}
For two distinct ordinary saddles above the same $(z,\alpha)$, these two equations are equivalent to
\begin{equation}\label{eq:balance}
 |P(w_1)|=|P(w_2)|,\qquad |P'(w_1)|=|P'(w_2)|.
\end{equation}
Consequently, the observation point lies on the perpendicular bisector of $w_1$ and $w_2$.
\end{theorem}
\begin{proof}
If $\partial_\alpha g_{12}(z_*,\alpha_*)\ne0$, the implicit-function theorem solves the
switching equation as $\alpha=\alpha(z)$ throughout a neighborhood of $z_*$. By the
switching-chart hypothesis, the corresponding fixed-ratio interface belongs to the
fixed-ratio support for every nearby $z$. Its projection therefore contains a neighborhood
of $z_*$ in $\Sh_P$, contrary to $z_*\in\partial\Sh_P$.

For the algebraic reduction, set
\[
 \mathfrak d_0:=\log\left|\frac{P(w_1)}{P(w_2)}\right|,
 \qquad
 \mathfrak d_1:=\log\left|\frac{P'(w_1)}{P'(w_2)}\right|.
\]
By \eqref{eq:theta}, the switching equation is
$(1-\alpha)\mathfrak d_0+\alpha\mathfrak d_1=0$, while stationarity is
$\mathfrak d_1-\mathfrak d_0=0$. The determinant of this linear system is one, so $\mathfrak
d_0=\mathfrak d_1=0$. The converse is immediate. Equation~\eqref{eq:g-alpha} then gives
$|w_1-z|=|w_2-z|$.
\end{proof}

This result concerns regular preimages with $0<\alpha<d$.
A global description of the boundary must also include saddle collisions,
endpoint ratios, and \emph{multiple-height transitions}, at which more than
two connecting saddles have the same critical height.

\subsection{Caustics and balanced secants}
\begin{figure}[t]
\centering
\includegraphics{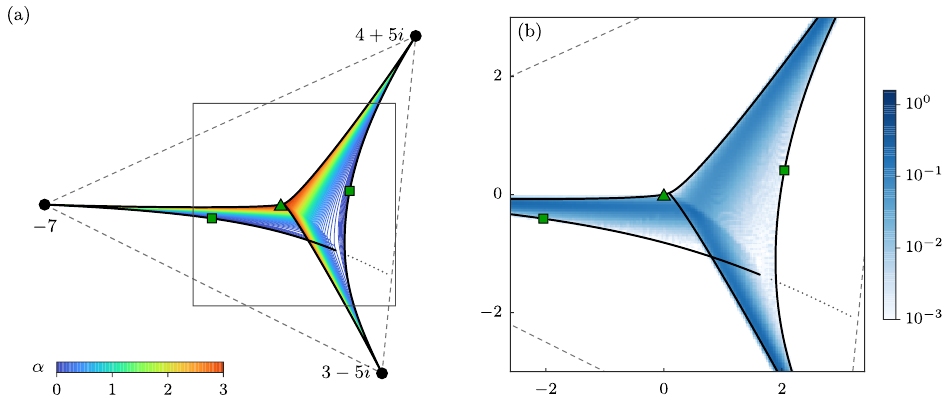}
\caption{(a)~The supports $\FixSh_{\alpha,P}$, $0<\alpha<3$, for the
cubic of Figure~\ref{fig:regimes}(a), colored by $\alpha$. Black:
critical values of $F_{P,\alpha}$ that are endpoints of
$\FixSh_{\alpha,P}$; dotted: the remaining caustic candidates; dashed:
$\partial\Hull_P$. The boundary of $\Sh_P$ lies on the black curves,
some of which also cross its interior.
(b)~The density of $\bar\nu_{P,n}$, $n=200$, in the box of~(a), on a
logarithmic scale; it changes abruptly across some caustic arcs inside
the core.}
\label{fig:candidates}
\end{figure}

The spatial projection of the ordinary saddle-collision locus is
called a \emph{caustic}. Its incidence equations are
\begin{equation}\label{eq:caustic-inc}
 z=F_{P,\alpha}(w),\qquad 1-\alpha h_P'(w)=0,
 \qquad 0<\alpha<d.
\end{equation}
It is parametrized by
\begin{equation}\label{eq:caustic-param}
 \alpha=\frac1{h_P'(w)},\qquad z=w-\frac{h_P(w)}{h_P'(w)}.
\end{equation}
Only values for which $1/h_P'(w)$ is real and belongs to $(0,d)$ are retained. The closure
of the resulting $z$-image is the \emph{caustic candidate set}. Its points are possible
collision values, not automatically selected endpoints; see
Figure~\ref{fig:candidates}(a). Since $h_P=\RedQ_P/\RedN_P$, the reduced
critical equation is
\begin{equation}\label{eq:reduced-critical}
 \RedN_P(w)^2-\alpha\bigl(\RedQ_P'(w)\RedN_P(w)-\RedQ_P(w)\RedN_P'(w)\bigr)=0,
\end{equation}
on $\RedQ_P(w)\RedN_P(w)\ne0$.

For a monic polynomial with roots $b_1,\ldots,b_N$, the discriminant
$\prod_{i<j}(b_i-b_j)^2$ is a polynomial in its coefficients and vanishes
exactly when the polynomial has a repeated root. Working with the reduced
rational function $h_P$ is essential when $P$ has multiple roots:
an unreduced polynomial discriminant can acquire components unrelated to
collisions of branches. Related geometry of critical values is developed in \cite{HSS}.

For distinct $w_1,w_2$ satisfying $h_P(w_1)\ne h_P(w_2)$, define the \emph{secant ratio} and
\emph{secant image} by
\begin{equation}\label{eq:secants}
 \alpha_{\mathrm{sec}}(w_1,w_2):=\frac{w_1-w_2}{h_P(w_1)-h_P(w_2)},\qquad
 z_{\mathrm{sec}}(w_1,w_2):=\frac{w_2h_P(w_1)-w_1h_P(w_2)}{h_P(w_1)-h_P(w_2)}.
\end{equation}
The equal-image condition
$F_{P,\alpha}(w_1)=F_{P,\alpha}(w_2)=z$ is equivalent to
$\alpha=\alpha_{\mathrm{sec}}(w_1,w_2)$ and
$z=z_{\mathrm{sec}}(w_1,w_2)$. The
\emph{off-diagonal balance candidate set} is the $z_{\mathrm{sec}}$-image of
\begin{equation}\label{eq:balanced-secants}
 \left\{(w_1,w_2):\begin{aligned}
 &w_1\ne w_2,\quad P(w_j)P'(w_j)\ne0,\quad
 \alpha_{\mathrm{sec}}(w_1,w_2)\in(0,d),\\
 &|P(w_1)|=|P(w_2)|,\qquad |P'(w_1)|=|P'(w_2)|
 \end{aligned}\right\}.
\end{equation}
After denominators are cleared while retaining their nonvanishing conditions, every
condition in \eqref{eq:balanced-secants} is semialgebraic. The candidate set is therefore
semialgebraic. Selection and exposure nevertheless remain separate requirements: the
candidate set is generally larger than the boundary of the core.

For unconstrained secants, the diagonal limit $w_1,w_2\to w$ gives
\[
 \alpha_{\mathrm{sec}}(w_1,w_2)\longrightarrow\frac1{h_P'(w)},\qquad
 z_{\mathrm{sec}}(w_1,w_2)\longrightarrow w-\frac{h_P(w)}{h_P'(w)}.
\]
These are the tangent limits of the secant formulas for the rational
graph of $h_P$. The equal-modulus constraints may prevent the limit
from being attained by balanced pairs. A caustic point therefore need
not be a limit of off-diagonal balanced secants.

\section{Semialgebraicity and the algebraic boundary}
\label{sec:algebraic}

The calculation at regular boundary points gives algebraic candidates for
the boundary, but a candidate need not be selected. We first extend the
calculation to all boundary points, including multiple-height transitions.
The resulting finite algebraic graph contains the boundary. The conservation
law then shows that the core contains either a whole root-free edge of
this graph or none of it. Finally, continuity of the core in Hausdorff
distance makes this selection constant on the pieces of a compatible
trivialization over coefficient space.

\begin{theorem}[The global boundary alternative]\label{thm:global-boundary}
Put
\[
 \mathscr D_P:=Z(P)\cup Z(P')\cup\{\centroid_P\}.
\]
For every $z\in\partial\Sh_P\setminus\mathscr D_P$ there is an $\alpha\in(0,d)$ for which at
least one of the following holds:
\begin{enumerate}[label=\textup{(\roman*)},leftmargin=2.2em]
\item An ordinary connecting saddle $w$ satisfies
\[
 z=w-\alpha h_P(w),\qquad 1-\alpha h_P'(w)=0.
\]
\item Two distinct simple ordinary connecting saddles $w_1,w_2$ have the same image $z$ and satisfy
\[
 |P(w_1)|=|P(w_2)|,\qquad |P'(w_1)|=|P'(w_2)|.
\]
\end{enumerate}
In the second case, all the connecting saddles have the same distance from $z$. Thus a
multiple-height transition produces no additional, nonalgebraic boundary mechanism.
\end{theorem}
\begin{proof}
By \eqref{eq:core-union}, choose $z_j\in\FixSh_{\alpha_j,P}$ with $z_j\to z$, and pass
to a convergent subsequence of the ratios. The endpoint limits in
Proposition~\ref{prop:global-facts} exclude the limits $0$
and $d$ because $z\notin\mathscr D_P$. Lemma~\ref{lem:incidence-closed} gives
$z\in\FixSh_{\alpha,P}$ at the limiting ratio $\alpha\in(0,d)$.

If a connecting saddle is multiple, (i) holds. Otherwise
Proposition~\ref{prop:selection-facts}\textup{(ii)}
writes the barrier near $(z,\alpha)$ as
\[
 \Blev_P(\zeta,\beta)=\max_{1\le j\le k}\psi_j(\zeta,\beta),
\]
where $\psi_1,\ldots,\psi_k$ are the critical heights of the simple
connecting saddles. There must be at least two, since a single branch would
make the potential harmonic near $z$. At the base point all $\psi_j$ have the same value.
Set
\[
 b_j:=\partial_\beta\psi_j(z,\alpha)
          =-\log|w_j-z|.
\]
Suppose these numbers are not all equal. For a sufficiently small fixed $\eps>0$, the
maximum at $(z,\alpha+\eps)$ is attained only among the indices having largest $b_j$,
and at $(z,\alpha-\eps)$ only among those having smallest $b_j$. The two index groups
are disjoint. The strict inequalities excluding the other indices persist for all $\zeta$ in
a small disk about $z$.

For each such $\zeta$, as $\beta$ crosses this interval, the maximizing index cannot remain
unique throughout: a unique maximizing index is locally constant, whereas the endpoint
groups are disjoint. At some ratio two distinct branches therefore attain the maximum. Their
spatial gradients are distinct, since
\[
 2\partial_\zeta\psi_i-2\partial_\zeta\psi_j
   =\beta\left(\frac1{w_i-\zeta}-\frac1{w_j-\zeta}\right)\ne0.
\]
The maximum is not harmonic at that spatial point. Indeed, if it were harmonic on a
neighborhood, its difference from either touching harmonic branch would be a nonnegative
harmonic function with an interior zero; both branches would then agree, contrary to the
displayed gradient inequality. Thus $\zeta$ belongs to the fixed-ratio support at that
ratio. The entire disk about $z$ would then be contained in $\Sh_P$, contradicting
$z\in\partial\Sh_P$.

All $b_j$ are consequently equal. All distances $|w_j-z|$ are equal, and equality of the
critical heights gives equality of $|P(w_j)|$. The saddle equations then give equality of
$|P'(w_j)|$. Any two distinct connecting saddles supply (ii).
\end{proof}

\subsection{A finite algebraic candidate graph}
Let $\mathscr G_P^{\mathrm{col}}$ be the $z$-image of the solutions of the ordinary collision
equations \eqref{eq:caustic-inc} with $P(w)P'(w)\ne0$. Let $\mathscr G_P^{\mathrm{bal}}$ be the
off-diagonal balance candidate set, that is, the $z_{\mathrm{sec}}$-image of
\eqref{eq:balanced-secants}. Define
\begin{equation}\label{eq:candidate-graph}
 \mathscr G_P:=\partial \Hull_P\ \cup\ \mathscr D_P\ \cup\
 \left(\Hull_P\cap\overline{\mathscr G_P^{\mathrm{col}}
                                \cup\mathscr G_P^{\mathrm{bal}}}\right).
\end{equation}
Closure here is taken in the observation plane for the fixed polynomial. The added hull
boundary is convenient for a finite planar decomposition; it is not asserted to be part of
the boundary of the core.

\begin{proposition}[One-dimensional algebraic candidates]\label{prop:algebraic-graph}
For fixed degree the family $\{(P,z):z\in\mathscr G_P\}$ is semialgebraic. Each fiber is
compact, has dimension at most one, and contains $\partial\Sh_P$. There is an integer $D_d$,
depending only on $d$, such that each $\mathscr G_P$ is contained in the zero set of a
nonzero real polynomial of degree at most $D_d$.
\end{proposition}
\begin{proof}
The collision and balanced-pair equations are polynomial equalities and inequalities in the
real and imaginary parts of their variables after nonzero denominators have been cleared.
This assertion uses the original $P,P'$ only at ordinary saddles, so collisions of roots
introduce no spurious ordinary solutions. Projection, fiberwise closure, and intersection
preserve semialgebraicity. The root convex hull and its boundary are semialgebraic families
as well: one can quantify over a list of roots and their convex-combination coefficients.
These facts prove the first assertion; compactness is immediate from
\eqref{eq:candidate-graph}. Theorem~\ref{thm:global-boundary} proves the boundary inclusion.

For fixed $\alpha$, the rational map $F_{P,\alpha}$ has finitely many critical values, so
the collision incidence in $(z,\alpha)$ has dimension at most one. Its spatial image has
dimension at most one. It remains to check that a symmetry cannot make the balanced image
two-dimensional.

Away from collisions, choose two simple inverse branches locally, and let $g$ be the
difference of their critical heights. The set of $(z,\alpha)$ at which these branches form a
balanced pair is contained in
\[
 g(z,\alpha)=0,\qquad \partial_\alpha g(z,\alpha)=0,
 \qquad \nabla_z g\ne0,
\]
by Proposition~\ref{prop:dictionary} and the equivalence of \eqref{eq:stationary} and
\eqref{eq:balance}. The regular level set $g=0$ is a smooth surface. At each stationary point,
the vertical ratio direction lies in its tangent plane, so the derivative
of the spatial projection has rank at most one. By Sard's theorem in local
charts of the surface, the image of the stationary points has zero area. Countably many
such charts cover all balanced pairs of simple saddles. Pairs
involving a multiple saddle project into the collision image already
considered. Thus the semialgebraic balanced image has empty interior and
dimension at most one. Closure preserves semialgebraic dimension.

For the degree bound, quantifier elimination gives a fixed-degree quantifier-free
description of the family \eqref{eq:candidate-graph} using finitely many real polynomials in
the coefficient parameters and $x,y$. Fix $P$ and discard the polynomials that specialize
identically to zero. At any point of $\mathscr G_P$, at least one remaining polynomial
vanishes; otherwise all signs in that description would be locally constant and the fiber
would have interior. The product of the remaining nonconstant polynomials is nonzero and
vanishes on the fiber. Their number and degrees are bounded in terms of the original degree
$d$ only. The fiber is nonempty and not the whole plane, so a nonconstant factor is present.
This proves the assertion with their total degree as a uniform bound. Standard quantifier
elimination and dimension facts are recalled in \cite{BCR}.
\end{proof}

\subsection{Continuation of selection along an algebraic edge}
An arc contained in an open set is a \emph{crossing arc} of that set if
there is a homeomorphism from the set onto an open rectangle that maps the
arc onto the horizontal middle interval of the rectangle.
This is a topological condition; no differentiability is required.

\begin{lemma}[Continuation along a root-free candidate arc]\label{lem:whole-edge}
Let $e$ be a connected open semialgebraic arc of $\mathscr G_P$,
disjoint from $Z(P)$. Suppose every point of $e$ has a simply connected
neighborhood whose intersection with $\mathscr G_P$ is a single crossing arc
of that neighborhood. Then either $e\subseteq\Sh_P$ or
$e\cap\Sh_P=\varnothing$. The arc need not be differentiable.
\end{lemma}
\begin{proof}
Since $\Sh_P$ is closed, it suffices to prove that $e\cap\Sh_P$
is relatively open in $e$. Take $z\in e\cap\Sh_P$ and a small root-free
neighborhood $D$ of $z$ as in the statement. Each of the two
components of $D\setminus\mathscr G_P$ lies wholly in the core or
wholly outside it, since $\partial\Sh_P\subseteq\mathscr G_P$.
If one of them lies in the core, then so does the crossing arc, because the
core is closed.
Otherwise $D\cap\Sh_P$ is contained in the crossing arc.
If some point of that arc were absent from the core, the gap around it would
join the two sides, so $D\setminus\Sh_P$ would be connected. On the simply
connected root-free domain $D$, choose a holomorphic primitive $\mathscr P_D$
of $-\rootC_P^2$. The conservation law would give
$C_{\rho_P}=\mathscr P_D+c$ on
this connected complement for one constant $c$. The crossing
arc has planar area zero, since it is one-dimensional and
semialgebraic; no smoothness is needed. The equality would thus hold
almost everywhere in $D$, and applying the distributional
$\partial_{\bar z}$ would give $\rho_P(D)=0$, contradicting
$z\in\supp\rho_P$. Hence the entire crossing arc lies in the core. This proves
relative openness, and the assertion follows by connectedness.
\end{proof}

\begin{theorem}[Semialgebraicity of the core]\label{thm:semialgebraic}
For every fixed degree $d$, the family
\[
 \{(P,z):P\text{ monic of degree }d,\ z\in\Sh_P\}
\]
is semialgebraic. Each boundary $\partial\Sh_P$ consists of finitely many real-algebraic arcs and
points and lies on a real algebraic curve whose degree is bounded in terms of $d$ only. These
conclusions concern $\Sh_P$, not the possibly infinite set of transients in $\AllSh_P$.
\end{theorem}
\begin{proof}
We begin with a fixed polynomial. Decompose the compact semialgebraic
graph $\mathscr G_P$ into finitely many vertices and open embedded
$C^1$ edges. Include $\mathscr D_P$ among the vertices, and split
a closed loop by adding a vertex if necessary. The complement of the graph has finitely many
connected
open semialgebraic components, which we call \emph{faces}. Each face is either entirely
contained in $\Sh_P$ or disjoint from it, since $\partial\Sh_P$ lies in the graph. By
Lemma~\ref{lem:whole-edge}, the same is true of each open edge. The core is therefore the
union of finitely many faces, open edges, and vertices, each of which is semialgebraic.
Thus each core is semialgebraic, although the comparisons of critical
heights used to select it involve logarithms.

For uniformity in the coefficients, first compress each ambient plane into the
open square $(-1,1)^2$ by applying $x\mapsto x/(1+|x|)$ to each coordinate.
The candidate graphs and distinguished root sets remain compact semialgebraic
families. Semialgebraic triviality for compact families \cite{BCR,Hardt},
applied compatibly with these families and with the closed
square and its boundary, gives a finite partition of
coefficient space into connected semialgebraic pieces. On each piece, the
trivialization preserves the boundary of the square and therefore restricts
to the open square; conjugating by the
compression gives a semialgebraic family of homeomorphisms of the plane over
the piece.
Use these homeomorphisms to transport the finite cell decomposition of one
sample fiber, with its roots among the vertices, to all fibers over the
piece. The inclusion $\partial\Sh_P\subseteq\mathscr G_P$ and
Lemma~\ref{lem:whole-edge} make every transported cell wholly included in the
core or wholly excluded from it. Crossing arcs are preserved by
homeomorphisms of the plane, so no
smoothness of the transported edges is required.

Only finitely many unions of cells of the sample fiber are possible; we call
them \emph{models}. The inverse image
of each core under the homeomorphism of its fiber is compact, so retain
only the compact models and identify any repeated ones. Their transported
images are continuous families of compact sets, pairwise distinct at every
parameter because the homeomorphisms are injective. The core itself depends
continuously on $P$ in Hausdorff distance by Proposition~\ref{prop:global-facts}. For each
model, equality with the core therefore defines a closed set of parameters.
These sets form a finite disjoint cover of the piece; the complement of each
is a finite union of closed sets, so each of them is also open. By
connectedness, one fixed model equals the core throughout the piece. That model
is a finite union of semialgebraic cells, and its image under the
semialgebraic trivialization is a semialgebraic family. Taking the finite
union over the pieces proves the assertion about the family.

Finally, Proposition~\ref{prop:algebraic-graph} supplies the degree-bounded algebraic curve
containing the boundary. A one-dimensional semialgebraic set admits a finite decomposition
into real-analytic algebraic arcs and points. Applying this to the boundary proves the last
geometric assertion.
\end{proof}

\begin{corollary}[Finite topology in each degree]\label{cor:finite-topology}
For fixed $d$, the cores have finitely many homeomorphism types and
uniformly bounded Betti numbers. Each core is homeomorphic to a finite
simplicial complex and has finitely many complementary components.
\end{corollary}
\begin{proof}
Apply semialgebraic triviality to the coefficient family in
Theorem~\ref{thm:semialgebraic}, and triangulate its compact fibers
\cite{BCR,Hardt}. Apply the same finite decomposition to the complements.
\end{proof}

\begin{corollary}[A degree-uniform boundary-length bound]\label{cor:uniform-length}
Let $D_d$ be as in Proposition~\ref{prop:algebraic-graph}. If $\Hull_P$ lies in a closed disk of
radius $R$, then
\[
 \Haus^1(\partial\AllSh_P)=\Haus^1(\partial\Sh_P)
       \le \pi D_d R.
\]
In particular, once the diameter of the root set is bounded, the length of the boundary is
bounded uniformly in each degree, including at degenerate polynomials.
\end{corollary}
\begin{proof}
The boundary $\partial\Sh_P$ is rectifiable and contained in an algebraic curve of degree at most
$D_d$. For almost every line it has at most $D_d$ intersection points, by restricting the
defining polynomial to that line; lines contained in the curve form a null exceptional
family. The planar Cauchy--Crofton formula \cite{Federer}, integrating unoriented normal
angle over $[0,\pi)$ and signed offset over $[-R,R]$, bounds its length by
$\tfrac12\pi(2R)D_d$. The transient set is countable and has zero $\Haus^1$ measure, so the
boundary decomposition in Proposition~\ref{prop:global-facts} gives the same bound for the full
shadow. A singleton hull is immediate.
\end{proof}

\begin{corollary}[Local H\"older stability of the core]\label{cor:holder-shape}
Write $P(z)=z^d+\sum_{k=0}^{d-1}c_kz^k$ and use the Euclidean norm of
$(c_0,\ldots,c_{d-1})\in\C^d$ as the coefficient norm $\|P\|_{\mathrm{coef}}$. For each
$R>0$ there are $C_R>0$ and $0<\gamma_R\le1$, depending also on $d$, such that
\[
 d_{\mathrm H}(\Sh_P,\Sh_Q)
       \le C_R\|P-Q\|_{\mathrm{coef}}^{\gamma_R}
 \quad\text{if }\|P\|_{\mathrm{coef}},\|Q\|_{\mathrm{coef}}\le R.
\]
\end{corollary}
\begin{proof}
Hausdorff distance between compact semialgebraic fibers is a semialgebraic function: its
defining distance inequalities involve only quantification over the two compact fibers. By
Proposition~\ref{prop:global-facts}, this function is continuous. On the compact semialgebraic
product of two coefficient balls, apply the semialgebraic \L ojasiewicz inequality
\cite{BCR} to this nonnegative function and $\|P-Q\|_{\mathrm{coef}}$. The latter vanishes
only on the diagonal, where the former vanishes too. Thus a positive integer power of the
Hausdorff distance is bounded by a constant times the coefficient distance, yielding the
assertion.
\end{proof}

\begin{remark}[Selection and degree bounds]
The equations for the candidate graph have degree controlled by $d$,
independently of the power level $n$. The core is a union of whole faces,
edges, and vertices of that graph. This description supplies neither
a small explicit bound for $D_d$ nor a list of the included cells;
selection is still determined by the barrier. In particular, the union
of all collision and balanced-pair candidates is not in general the
boundary; in Figure~\ref{fig:candidates}(a), some caustic arcs cross
the interior of the core.
\end{remark}

\subsection{Consequences for the full shadow}
\begin{corollary}[Boundary size of the full shadow]\label{cor:literal-boundary}
For every nonconstant polynomial,
\[
 \Haus^1(\partial\AllSh_P)=\Haus^1(\partial\Sh_P)<\infty.
\]
The boundary of the core is a finite union of points and bounded real-algebraic arcs. The
boundary of the full shadow adds at most countably many points, all isolated and outside $\Sh_P$.
\end{corollary}
\begin{proof}
Combine Theorem~\ref{thm:semialgebraic} with Proposition~\ref{prop:global-facts}; a countable set has
zero $\Haus^1$ measure. Corollary~\ref{cor:uniform-length} gives a degree-uniform
quantitative bound.
\end{proof}

\begin{corollary}[When the full shadow is connected or semialgebraic]
\label{cor:literal-topology}
The full shadow $\AllSh_P$ is connected if and only if $T_P$ is empty.
It is semialgebraic if and only if $T_P$ is finite. Equivalently,
it is definable, with real parameters, in some o-minimal expansion
of the real field if and only if $T_P$ is finite. Its connected
components are $\Sh_P$ and the singletons $\{z\}$ for $z\in T_P$.
\end{corollary}
\begin{proof}
The core is nonempty and connected, whereas every point of
$T_P$ is isolated in the full shadow. This gives the component
statement. If $T_P$ is finite, its union with the semialgebraic set
$\Sh_P$ is semialgebraic, hence definable in every o-minimal expansion of
the real field. Conversely, if the full shadow is
definable in some o-minimal expansion, then
$T_P=\AllSh_P\setminus\Sh_P$ is a definable discrete set.
Such a set has finitely many connected components and is therefore
finite.
\end{proof}

\section{Reflection rigidity and the planar body}\label{sec:rigidity}
We first use the conservation law to show that an isolated arc of the core
is straight. A stronger statement holds for a support arc that remains
stationary over an interval of differentiation ratios: it forces a
reflection symmetry of the entire root multiset.

For $s\ge2$, let $\mathscr A_P$ be the union of the reflection axes of
$P$, with $\mathscr A_P=\varnothing$ when there are none. There are only
finitely many axes, and all pass through $\centroid_P$. Indeed, a
reflection preserving the root measure fixes its centroid; it acts on the
finite root set, which has at least two points, by one of finitely many
permutations, and two such reflections with the same action coincide.

\begin{definition}[Regular filament]\label{def:filament}
A \emph{regular filament} is a set $F=D\cap\Sh_P$, where $D$ is an open disk disjoint
from $Z(P)$ and $F$ is a single embedded regular $C^1$ arc crossing $D$ and separating it
into two components.
\end{definition}
Here and below, an arc is \emph{embedded} if, in suitable local coordinates of the plane, it
is a coordinate axis near each of its interior points and a closed half-axis near each of its
endpoints.

\subsection{Straight filaments of constant density}
\begin{theorem}[Filament rigidity]\label{thm:filament}
Every regular filament is a straight segment and carries a constant positive density with
respect to arclength. An embedded arc contained in the core cannot have a free endpoint away
from the roots. Here a \emph{free endpoint} of such an arc is an endpoint contained in an open
disk, disjoint from $Z(P)$, in which
the core consists only of that terminating arc.
\end{theorem}
\begin{proof}
On the root-free disk, choose a holomorphic primitive $\mathscr P_D$ of
$-\rootC_P^2$. On the two sides of the arc, Theorem~\ref{thm:Cprime} gives
\[
 C_{\rho_P}=\mathscr P_D+\kappa_+
 \quad\text{and}\quad
 C_{\rho_P}=\mathscr P_D+\kappa_-
\]
for constants $\kappa_+,\kappa_-$. They are distinct. Otherwise
$C_{\rho_P}$ would agree almost everywhere with one holomorphic function throughout $D$, and
its distributional $\partial_{\bar z}$ would vanish there, contrary to $\partial_{\bar
z}C_{\rho_P}=\pi\rho_P$ and $F\subseteq\Sh_P=\supp\rho_P$.

Across the oriented arc, the distributional jump of $C_{\rho_P}$ gives the measure
\[
 \frac{\kappa_--\kappa_+}{2\pi i}\dd z,
\]
up to reversal of orientation. Positivity forces the tangent direction to be constant. Hence
the arc is straight and has arclength density
$|\kappa_--\kappa_+|/(2\pi)$. More explicitly, near each $z_0\in F$ the support of $\rho_P$
coincides with the line $\Ree((\kappa_+-\kappa_-)(z-z_0))=0$.
By connectedness, $\Ree((\kappa_+-\kappa_-)z)$ is constant on $F$, so this is the same line
for every $z_0$. The filament is relatively closed
and locally open in this line inside $D$, so $F$ is the whole chord of $D$ on this line.

Suppose now that an embedded arc in the core has a free endpoint, and choose a
root-free disk in which the core consists only of the terminating arc. The complement of the
arc in the
disk is connected, so there is only one additive constant. Then $C_{\rho_P}$ would agree
almost everywhere on the entire
disk with a holomorphic function and its distributional $\partial_{\bar z}$ would vanish,
which is impossible because the arc belongs to $\supp\rho_P$.
\end{proof}

\subsection{Stationary arcs and reflection}
\begin{theorem}[A stationary arc forces reflection]\label{thm:stationary-reflection}
Suppose a regular switching chart contains a product $\Gamma\times I$ in its switching
sheet, where $\Gamma$ is an embedded nonconstant $C^1$ spatial arc and $I\subset(0,d)$ is an
open interval. Then $\Gamma$ lies on a reflection axis of the root multiset.
\end{theorem}
\begin{proof}
Fix $z\in\Gamma$ and a unit tangent $\mathbf v\in\C$ at $z$. Write
$\xi_j=w_j(z,\alpha)-z$ and $g=\psi_1-\psi_2$.
Since $g=0$ on $\Gamma\times I$, its derivative in $\alpha$ vanishes there.
Equation~\eqref{eq:g-alpha} gives $|\xi_1|=|\xi_2|=:r>0$. Its tangential derivative also
vanishes, and \eqref{eq:g-z} gives
\[
 0=\Ree\bigl((2\partial_zg)\mathbf v\bigr)
   =\frac{\alpha}{r^2}\Ree\bigl(\overline{w_1-w_2}\,\mathbf v\bigr).
\]
Thus $w_1-w_2$ is perpendicular to the tangent line through $z$, while $z$ is equidistant
from the two saddles. Reflection in that line,
\[
 \mathfrak r(w)=z+\mathbf v^2\overline{w-z},
\]
therefore interchanges $w_1$ and $w_2$ for every $\alpha\in I$.

The saddle equations $\xi_j=\alpha h_P(w_j)$ imply
\begin{equation}\label{eq:h-reflection-curve}
 h_P(\mathfrak r(w_1(z,\alpha)))=
 \mathbf v^2\overline{h_P(w_1(z,\alpha))}.
\end{equation}
The branch $\alpha\mapsto w_1(z,\alpha)$ is injective: equality at two ratios,
together with $w_1-z=\alpha h_P(w_1)$ and $h_P(w_1)\ne0$, forces the ratios to be equal.
Its image is therefore infinite. Set
$Q(w):=\overline{P(\mathfrak r(w))}$, which is a polynomial in $w$.
Conjugating \eqref{eq:h-reflection-curve} and clearing denominators gives
$Q'P=QP'$ on this infinite image, hence identically as a polynomial identity.
In particular,
\begin{equation}\label{eq:h-reflection-global}
 h_P(\mathfrak r(w))=\mathbf v^2\overline{h_P(w)}
\end{equation}
identically as a rational identity. Moreover, $(Q/P)'=0$, so $Q/P$ is a nonzero constant.
Thus reflection preserves the entire root multiset, including multiplicities.
The tangent line at $z$ is a reflection axis.

The nonconstant arc lies in $\Hull_P$, so there are at least two distinct roots.
A map of the form $w\mapsto u\overline w+b$ is determined by its values on two
such roots; preserving the finite root set leaves only finitely many possibilities.
As $z$ varies along the connected arc, the parameters of the tangent reflection vary
continuously among these finitely many possibilities and are therefore constant. Hence the
entire arc, including its endpoints
by continuity, lies on one reflection axis.
\end{proof}

\begin{theorem}[All curve mass lies on reflection axes]\label{thm:axes}
For any polynomial with at least two distinct roots, the curve part of $\rho_P$ is supported
on $\mathscr A_P$. In particular, if the root multiset has no reflection axis, then
\begin{equation}\label{eq:no-axis-ac}
 \rho_P=f_P\dd A+\sum_{j=1}^s p_j^2\delta_{a_j},
 \qquad
 \int_\C f_P\dd A=1-\sum_{j=1}^s p_j^2>0.
\end{equation}
In this case $\Sh_P$ is the closure of its interior.
\end{theorem}
\begin{proof}
By Theorem~\ref{thm:decomposition}, every contribution to the curve part comes from a
two-dimensional
stratum of $\Inc_P$ with nonzero ratio differential and spatial rank one. In the coordinates
used there, such a stratum is locally the product of a regular stationary support arc with an
open ratio
interval. Fix an interior ratio $\alpha$ and discard the finitely many roots and critical
values of
the rational map $F_{P,\alpha}$. Near each remaining observation point, the connecting saddles
continue analytically and the barrier is
the maximum of their finitely many critical heights. At every support point at least two of
these heights attain the maximum: a unique
maximizer would make the potential locally harmonic. The finitely many closed
sets on which two of the heights are tied therefore cover the product, and the union of their
interiors is dense by the
Baire category theorem. On each open patch on which two heights are tied, the calculation in
the proof of
Theorem~\ref{thm:stationary-reflection}, which uses the saddle equations and equality of
heights, gives a tangent reflection preserving the entire root multiset. The finite union
of axes is closed, so it contains the whole arc, including the discarded observation points.
The remaining lower-dimensional strata of $\Inc_P$ carry no mass of the nonatomic part
$\Omega_P^\circ$ of the incidence measure.
This proves the assertion about the support of the curve part.

Without axes the curve part vanishes, giving \eqref{eq:no-axis-ac}; the planar mass
is positive since $0<p_j<1$ and $\sum_jp_j^2<\sum_jp_j=1$.
The boundary of the semialgebraic core has zero area. A point of the core that is not a root
and lies outside the closure
of the interior would have a root-free neighborhood disk meeting the core only in its boundary;
absolute continuity would give the disk zero shadow mass, a contradiction. A connected set
containing two distinct roots has no isolated points, so roots are limits of non-root
points of the core as well. The reverse inclusion follows from closedness.
\end{proof}

Recall the planar body $\mathcal B_P^{\mathrm{pl}}$ and the axial set
$\mathcal L_P$ defined in the introduction.
The next theorem proves \eqref{eq:full-body-intro} and describes the
geometry and the measure of the part of the core outside its planar body.

\begin{theorem}[The planar body and all one-dimensional appendages]\label{thm:body-segments}
Suppose $P$ has at least two distinct roots, and consider the closed sets
\[
 \mathcal B_P^{\mathrm{pl}}:=\overline{\operatorname{int}\Sh_P}
       =\overline{\operatorname{int}\AllSh_P},\qquad
 \mathcal L_P:=\overline{\Sh_P\setminus\mathcal B_P^{\mathrm{pl}}}.
\]
The planar body is exactly the support of the absolutely continuous part:
\[
 \mathcal B_P^{\mathrm{pl}}=\supp(f_P\dd A).
\]
The set $\mathcal L_P$ is empty or a finite union of nondegenerate closed
straight segments lying on reflection axes of the root multiset.
The intersection $\mathcal L_P\cap\mathcal B_P^{\mathrm{pl}}$ is finite, and
\begin{equation}\label{eq:body-segment-decomposition}
 \Sh_P=\mathcal B_P^{\mathrm{pl}}\cup\mathcal L_P,
 \qquad
 \AllSh_P=\mathcal B_P^{\mathrm{pl}}\cup\mathcal L_P\cup T_P.
\end{equation}
The segments can be subdivided so that their endpoints belong to
$\mathcal B_P^{\mathrm{pl}}\cup Z(P)\cup\{\centroid_P\}$. On every open segment away from
this set the arclength density of $\rho_P$ is constant and positive. Every root outside the
union of reflection axes belongs to the planar body. If there is no reflection axis,
$\mathcal L_P=\varnothing$.
\end{theorem}
\begin{proof}
Both $\Sh_P$ and $\mathcal B_P^{\mathrm{pl}}$ are semialgebraic.
Every disk of positive $f_P\dd A$-mass meets the semialgebraic set
$\Sh_P$ in positive area, and hence meets its interior. Thus
$\supp(f_P\dd A)\subseteq\mathcal B_P^{\mathrm{pl}}$.
Conversely, a disk meeting $\operatorname{int}\Sh_P$ contains a smaller
disk in that interior avoiding the finitely many reflection axes and
the roots. As a neighborhood of a support point it has positive
$\rho_P$-mass, which is entirely planar by Theorem~\ref{thm:axes}.
Taking closures proves the reverse inclusion.
If a disk avoids the planar body, its intersection with $\Sh_P$ has empty interior and hence
zero planar area. On such a disk Theorem~\ref{thm:axes} therefore places the measure on the
reflection axes and the finitely many roots. A support point outside those axes could only
be a root isolated in the support, which connectedness excludes. Consequently
\[
 \Sh_P\setminus\mathcal B_P^{\mathrm{pl}}\subseteq\mathscr A_P.
\]
This difference is a bounded semialgebraic subset of finitely many lines. It has no isolated
points: it is relatively open in $\Sh_P$, which has none. Its closure is thus a finite union
of nondegenerate closed line segments, unless empty. This proves
\eqref{eq:body-segment-decomposition}, using Proposition~\ref{prop:global-facts} for the full shadow.

Remove the centroid and intersect $\Sh_P\setminus\mathcal B_P^{\mathrm{pl}}$ with each axis.
Each of these sections has finite
boundary by semialgebraicity, so its closure adds only finitely many points; distinct axes
meet only at the centroid. Consequently $\mathcal L_P\cap\mathcal B_P^{\mathrm{pl}}$ is finite.
Split at these points, at roots, and at the centroid. Near any remaining support point, a
small root-free rectangle avoiding the body and centroid has support confined to one axis.
The conservation law gives a common analytic primitive for its Cauchy-transform derivative.
A gap on the middle line would connect the slit complement; analytic continuation would
then exclude the support point, a contradiction. Thus the local support is the whole line
section, excluding further endpoints. Theorem~\ref{thm:filament}, or the argument in its
proof, gives constant positive arclength density. Constants agree on
overlapping neighborhoods because their common support neighborhood has positive finite
shadow mass, so the density is constant on each connected open segment. The statements
about roots and the absence of axes follow from the displayed inclusion.
\end{proof}

\begin{corollary}[Generic absolute continuity in every degree]\label{cor:generic-ac}
For $d\ge3$, the monic degree-$d$ polynomials with no reflection axis form an open dense
subset of coefficient space whose complement has Lebesgue measure zero. On that set
\eqref{eq:no-axis-ac} holds, and the two shadows have the same nonempty
interior.
\end{corollary}
\begin{proof}
Identify coefficient space with $\R^{2d}$. Write an axis as
$e^{i\theta}(\R+i\tau)$, with $\theta,\tau\in\R$. Reflection symmetry in this axis is
equivalent
to reality of the coefficients of the monic polynomial
$\widetilde P(w)=e^{-id\theta}P(e^{i\theta}(w+i\tau))$. Conversely,
$P(z)=e^{id\theta}\widetilde P(e^{-i\theta}z-i\tau)$ reconstructs every symmetric polynomial.
Thus the symmetric locus is the image of a smooth map with $d+2$ real
parameters: the angle, the normal offset, and the $d$ lower coefficients of $\widetilde P$.
On countably many compact parameter sets the derivatives are bounded, so the
image has $2d$-dimensional Lebesgue measure zero because $d+2<2d$ for $d\ge3$.
The locus is closed: every axis passes through the root centroid, centroids
converge under coefficient convergence, and a subsequence of unit axis
directions converges. Reality of the transformed coefficients passes to the
limit. The complement is therefore open, dense, and of full measure.
A polynomial without an axis has at least two distinct roots. By
Theorem~\ref{thm:axes}, its nonempty core is the closure of its interior, so that
interior is nonempty; Proposition~\ref{prop:global-facts} gives the same interior for the
full shadow.
\end{proof}

\begin{corollary}[The hole examples have genuine planar bodies]
\label{cor:hole-planar-body}
For every sufficiently large $N$, the core of
$P_N(z)=(z^2-1)^N(z-2-i)(z-3-i)$ is the closure of its nonempty
interior and has a hole. The same conclusions hold on a sufficiently
small open coefficient neighborhood consisting of polynomials with
no reflection axis.
\end{corollary}
\begin{proof}
For every sufficiently large $N$, the point
$q_\diamond=\tfrac23+\tfrac{i}{30}$ lies in a bounded component of
$\C\setminus\Sh_{P_N}$, and this separation persists on an open coefficient
neighborhood \cite[Theorem~5.4]{HS2}. These polynomials have no reflection axis.
Theorem~\ref{thm:axes} gives
$\Sh_{P_N}=\overline{\operatorname{int}\Sh_{P_N}}$.
The core is nonempty, so its interior is nonempty.
Intersect the neighborhood on which the hole persists with the open set of
polynomials without reflection axis from Corollary~\ref{cor:generic-ac}. The
intersection still contains the original polynomial, and both conclusions
hold for each of its members.
\end{proof}

\section{Segments, radial stars, and zero-area classification}
\label{sec:models}

The chord law gives explicit shadow measures for collinear and for radially
symmetric root configurations. In the first case, moments determine the
measure on the root segment. In the second, differentiation preserves
a finite radial star, and a sampling formula determines the distribution of
mass on it. In both cases every zero of every descendant already belongs to
the core.

\subsection{Collinear roots}
\begin{theorem}[The segment formula]\label{thm:collinear}
Suppose all roots of $P$ lie on one line, with distinct roots $a_1,\dots,a_s$ and weights
$p_j=r_j/d$. Then
\begin{equation}\label{eq:collinear}
 \rho_P=\sum_{j=1}^s p_j^2\delta_{a_j}
       +2\sum_{i<j}p_ip_j\Unif[a_i,a_j].
\end{equation}
Here $\Unif[a_i,a_j]$ is the normalized arclength measure on the segment $[a_i,a_j]$. If
$s\ge2$, then
\[
 \AllSh_P=\Sh_P=[a_-,a_+],
 \qquad \mathcal Z_{P,n}\longrightarrow[a_-,a_+]\quad\text{in Hausdorff distance},
\]
where $a_-,a_+$ are the extreme roots along the line. If $s=1$, both shadows equal the single root.
\end{theorem}
\begin{proof}
By affine covariance, reduce to real roots. Every descendant has only real roots, all in
$[a_-,a_+]$,
by repeated application of Rolle's theorem, so $\rho_P$ is supported in this interval. The
chord measure $\sigma_P$
is supported in that interval and has the same moments as $\rho_P$. Polynomials are
uniformly dense in the continuous functions on a compact interval; hence the two measures
coincide. Conditioning on the ordered pair of root endpoints in the chord law gives
\eqref{eq:collinear}. The pair of extreme roots gives positive density on the whole
interval, so this interval is the support. Every finite set $\mathcal Z_{P,n}$ is already
contained in the interval, while Proposition~\ref{prop:global-facts} shows that every point of the
interval is approached uniformly by $\mathcal Z_{P,n}$. This proves Hausdorff convergence
and equality with the full shadow.
\end{proof}

The density of the nonatomic part $\rho_P^{\mathrm{na}}$ of $\rho_P$ with respect to
arclength on the root line is piecewise constant:
\[
 \frac{\dd\rho_P^{\mathrm{na}}}{\dd(\Haus^1\!\restriction[a_-,a_+])}(z)
 =2\sum_{i<j}\frac{p_ip_j}{|a_i-a_j|}\ind_{[a_i,a_j]}(z),
\]
up to an irrelevant choice of values at the finitely many roots. For $P=z^2-1$ this becomes
\begin{equation}\label{eq:quadratic}
 \rho_P=\tfrac14\delta_{-1}+\tfrac14\delta_1
          +\tfrac14\ind_{[-1,1]}(x)\dd x.
\end{equation}
This constant density is an average over all differentiation ratios; it is
not the limiting density of zeros at a single fixed ratio.

\subsection{Preservation of signed radial stars}
Fix $q\ge2$, set $\omega_q=e^{2\pi i/q}$ and $\zeta_q=e^{\pi i/q}$, and suppose
\begin{equation}\label{eq:radial-class}
 P(z)=z^{r_0}H(z^q),\qquad r_0\ge0,
\end{equation}
where $H$ is monic and has only nonzero real roots. Roots of both signs are allowed.
Zero roots can be absorbed into $z^{r_0}$; if $H$ is constant, then $P$ is a power of $z$ and
has a singleton shadow. Assume henceforth that
$H$ is nonconstant, and put
\[
 b_+=\max\bigl(\{0\}\cup Z(H)\bigr),\qquad
 b_-=\min\bigl(\{0\}\cup Z(H)\bigr),\qquad
 R_+=b_+^{1/q},\quad R_-=(-b_-)^{1/q}.
\]
Define the \emph{signed radial star}
\begin{equation}\label{eq:signed-star}
 \mathscr R_{q,R_+,R_-}:=
 \bigcup_{j=0}^{q-1}\omega_q^j[0,R_+]
 \ \cup\
 \bigcup_{j=0}^{q-1}\zeta_q\omega_q^j[0,R_-].
\end{equation}
We call the segments in \eqref{eq:signed-star} the \emph{spokes} of the star.
A zero radius contributes only the origin. Positive roots of $H$ produce one
family of $q$ spokes; negative roots produce the interlaced family, rotated by $\pi/q$.

\begin{lemma}[Preservation of the signed radial star]\label{lem:star}
Every zero of every descendant of a polynomial \eqref{eq:radial-class},
including the descendants of order $m=0$, lies in $\mathscr R_{q,R_+,R_-}$.
\end{lemma}
\begin{proof}
Write a polynomial in the class as $z^\ell G(z^q)$, absorbing any zero roots of $G$ into
$z^\ell$, so that $G(0)\ne0$ unless $G$ is constant. Its derivative is
\[
 \begin{cases}
 z^{\ell-1}\bigl(\ell G(y)+qyG'(y)\bigr),&\ell>0,\\
 qz^{q-1}G'(y),&\ell=0,
 \end{cases}
 \qquad y=z^q.
\]
The roots of $G'$ remain real in the convex hull of the roots of $G$. For $\ell>0$,
apply Gauss--Lucas to $F(y):=y^\ell G(y)^q$. Its roots lie in
$[\min(\{0\}\cup Z(G)),\max(\{0\}\cup Z(G))]$, and
\[
 F'(y)=y^{\ell-1}G(y)^{q-1}\bigl(\ell G(y)+qyG'(y)\bigr).
\]
Every root of $\ell G+qyG'$ is therefore a root of $F'$ and lies in the same
real interval. Constant $G$ causes no difficulty. Thus the real-rooted class
and its outer bounds are preserved by differentiation, and they are also
preserved by taking powers, proving the claim.
\end{proof}

The zero-counting measure of every descendant, and hence $\rho_P$, is invariant under
multiplication by
$\omega_q$: the polynomial changes by a scalar factor under this rotation, and taking powers
and derivatives preserves that property.

\subsection{Sampling on a signed radial star}
Let
\[
 \eta_q:=(z\mapsto z^q)_*\rootmu_P,
\]
a probability measure on $[b_-,b_+]$. Let $Y_1,Y_2$ be independent with law $\eta_q$, and
let $\mathsf U_q$ be uniform on $[0,1]$ and independent of them. Set
\[
 V_q=(1-\mathsf U_q)Y_1+\mathsf U_qY_2.
\]
Let $\mathsf R_q$ and $J_q$ be mutually independent and independent of these variables, with
\begin{equation}\label{eq:Tq}
 \Law(\mathsf R_q)=\frac1q\delta_1+\frac{q-1}{q}\Unif[0,1],
 \qquad J_q\text{ uniform on }\{0,\ldots,q-1\}.
\end{equation}
Choose the following explicit $q$th root of a real number:
\[
 \chi_q(v)=
 \begin{cases}
 v^{1/q},&v\ge0,\\
 e^{\pi i/q}|v|^{1/q},&v<0.
 \end{cases}
\]
In particular, $\chi_q(v)^q=v$.

\begin{theorem}[Signed radial sampling formula]\label{thm:radial}
For a polynomial \eqref{eq:radial-class} with $H$ nonconstant,
\begin{equation}\label{eq:radial-law}
 \rho_P=\Law\bigl(\omega_q^{J_q}\mathsf R_q\chi_q(V_q)\bigr).
\end{equation}
Moreover,
\[
 \AllSh_P=\Sh_P=\mathscr R_{q,R_+,R_-},\qquad
 \mathcal Z_{P,n}\longrightarrow\mathscr R_{q,R_+,R_-}
 \quad\text{in Hausdorff distance}.
\]
All zeros of all descendants are already contained in this star.
\end{theorem}
\begin{proof}
Write $M_{q,k}=\int y^k\dd\eta_q(y)$. Rotational invariance gives
$\mom_k(P)=0$ unless $q$ divides $k$, and $\mom_{qk}(P)=M_{q,k}$. The moment formula
\eqref{eq:moments-intro} gives
\begin{equation}\label{eq:radial-moments}
 \int z^{qk}\dd\rho_P(z)=\frac1{qk+1}\sum_{j=0}^kM_{q,j}M_{q,k-j},
\end{equation}
while all other moments vanish. Direct integration of the real chord and of $\mathsf R_q$ gives
\[
 \mathbb E[V_q^k]=\frac1{k+1}\sum_{j=0}^kM_{q,j}M_{q,k-j},
 \qquad
 \mathbb E[\mathsf R_q^{qk}]=\frac{k+1}{qk+1}.
\]
The law in \eqref{eq:radial-law} therefore has the same holomorphic moments and rotational
invariance as $\rho_P$.

Both measures are supported on \eqref{eq:signed-star}. Their pushforwards by $z\mapsto z^q$
are measures on the compact real interval $[b_-,b_+]$ with identical moments, hence
coincide. A rotation-invariant measure on this star is determined by its pushforward: away
from zero it divides mass equally among the $q$ preimages, and zero has a single preimage.
This proves the sampling formula.

By Proposition~\ref{prop:global-facts}, the core $\Sh_P$ is connected and contains the
roots of $P$, hence the tip of every spoke of positive length. A connected subset of the
signed star containing these tips contains
every spoke: the signed radius, positive on one spoke and negative on the others,
is continuous on the star, so the intermediate value theorem forces the center
and every intermediate radius. Thus $\Sh_P$ is the whole star. Lemma~\ref{lem:star}
and the convergence $\mathcal Z_{P,n}\to\Sh_P$ in \eqref{eq:core-union} give the equality
with the full shadow, the inclusion of all zeros of descendants, and Hausdorff convergence.
\end{proof}

For example, $(z^q-a^q)(z^q-b^q)$ with $0<a<b$ has $q$ spokes of length $b$; the inner ring
of roots changes the measure but not the support. In contrast, $(z^q-a^q)(z^q+b^q)$ with $a,b>0$ has
two interlaced families of spokes of lengths $a$ and $b$. When $q=2$ this gives the unequal
cross with endpoints $\pm a,\pm ib$. Arbitrary multiplicities are allowed, and affine covariance gives
the corresponding formulas after any affine change of coordinates.

\begin{corollary}[Uniform spokes for a regular polygon]\label{cor:regular-polygon}
For $P(z)=z^q-1$, $q\ge2$,
\begin{equation}\label{eq:uniform-spokes}
 \rho_P=\frac1{q^2}\sum_{j=0}^{q-1}\delta_{\omega_q^j}
 +\frac{q-1}{q^2}\sum_{j=0}^{q-1}
       (r\mapsto\omega_q^jr)_*(\dd r\!\restriction[0,1]).
\end{equation}
Thus the nonatomic part has constant arclength density on each spoke, its total mass is
$1-1/q$, and the mass at each root is $1/q^2$.
\end{corollary}
\begin{proof}
Here $\eta_q=\delta_1$, so $V_q=1$. Substitute \eqref{eq:Tq} into \eqref{eq:radial-law}
and distribute the mass equally over the rotations. Equivalently, the nonzero moments are
$\int z^{qk}\dd\rho_P=(k+1)/(qk+1)$, which are exactly the moments of \eqref{eq:uniform-spokes}.
\end{proof}

\begin{figure}[t]
\centering
\includegraphics{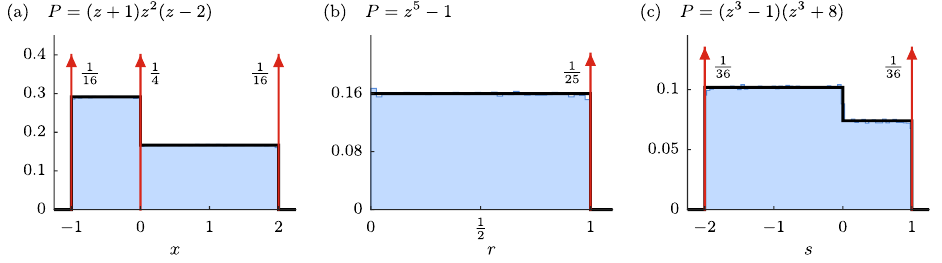}
\caption{Histograms of $\bar\nu_{P,n}$ (shaded) and the limiting
densities (lines); arrows: atoms, labeled by their masses.
(a)~Density on $\R$ for $P=(z+1)z^2(z-2)$, $n=100$, from
\eqref{eq:collinear}.
(b)~Arclength density on each spoke for $P=z^5-1$, $n=80$, from
\eqref{eq:uniform-spokes}.
(c)~Arclength density on each spoke for $P=(z^3-1)(z^3+8)$, $n=60$,
against the signed radius $s$, positive on the spokes through the cube
roots of unity; \eqref{eq:radial-law} gives the values $\tfrac2{27}$
and $\tfrac{11}{108}$.}
\label{fig:densities}
\end{figure}

Figure~\ref{fig:densities} compares \eqref{eq:collinear},
\eqref{eq:uniform-spokes} and \eqref{eq:radial-law} with the pooled
zero measures at finite power levels.

\begin{example}[A spurious endpoint tie]
For $P(w)=w^3-1$, the formal saddle equation at $\alpha=0$ is $3w^2(w-z)=0$. Counting the
double root $w=0$ twice and taking the two largest limiting real heights would declare a tie
whenever $|z^3-1|<1$. This includes an open set in the root triangle.
By Lemma~\ref{lem:star}, however, every zero of every descendant lies on the three
spokes of Corollary~\ref{cor:regular-polygon}. At the endpoint $\alpha=0$, the formal equation
no longer represents two
distinct contributing saddles. The barrier is instead defined at
interior ratios, and its endpoint behavior is determined by a limit.
\end{example}
\begin{proof}[Verification]
The degree-normalized endpoint heights are $0,0,\frac13\log|z^3-1|$, with $-\infty$
at an original root. Off these roots, a saddle tending to $z\ne0$ satisfies
$w^2(w-z)=\frac\alpha3(w^3-1)$, so $\alpha\log|w-z|\to0$; a saddle tending to zero
has height limit zero by continuity. The point $z_0=\frac14+\frac i8$ lies strictly
inside the root triangle, satisfies $|z_0^3-1|<1$, and has $\Imm(z_0^3)\ne0$.
These strict conditions persist on an open neighborhood, which therefore misses all
three spokes while retaining the false maximal tie.
\end{proof}

\subsection{The classification}\label{sec:dimension}
\begin{theorem}[Complete zero-area classification]\label{thm:dimension}
Let $P$ be a nonconstant polynomial with at least two distinct roots. The following
conditions are equivalent:
\begin{enumerate}[label=\textup{(\roman*)},leftmargin=2.2em]
\item $\Area(\Sh_P)=0$.
\item $\Area(\AllSh_P)=0$.
\item Every root of $P$ lies on a reflection axis of its root multiset.
\item There are $a,b,c\in\C$ with $a c\ne0$, an integer $q\ge1$, an integer $r_0\ge0$, and a
nonconstant polynomial $H$ with only nonzero real roots, such that
\begin{equation}\label{eq:dimension-normalform}
 cP(az+b)=z^{r_0}H(z^q).
\end{equation}
\end{enumerate}
In these cases $\AllSh_P=\Sh_P$: the common set is a segment if $q=1$, and the finite star
described in Theorem~\ref{thm:radial} if $q\ge2$, after the corresponding affine change of
coordinates. In every other case both shadows have positive area and the same nonempty
interior. A polynomial with a single distinct root has a singleton shadow.
\end{theorem}
\begin{proof}
The equivalence of (i) and (ii) follows from Proposition~\ref{prop:global-facts}, because
$T_P$ is countable. If (i) holds, the planar
part of $\rho_P$ has zero mass. By Theorem~\ref{thm:axes}, the support is contained in the
union of $\mathscr A_P$ and the finite root set. A root outside $\mathscr A_P$ would be
isolated in the support, contradicting Proposition~\ref{prop:global-facts}. This proves (iii).

Suppose (iii). Center at the root centroid and choose one of the reflection axes as the real
line; an affine change of variable and normalization then give a monic
real polynomial $Q$ of degree $d$. Write $[z^k]Q$ for the coefficient of $z^k$ in $Q$.
Let $q=\gcd\{d-k:k<d,\,[z^k]Q\ne0\}$ and let $r_0$ be the least exponent $k$ with $[z^k]Q\ne0$.
The set of gaps $d-k$ is nonempty, since $Q$ has two distinct roots, so $q\ge1$.
For $\lambda\ne0$, rotation by $\lambda$ preserves the root multiset exactly when
$Q(\lambda z)=\lambda^dQ(z)$; coefficient comparison makes this equivalent to
$\lambda^q=1$. Every reflection axis passes through zero. If a reflection fixes a
root $z$, write $z=\lambda\overline z$. Composing with conjugation gives a rotation
preserving the multiset, so $\lambda^q=1$ and $z^q\in\R$.
Every exponent with a nonzero coefficient is congruent to $r_0$ modulo $q$, hence
$Q(z)=z^{r_0}H(z^q)$ for a monic nonconstant $H$ with $H(0)\ne0$.
If $H(t)=0$, choose $z$ with $z^q=t$; then $Q(z)=0$, so $t$ is real and nonzero. This proves (iv), including $q=1$.

Finally, for $q=1$, repeated application of Rolle's theorem confines every zero of every
descendant to the root
segment. For $q\ge2$, Lemma~\ref{lem:star} confines every zero of every descendant to a
radial star; it allows both positive and negative roots of $H$. Thus (iv) implies
(ii). Theorems~\ref{thm:collinear} and~\ref{thm:radial} identify the support with the whole
indicated segment or star, which also proves equality of the full shadow and the core. In
the remaining cases the area is positive, and definability gives nonempty interior. The
one-root assertion follows directly by differentiation.
\end{proof}

This criterion involves only the root multiset. It can be
checked from the roots or from the normal form
\eqref{eq:dimension-normalform}, without selecting saddles.
Reflection symmetry alone is not sufficient: every root must lie on
an axis. When some roots lie off the axes, planar and curve mass may
coexist.

\subsection{Finite coefficient tests}
\begin{proposition}[Reflection and zero area from coefficients]\label{prop:coefficient-test}
Center $P$ at its root centroid and write
\[
 Q(w):=P(\centroid_P+w)
       =w^d+\sum_{j=2}^d\beta_j w^{d-j},
 \qquad \mathcal J:=\{j:\beta_j\ne0\}.
\]
If $\mathcal J=\varnothing$, both shadows are the singleton $\{\centroid_P\}$.
Otherwise choose any $j_0\in\mathcal J$. All reflection axes are obtained from the at most
$j_0\le
d$ solutions of
\begin{equation}\label{eq:coefficient-reflection}
 \lambda^{j_0}=\frac{\beta_{j_0}}{\overline{\beta_{j_0}}}
 \quad\text{that also satisfy}\quad
 \beta_j=\overline{\beta_j}\lambda^j\quad(j\in\mathcal J).
\end{equation}
Every candidate has modulus one. A retained $\lambda$ corresponds to the reflection
$w\mapsto\lambda\overline w$ in centered coordinates, so its axis in the observation plane
passes through $\centroid_P$.

If no candidate is retained, both shadows have positive area. Otherwise put
$q:=\gcd\mathcal J$, choose
a retained $\lambda=e^{2i\theta}$, and set $r_0:=d-\max\mathcal J$. Then
\begin{equation}\label{eq:coefficient-real-form}
 e^{-id\theta}Q(e^{i\theta}x)=x^{r_0}H(x^q),
\end{equation}
where $H$ is a monic real polynomial with $H(0)\ne0$. Both shadows have zero area if and
only if $H$ is real-rooted. Thus the zero-area classification requires only centering,
finitely many coefficient identities, an integer greatest common divisor, and a
real-rootedness test; no saddle-selection test is needed.
\end{proposition}
\begin{proof}
Reflection in $e^{i\theta}\R$ preserves the root multiset of $Q$ if and only if the monic
rotated polynomial $e^{-id\theta}Q(e^{i\theta}x)$ has real coefficients. Its coefficient
indexed by $j$ is $\beta_j e^{-ij\theta}$, whose reality is exactly
$\beta_j=\overline{\beta_j}e^{2ij\theta}$. This proves \eqref{eq:coefficient-reflection},
including completeness of the finite list.

Fix one solution $\lambda_0$. The other solutions are precisely
$\lambda_0\omega$ with $\omega^j=1$ for all $j\in\mathcal J$, or equivalently $\omega^q=1$. In
coordinates rotated by $\theta$, the axes are therefore the $q$ lines of angles $\pi\ell/q$,
$0\le\ell<q$. All exponents of the rotated polynomial are congruent to $d$ modulo $q$, and
their smallest is $r_0$. This proves \eqref{eq:coefficient-real-form} with $H$ real and
$H(0)\ne0$. A nonzero point lies on one of the $q$ axes exactly when its $q$th power is
real. Hence all roots lie on the axes exactly when all roots of $H$ are real.
Theorem~\ref{thm:dimension} completes the proof.
\end{proof}

For algebraic coefficients these tests can be performed exactly. To test real-rootedness
even when $H$ has multiple roots, take its square-free part $G:=H/\gcd(H,H')$. Form $G_0=G$,
$G_1=G'$, and successively $G_{k+1}=-\operatorname{rem}(G_{k-1},G_k)$ stopping as soon as
the last polynomial is constant, where $\operatorname{rem}$ is the remainder in polynomial
division. Sturm's theorem states that the difference between the numbers of sign variations
of this sequence at $-\infty$ and $+\infty$ counts the distinct real roots. It equals $\deg
G$ exactly when $H$ is real-rooted; see \cite{BCR}.

\section{Analytic regularity of the planar density}\label{sec:planar-stability}
The decomposition of the shadow measure can be sharpened away from an
exceptional set consisting of finitely many arcs and points. Over its
complement, the sheets of the incidence space project locally invertibly, so
the densities that they contribute are real-analytic and strictly positive.

\begin{theorem}[Regularity of the planar measure]\label{thm:planar-regularity}
For every polynomial $P$ there is a compact definable set
$\mathscr E_P\subseteq\Sh_P$, consisting of finitely many points and
bounded $C^1$ arcs of finite length, such that
$\partial\Sh_P\subseteq\mathscr E_P$ and
\[
 \rho_P\!\restriction(\C\setminus\mathscr E_P)
   =f_P\,\dd A\!\restriction(\C\setminus\mathscr E_P).
\]
On each component of $\operatorname{int}\Sh_P\setminus\mathscr E_P$,
the density $f_P$ has a strictly positive real-analytic representative.
In particular
\begin{equation}\label{eq:positive-planar-density}
 f_P(z)>0\quad\text{for almost every }z\in\operatorname{int}\Sh_P,
 \qquad
 f_P(z)=0\quad\text{for almost every }z\notin\operatorname{int}\Sh_P.
\end{equation}
Thus the planar part of $\rho_P$ and Lebesgue measure restricted to
$\operatorname{int}\Sh_P$ have exactly the same null sets.
\end{theorem}
\begin{proof}
Use the compact definable incidence space $\Inc_P$ and a
$C^2$ stratification of it satisfying the frontier condition, refined so that the spatial
projection has constant rank on each stratum. Away from the root lines, the
endpoint fibers and the ordinary saddle-collision image
$\mathscr G_P^{\mathrm{col}}$, a two-dimensional rank-two stratum is locally
a $C^2$ graph $\alpha=\alpha(z)$. By the frontier condition, the incidence space coincides
with this graph near each of its points. The graph therefore divides a small tube about it
into two connected
chambers that do not meet $\Inc_P$. On each chamber the selected saddle is unique and
simple, and its critical height is one fixed analytic branch;
by continuity, the maximum of these two branches equals the barrier throughout
the tube. The two branches are different, because the graph belongs to the support.
Their difference $g$ satisfies $g(z,\alpha(z))=0$ and
$\nabla_zg\ne0$ by \eqref{eq:g-z}. Differentiating the graph identity gives
$\partial_\alpha g\ne0$, so the analytic implicit-function theorem makes the
graph real-analytic. Additional heights may be tied on it without changing this pair.
Let $\mathscr E_P$ be the closure of the union of $\partial\Sh_P$, the
finite set $\mathscr D_P=Z(P)\cup Z(P')\cup\{\centroid_P\}$,
$\mathscr G_P^{\mathrm{col}}\cap\Sh_P$, and the spatial images of all
lower-dimensional or rank-deficient strata. Include the root lines and endpoint
fibers among the exceptional strata.

The image of each rank-deficient stratum has dimension at most one,
as does the image of a stratum of dimension at most one. Closure preserves
the dimension of definable sets, and the incidence space is compact.
Consequently $\mathscr E_P$ is compact and definable of dimension at
most one. The monotone-arc argument of Section~\ref{subsec:tame} gives
its stated decomposition and finite length. It contains every
possible curve or atomic contribution to $\rho_P$.

Fix $z\in\Sh_P\setminus\mathscr E_P$. The incidence fiber over
$z$ is nonempty and compact, and all its points have projection rank
two. Each such point
is isolated in the fiber by the inverse-function theorem, so the
fiber is finite: an infinite fiber would have an accumulation point
in that same compact fiber, contradicting local invertibility.
There is a disk $D$ about $z$ over which all these points continue as
finitely many real-analytic graphs $\alpha=\alpha_j(\zeta)$.
There are no additional preimages after shrinking $D$. Otherwise
compactness would give a limiting preimage over $z$ outside the
chosen neighborhoods. In particular $D\subseteq\Sh_P$.

For each graph choose the pair of branches that switch across it, as above, and write
$g_j$ for
the corresponding difference of critical heights. Continuous cutoffs with disjoint compact
supports,
equal to one near these preimages, partition every nearby incidence fiber.
Apply \eqref{eq:density} to each cutoff and add, counting each geometric
preimage once. On $D$ this gives
\begin{equation}\label{eq:analytic-planar-density}
 f_P(\zeta)=\frac1{\pi d^2}
    \sum_j
    \frac{|\nabla_zg_j(\zeta,\alpha_j(\zeta))|^2}
         {|\partial_\alpha g_j(\zeta,\alpha_j(\zeta))|}.
\end{equation}
Every numerator is positive and every denominator is nonzero. The
sign of each denominator is locally constant before taking its
absolute value, so every summand is real-analytic. This proves local
analyticity and strict positivity. The local representatives agree
on overlaps, since they represent the same measure, proving the
componentwise assertion. Finally $\mathscr E_P$ and the boundary have
zero planar area, and the measure is supported on $\Sh_P$.
This proves \eqref{eq:positive-planar-density} and the null-set claim.
\end{proof}

Figure~\ref{fig:candidates}(b) approximates the planar density of a
cubic near its centroid. The caustic arcs inside the core lie in the set
$\mathscr E_P$ constructed in the proof.

\section{Finite atomic laws and vanishing root masses}
\label{sec:weighted}
The shadow measure and the core of $P$ depend only on the normalized root
measure $\rootmu_P$, and their construction extends to arbitrary finite
atomic probability measures \cite{HS2}. In this setting the root masses
vary continuously and may, in particular, tend to zero. We recall the
framework and extend the conservation law, the decomposition of the shadow
measure, and reflection rigidity to it.

\subsection{Weighted potentials and supports}
Let
\begin{equation}\label{eq:weighted-data}
 \eta=\sum_{a\in A}p_a\delta_a,\qquad
 A\subset\C\text{ finite},\quad p_a>0,\quad \sum_{a\in A}p_a=1.
\end{equation}
We call such a measure a \emph{finite atomic law} and the points of $A$ its
\emph{sites}. The masses $p_a$ are arbitrary positive real numbers: unlike
the weights $p_j=r_j/d$ of a polynomial, they need not be rational. Put
\[
 \begin{gathered}
 C_\eta(w)=\sum_{a\in A}\frac{p_a}{w-a},\qquad
 h_\eta(w)=\frac1{C_\eta(w)},\qquad
 \centroid_\eta=\sum_{a\in A}p_aa,\\
 F_{t,\eta}(w)=w-t h_\eta(w),\qquad 0<t<1.
 \end{gathered}
\]
For $z\notin A$, the weighted landscape and its local holomorphic phase are
\[
 L_{\eta,z,t}(w)=U_\eta(w)-t\log|w-z|,\qquad
 \Phi_{\eta,z,t}(w)=\sum_{a\in A}p_a\log(w-a)-t\log(w-z).
\]
Define
\[
 E_\ell(\eta,z,t)=
 \{w\in\C\setminus(A\cup\{z\}):L_{\eta,z,t}(w)\ge\ell\}
       \cup\{z,\infty\},
\]
and let $\Blev_\eta(z,t)$ be the supremum of the levels $\ell$ for which the two
marked ends belong to the same connected component of this set. At the
barrier level, a saddle is selected precisely when it belongs to the
closures of both components of the strict superlevel set that contain the marked ends.
The reduced saddle equation is
\begin{equation}\label{eq:weighted-saddle}
 (w-z)\RedN_\eta(w)-t\RedQ_\eta(w)=0,
\end{equation}
where
\[
 \RedQ_\eta(w)=\prod_{a\in A}(w-a),\qquad
 \RedN_\eta(w)=\sum_{a\in A}p_a\prod_{b\in A\setminus\{a\}}(w-b).
\]
It has degree $|A|$ and leading coefficient $1-t$, and
$C_\eta=\RedN_\eta/\RedQ_\eta$.

\begin{proposition}[Weighted measures and selection]
\label{prop:weighted-facts}
\cite[Propositions~2.1 and~2.5]{HS2}
There is a unique probability measure $\mu_{t,\eta}$ with potential
\begin{equation}\label{eq:weighted-potential}
 U_{\mu_{t,\eta}}(z)=
 \frac{\Blev_\eta(z,t)-\mathfrak h_1(t)}{1-t},\qquad
 \mathfrak h_1(t)=-t\log t-(1-t)\log(1-t),
\end{equation}
with subharmonic extension at the atoms. Set
\begin{equation}\label{eq:weighted-supports}
 S_{t,\eta}=\supp\mu_{t,\eta},\qquad
 \rho_\eta=2\int_0^1(1-t)\mu_{t,\eta}\dd t,\qquad
 S_\eta=\supp\rho_\eta=
       \overline{\bigcup_{0<t<1}S_{t,\eta}}.
\end{equation}
All these measures and sets are supported in $\conv A$. Their atomic
parts are
\begin{equation}\label{eq:weighted-atoms}
 \mu_{t,\eta}^{\mathrm{at}}=
       \sum_{a\in A}\frac{(p_a-t)_+}{1-t}\delta_a,
 \qquad \rho_\eta^{\mathrm{at}}=\sum_{a\in A}p_a^2\delta_a.
\end{equation}
A site $a$ belongs to $S_{t,\eta}$ whenever $0<t\le p_a$
and $t<1$. The measures and supports depend continuously on the ratio
$t\in(0,1)$, and
\begin{align}\label{eq:weighted-endpoints}
 \mu_{t,\eta}&\Longrightarrow\eta,
 &S_{t,\eta}&\longrightarrow A\cup Z(C_\eta)
        &&(t\downarrow0),\notag\\
 \mu_{t,\eta}&\Longrightarrow\delta_{\centroid_\eta},
 &S_{t,\eta}&\longrightarrow\{\centroid_\eta\}
        &&(t\uparrow1),
\end{align}
where $Z(C_\eta)$ is the set of finite zeros of the rational function
$C_\eta$, and the set limits are in Hausdorff distance.

Off $S_{t,\eta}$, the Cauchy transform $C_{\mu_{t,\eta}}$ is nonzero and
\begin{equation}\label{eq:weighted-inverse}
 \begin{split}
 W_{t,\eta}(z)&=z+\frac{t}{(1-t)C_{\mu_{t,\eta}}(z)},\\
 F_{t,\eta}(W_{t,\eta}(z))&=z,\qquad
 C_\eta(W_{t,\eta}(z))=(1-t)C_{\mu_{t,\eta}}(z).
 \end{split}
\end{equation}
At observation points $z\notin A$ this is the selected simple saddle;
$F_{t,\eta}'(W_{t,\eta})\ne0$. It is holomorphic in $z$ and locally
real-analytic in $t$. The separating-contour, local-maximum, and harmonic
selection assertions of Proposition~\ref{prop:selection-facts} hold with
$L_{\eta,z,t}$ and $\Blev_\eta$ in place of their polynomial counterparts,
jointly under changes of the positive masses and distinct locations.
\end{proposition}

When $\eta=\rootmu_P$, the normalizations agree through
\begin{equation}\label{eq:polynomial-weighted-bridge}
 \mu_{t,\eta}=\mu_{dt,P},\quad
 S_{t,\eta}=\FixSh_{dt,P},\quad
 \rho_\eta=\rho_P,\quad S_\eta=\Sh_P,\quad
 \Blev_\eta(z,t)=d^{-1}\Blev_P(z,dt).
\end{equation}
Thus $S_\eta$ is the counterpart of the core $\Sh_P$, and not of the full
shadow $\AllSh_P$; we call it the \emph{core} of $\eta$, and we call
$W_{t,\eta}$ the selected inverse.

\begin{proposition}[Definability, connectedness, and uniform complexity]
\label{prop:weighted-tame}
\cite[Theorems~4.1--4.2, Lemma~4.7, and Corollaries~4.4, 4.8--4.9]{HS2}
For a fixed bound $s$ on $|A|$, the weighted barriers, fixed-ratio
supports, and cores form definable families in the sites and the
positive masses. The incidence set
\[
 \Inc_\eta=\overline{\{(z,t):0<t<1,\ z\in S_{t,\eta}\}}
\]
is compact and has dimension at most two, with spatial projection $S_\eta$.
The set $S_\eta$ is connected. Each point of a fixed-ratio support can
be joined to $\centroid_\eta$ by a continuous definable path through
supports at increasing ratios.

There is an integer $M_s$, depending only on $s$, such that
$\partial S_\eta$ is a union of at most $M_s$ points and at most $M_s$
$C^1$ arcs on each of which both coordinates are monotone. In particular,
if $A\subset\{|z|\le R\}$, then
\[
 \Haus^1(\partial S_\eta)\le4M_sR.
\]
The numbers of components and the Betti numbers of the fibers are uniformly
bounded, and there are finitely many homeomorphism types.
\end{proposition}

\subsection{Finite-location limits and selected transport}
Let $B\supset A$ be another finite set. A point of $B\setminus A$
represents a root location whose limiting mass is zero. For such a point
$b$ with $C_\eta(b)\ne0$, define its \emph{transport} at ratio $t$ by
\[
 \theta_b(t)=b-\frac{t}{C_\eta(b)}.
\]
Only selected transports contribute to the Hausdorff limit of the zero sets. Set
\begin{equation}\label{eq:selected-transports}
 \begin{split}
 \mathcal O_t(\eta,B)=\bigl\{\theta_b(t):
  &b\in B\setminus A,\ C_\eta(b)\ne0,\\[-2mm]
  &\theta_b(t)\notin S_{t,\eta},\quad
       W_{t,\eta}(\theta_b(t))=b\bigr\},
 \end{split}
\end{equation}
and define the \emph{finite-location limit set}
\begin{equation}\label{eq:transport-limit-set}
 K_{\eta,B}:=S_\eta\cup B\cup
             \overline{\bigcup_{0<t<1}\mathcal O_t(\eta,B)}.
\end{equation}
This definition retains both the limiting masses and the locations at
which mass disappears.

\begin{proposition}[Finite-location limits of zero sets]
\label{prop:finite-location}
\cite[Theorem~3.1, Corollary~3.5, and Proposition~4.5]{HS2}
Let $G_j$ be monic of degree $N_j\to\infty$, with
\begin{equation}\label{eq:finite-location-input}
 \eta_j=N_j^{-1}\divz G_j\Longrightarrow\eta,
 \qquad Z(G_j)\longrightarrow B
       \quad\text{in Hausdorff distance},
\end{equation}
where $B$ is finite. If $m_j/N_j\to t\in(0,1)$, then
\[
 \frac{\divz G_j^{(m_j)}}{N_j-m_j}\Longrightarrow\mu_{t,\eta},
 \qquad
 Z(G_j^{(m_j)})\longrightarrow S_{t,\eta}\cup\mathcal O_t(\eta,B).
\]
The pooled measures converge to $\rho_\eta$, and
\[
 \bigcup_{m=0}^{N_j-1}Z(G_j^{(m)})\longrightarrow K_{\eta,B}.
\]
More generally, if finite atomic laws $\eta_j\Longrightarrow\eta$ and
finite sets $B_j\supset\supp\eta_j$ satisfy $B_j\to B$ with $B$ finite,
then $K_{\eta_j,B_j}\to K_{\eta,B}$ in Hausdorff distance.

The set $K_{\eta,B}$ is obtained from $S_\eta$ by adjoining finitely
many closed straight segments, each meeting $S_\eta$ at an endpoint.
It is contained in
\[
 S_\eta\cup B\cup
 \bigcup_{\substack{b\in B\setminus A\\C_\eta(b)\ne0}}
       \{\theta_b(t):0\le t\le1\}.
\]
In particular, $K_{\eta,B}$ is connected and has the same interior and
area as $S_\eta$. When $B=A$, it equals $S_\eta$.
\end{proposition}

\begin{corollary}[Polynomial masses and finite-location limits]
\label{cor:polynomial-transport-semialgebraic}
If $\eta=\rootmu_P$ and $B\supset Z(P)$ is finite, then $K_{\eta,B}$
is semialgebraic.
\end{corollary}
\begin{proof}
By Proposition~\ref{prop:finite-location}, it is the union of
$S_\eta=\Sh_P$ and finitely many closed straight segments. Apply
Theorem~\ref{thm:semialgebraic}.
\end{proof}

\subsection{Conservation and regularity for weighted laws}
\begin{proposition}[Weighted conservation]\label{prop:weighted-conservation}
For every finite atomic law \eqref{eq:weighted-data},
\begin{equation}\label{eq:weighted-conservation}
 C_{\rho_\eta}'=-C_\eta^2\qquad\text{on }\C\setminus S_\eta.
\end{equation}
For every closed piecewise smooth curve $\gamma$ in this complement,
\[
 \wind(\gamma,a)=\wind(\gamma,\centroid_\eta)\qquad(a\in A).
\]
\end{proposition}
\begin{proof}
On a disk $D$ whose closure is disjoint from $S_\eta$, put
$u_t=(1-t)C_{\mu_{t,\eta}}$. Differentiating
\eqref{eq:weighted-inverse} gives
\[
 u_t\,\partial_tu_t=\partial_zu_t,\qquad u_0=C_\eta,\quad u_1=0.
\]
Since $C_{\rho_\eta}=2\int_0^1u_t\dd t$, integration yields
\eqref{eq:weighted-conservation}. The bound
$|\partial_zu_t|\le(1-t)/\dist(D,S_\eta)^2$ on $D$
justifies differentiation under the integral.

For $a\in A$, the curve $W_{t,\eta}\circ\gamma$ avoids $a$.
The endpoint limits in \eqref{eq:weighted-endpoints} imply, uniformly on
$\gamma$,
\[
 W_{t,\eta}(z)\longrightarrow z\quad(t\downarrow0),\qquad
 \frac{1-t}{t}(W_{t,\eta}(z)-a)
       \longrightarrow z-\centroid_\eta\quad(t\uparrow1).
\]
The winding number is constant in $t$, giving the assertion.
\end{proof}

\begin{proposition}[Weighted measure decomposition and planar regularity]
\label{prop:weighted-regularity}
The measure $\rho_\eta$ has a decomposition
\[
 \rho_\eta=f_\eta\dd A+
       \vartheta_\eta\,\Haus^1\!\restriction\Gamma_\eta+
       \sum_{a\in A}p_a^2\delta_a,
\]
with nonnegative integrable densities and a finite union $\Gamma_\eta$
of bounded $C^1$ arcs of finite length. There is a compact definable set
$\mathscr E_\eta\subset S_\eta$ of dimension at most one, containing
$\partial S_\eta$, such that off $\mathscr E_\eta$ the measure is
$f_\eta\dd A$. On each component of
$\operatorname{int}S_\eta\setminus\mathscr E_\eta$, the density is
strictly positive and real-analytic. In particular, $f_\eta\dd A$ and
Lebesgue measure on $\operatorname{int}S_\eta$ have the same null sets.
\end{proposition}
\begin{proof}
As in the proof of Theorem~\ref{thm:decomposition}, all solutions of the
saddle equation, including multiple ones, are locally confined to a compact
set, and the weighted first-variation formula makes the barrier locally Lipschitz
in $(z,t)$ off $A$, for $0<t<1$. Circle means of the potentials and
averaging over $t$ therefore give
$\rho_\eta\!\restriction(\C\setminus A)\ll\Haus^1$.
This argument remains valid at saddle collisions.

Stratify the compact incidence set $\Inc_\eta$ by a $C^2$
stratification satisfying the frontier condition and refined by the rank of
the spatial projection, as in the proof of
Theorem~\ref{thm:planar-regularity}. Away from collisions, each
two-dimensional rank-two stratum is locally a smooth graph $t=t(z)$.
On the two chambers into which it divides a small tube, the selected simple
saddles have two fixed critical heights, and by
continuity their maximum equals the barrier throughout
the tube. Their difference $g$ has $\nabla_zg\ne0$, since distinct
saddles have distinct spatial derivatives $t/(w-z)$.
Differentiating the graph identity gives $\partial_tg\ne0$, so the
analytic implicit-function theorem shows that the graph is real-analytic.
Additional heights may be tied on it.

Let $\mathscr E_\eta$ be the closure, inside $S_\eta$, of the union of the
boundary $\partial S_\eta$,
the sites, the zeros of $C_\eta$, the centroid, the ordinary collision image,
and the spatial images of lower-dimensional or rank-deficient strata;
include the lines above the sites and the endpoint fibers among the
exceptional strata.
Definable stratification and finite-length curve decomposition make
$\mathscr E_\eta$ compact, of dimension at most one and zero area,
and the union of finitely many points and bounded regular $C^1$ arcs
of finite total length. Call the union of these arcs $\Gamma_\eta$.

Above any point of $S_\eta$ outside $\mathscr E_\eta$, compactness and local
invertibility show that the preimages form finitely many analytic
graphs $t=t_j(z)$, with no further preimages over a smaller disk.
Continuous cutoffs with disjoint compact supports, equal to one near these
preimages,
partition every nearby incidence fiber. For each cutoff, the fixed-ratio
measure on the interface is
$|\nabla_zg_j|\,\Haus^1/(2\pi(1-t))$, as in \eqref{eq:line-density}.
Multiplication by $2(1-t)$ and the local area formula thus give
\begin{equation}\label{eq:weighted-analytic-density}
 f_\eta(z)=\frac1\pi\sum_j
       \frac{|\nabla_zg_j(z,t_j(z))|^2}
            {|\partial_tg_j(z,t_j(z))|}.
\end{equation}
The local measure laws agree on overlaps and give the global planar
density, set to zero on $\mathscr E_\eta$ and outside $S_\eta$.
Each summand is positive and real-analytic.

The remaining measure is carried by $\mathscr E_\eta$.
After removing the sites, the absolute continuity with respect to arclength
proved above and
$\Haus^1(\mathscr E_\eta)<\infty$ give an integrable nonnegative
Radon--Nikodym density. Finitely many points have zero arclength, so this is
$\vartheta_\eta\,\Haus^1\!\restriction\Gamma_\eta$.
Formula~\eqref{eq:weighted-atoms} supplies precisely the remaining
$p_a^2\delta_a$ terms. Finiteness of $\rho_\eta$ makes both densities
integrable. The exceptional set has zero area and contains the boundary;
strict positivity of \eqref{eq:weighted-analytic-density} on the
remaining interior proves the asserted equality of null sets.
\end{proof}

\subsection{Weighted reflection and zero area}

The reflection argument also applies to arbitrary positive masses.
A \emph{reflection axis} of $\eta$ is a line
whose reflection preserves $\eta$, including the mass at each
site.

\begin{theorem}[Reflection and dimension for every finite atomic law]
\label{thm:weighted-dimension}
Let $\eta=\sum_{a\in A}p_a\delta_a$ be a probability measure with
finite support, $p_a>0$, and at least two support points.  The curve part of $\rho_\eta$ is
supported on reflection axes of $\eta$, and
\begin{equation}\label{eq:weighted-dimension}
 \Area(S_\eta)=0
 \quad\Longleftrightarrow\quad
 \text{every point of $A$ lies on a reflection axis of $\eta$}.
\end{equation}
In this case
\[
 S_\eta=\bigcup_{a\in A}[\centroid_\eta,a],
\]
where $\centroid_\eta=\int z\dd\eta(z)$ is the centroid of $\eta$; this set
is a segment or a radial star. If $\eta$ has no reflection axis, then
\[
 \rho_\eta=f_\eta\dd A+\sum_{a\in A}p_a^2\delta_a,
 \quad \int f_\eta\dd A=1-\sum_{a\in A}p_a^2,
 \quad S_\eta=\overline{\operatorname{int}S_\eta}.
\]
\end{theorem}
\begin{proof}
Proposition~\ref{prop:weighted-regularity} supplies the decomposition into
planar, curve, and atomic parts. On a regular arc carrying curve mass, the
rank-one stratum of the incidence set gives a stationary support arc over a ratio
interval. Apply the argument of Theorem~\ref{thm:stationary-reflection}
with $h_\eta=1/C_\eta$ in place of $h_P$: equality of the selected heights
forces the reflection in the tangent line to exchange two distinct saddle
branches. Their saddle
equations give the analogue of \eqref{eq:h-reflection-curve} for $h_\eta$ on
the infinite image of either
branch. Clearing denominators and comparing poles
therefore shows that this reflection preserves every mass of $\eta$.
There are only finitely many axes, so the curve part is carried
by their closed union. The conclusions in the absence of axes follow as in
Theorem~\ref{thm:axes}, using connectedness from
Proposition~\ref{prop:weighted-tame} and
$\int f_\eta\dd A=1-\sum_a p_a^2>0$.

If $S_\eta$ has zero area, the planar part of $\rho_\eta$ vanishes. Outside
the axes, the support $S_\eta$ could then contain only the finitely many sites.
A connected set with at least two points has no isolated points,
so the whole support lies on the
axes. Since $A\subset S_\eta$, this proves the forward implication
in \eqref{eq:weighted-dimension}.

Conversely, suppose every site lies on an axis. Approximate the masses by
integer multiplicities:
\[
 r_{a,j}=1+\lfloor jp_a\rfloor,\qquad
 P_j(w)=\prod_{a\in A}(w-a)^{r_{a,j}},\qquad
 \eta_j=\frac{1}{\sum_a r_{a,j}}\sum_a r_{a,j}\delta_a.
\]
The masses of $\eta_j=\rootmu_{P_j}$ converge to $p_a$. Equal masses receive
equal integers, so every reflection preserving $\eta$ preserves the entire
root multiset of every $P_j$. All roots of $P_j$ thus lie on reflection
axes of $P_j$. Theorem~\ref{thm:dimension} gives exactly
\[
 \Sh_{P_j}=\bigcup_{a\in A}[\centroid_{\eta_j},a],\qquad
 \centroid_{\eta_j}=\int z\dd\eta_j(z).
\]
For the segment case this uses $\centroid_{\eta_j}\in\conv A$; for
the signed star all tips of spokes of positive length are roots.
Now $\centroid_{\eta_j}\to\centroid_\eta$, and moving the center changes the
Hausdorff distance of these unions of segments by at most
$|\centroid_{\eta_j}-\centroid_\eta|$, as one sees by using the same
parameter on corresponding segments.
Proposition~\ref{prop:finite-location} also gives
$\Sh_{P_j}=S_{\eta_j}\to S_\eta$.
Uniqueness of the Hausdorff limit of compact sets identifies $S_\eta$ with
$\bigcup_{a\in A}[\centroid_\eta,a]$, a finite union of segments of zero area.
\end{proof}

\section{Area stability and convergence}\label{sec:area-stability}
The uniform bound on the complexity of the boundary $\partial S_\eta$ for
finite atomic laws controls the area
of small neighborhoods of the boundary. It yields quantitative area estimates
when the number of atoms is bounded. For qualitative convergence, a
persistence argument for regular interior points also allows the number
of root locations to grow.

\subsection{Boundary neighborhoods and area stability}
For a compact set $E\subset\C$ and $\eps\ge0$, define
$E^{[\eps]}=\{z:\dist(z,E)\le\eps\}$.
We write
$E\mathbin{\triangle}E'=(E\setminus E')\cup(E'\setminus E)$ for the symmetric
difference of two sets.

\begin{lemma}[A uniform boundary-neighborhood bound]\label{lem:shadow-tube}
There is a constant $C_s>0$, depending only on $s$, such that for any
atomic probability measure $\eta$ with at most $s$ atoms in a disk of
radius $R$, and every $\eps>0$,
\begin{equation}\label{eq:shadow-tube}
 \Area\bigl((\partial S_\eta)^{[\eps]}\bigr)
       \le C_s(R\eps+\eps^2).
\end{equation}
\end{lemma}
\begin{proof}
Use the at most $M_s$ coordinate-monotone arcs and $M_s$ points from
Proposition~\ref{prop:weighted-tame}. On each arc choose coordinate signs
so that $u(z)=\pm\Ree(z-c)\pm\Imm(z-c)$ increases along it, where $c$ is
the disk center. Then $|z-w|\le|u(z)-u(w)|$ for arc points and
$-2R\le u\le2R$. The occupied intervals of length $\eps$ in this
scalar range give at most $4R/\eps+1$ disks of radius $\eps$
covering the arc. Their closed versions also cover its closure, so the
closed $\eps$-neighborhood is covered by the same disks with radius
$2\eps$. Its area is at most
$16\pi R\eps+4\pi\eps^2$. Each point piece contributes
$\pi\eps^2$. Summing gives
$16\pi M_sR\eps+5\pi M_s\eps^2$, which proves the assertion
with $C_s=16\pi M_s$. Point pieces also cover a singleton hull.
\end{proof}

\begin{theorem}[Quantitative stability in planar area]\label{thm:area-stability}
Let $\eta,\widetilde\eta$ be probability measures with at most $s$
atoms each, all contained in one disk of radius $R$, and put
$\delta=d_{\mathrm H}(S_\eta,S_{\widetilde\eta})$.
For a constant $A_s$ depending only on $s$,
\begin{equation}\label{eq:area-symmetric-difference}
 \Area(S_\eta\mathbin{\triangle}S_{\widetilde\eta})
       \le A_s(R\delta+\delta^2).
\end{equation}
In particular the area is continuous under weak convergence and
Hausdorff convergence of finite atomic supports when the number of
atoms is bounded.

For polynomials $P,Q$ with at most $s$ distinct roots each, the same
bound holds with $\Area(\AllSh_P\mathbin{\triangle}\AllSh_Q)$ on the
left and $\delta=d_{\mathrm H}(\Sh_P,\Sh_Q)$ on the right.
The bound is independent of their degrees and multiplicities.
\end{theorem}
\begin{proof}
For nonempty compact sets $E,E'$ at Hausdorff distance $\delta>0$,
\[
 E\setminus E'\subseteq(\partial E')^{[\delta]},\qquad
 E'\setminus E\subseteq(\partial E)^{[\delta]}.
\]
Indeed a nearest point of $E'$ to a point outside $E'$ lies on
$\partial E'$. Apply Lemma~\ref{lem:shadow-tube} to $S_\eta$ and
$S_{\widetilde\eta}$ and
take $A_s=2C_s$. If $\delta=0$, the sets coincide and the assertion is
immediate. For the continuity statement, let $\eta_j\Longrightarrow\eta$ be
atomic probability measures with at most $s$ atoms, assume that the supports
themselves
converge to $\supp\eta$, and apply
Proposition~\ref{prop:finite-location} with
$B_j=\supp\eta_j$ and $B=\supp\eta$. There are then no additional
locations in the limiting data, so $K_{\eta,B}=S_\eta$. For polynomials, adding or
removing the countable sets $T_P,T_Q$ changes no planar area, by
Proposition~\ref{prop:global-facts}. This proves the last assertion.
\end{proof}

\begin{corollary}[H\"older area stability in coefficients]
\label{cor:holder-area}
For each fixed degree $d$ and each $R>0$, there are
constants $\widetilde C_R>0$ and $0<\widetilde\gamma_R\le1$ such that
\[
 \|\ind_{\AllSh_P}-\ind_{\AllSh_Q}\|_{L^1(\C)}
       \le\widetilde C_R\|P-Q\|_{\mathrm{coef}}^{\widetilde\gamma_R}
\]
whenever the monic degree-$d$ polynomials $P,Q$ have coefficient norms
at most $R$. The same holds for the cores. In particular,
$P\mapsto\Area(\AllSh_P)=\Area(\Sh_P)$ is locally H\"older continuous,
and the positive-area locus is open.
On every compact coefficient family containing only positive-area
shadows, their areas have a common positive lower bound.
\end{corollary}
\begin{proof}
The Cauchy root bound puts all these roots in the disk of radius
$1+R$. Combine Corollary~\ref{cor:holder-shape} with
Theorem~\ref{thm:area-stability}, applied with $s=d$ and this radius, using
$\delta\le2(1+R)$ to absorb the quadratic term into a constant times
$\delta$. The $L^1$ norm of the difference of two indicators is the
area of their symmetric difference. The last conclusions follow from
continuity and compactness.
\end{proof}

\subsection{Area convergence for varying root data}
\label{sec:unbounded-area}

The quantitative boundary estimate uses a bound on the number of
distinct roots. Qualitative area convergence needs no such bound when
the limiting location set is finite. The essential point is local:
a regular interior point persists under perturbation of the
logarithmic landscape, even when arbitrarily many new saddles appear
near the root clusters.

\begin{lemma}[Persistence of the regular interior]
\label{lem:regular-interior-persistence}
Let $\eta_j$ be finitely supported probability measures with
$\eta_j\Longrightarrow\eta$ and
$\supp\eta_j\to B$ in Hausdorff distance, where $B$ is finite.
There is a set $E\subset\C$ of planar area zero such that every
$z\in\operatorname{int}S_\eta\setminus E$ has a neighborhood
contained in $S_{\eta_j}$ for all sufficiently large $j$.
\end{lemma}
\begin{proof}
Write $A=\supp\eta\subset B$. Use
Proposition~\ref{prop:weighted-regularity}, and enlarge its exceptional
set by $B$ and the finitely many segments
\[
 \left\{b-\frac{t}{C_\eta(b)}:0\le t\le1\right\},\qquad
 b\in B,\quad C_\eta(b)\ne0.
\]
Their union $E$ has area zero. At a remaining interior point $z_0$
there is a regular switching sheet through a point $(z_0,t_0)$ whose
difference of critical heights $g$
satisfies $\nabla_zg\ne0$ and $\partial_tg\ne0$ at $(z_0,t_0)$.
Its two simple saddle branches are distinct. Neither of their values at
$(z_0,t_0)$ lies in $B$: a finite saddle $b\in B$ would give
$z_0=b-t_0/C_\eta(b)\in E$.

Choose two nearby ratios $t_-<t_+$ where opposite branches have
strictly larger height. Such ratios exist in every neighborhood of
$t_0$, since the nonzero derivative excludes a local extremum of $g$.
On either strict side the barrier equals that branch's critical
height locally. Its spatial differential, compared with
\eqref{eq:weighted-inverse}, identifies the selected inverse
with the winning branch. In particular $z_0\notin S_{t_\pm,\eta}$.
It also lies outside $\mathcal O_{t_\pm}(\eta,B)$, since all points of these
sets belong to the removed transport segments.

The continuity of the fixed-ratio limit sets \cite[Corollary~3.5]{HS2}
therefore gives a disk $D$ about $z_0$ disjoint from both supports
$S_{t_\pm,\eta_j}$ for all sufficiently large $j$. We also have
$\mu_{t_\pm,\eta_j}\Longrightarrow\mu_{t_\pm,\eta}$.
To justify this additional weak assertion, fix either endpoint and,
for each $j$, choose a positive integer polynomial approximation of
$\eta_j$ with exactly its sites and degree at least $j$. Choose it so
that both its input law and its rounded derivative law are within
$1/j$ of $\eta_j$ and $\mu_{t_\pm,\eta_j}$, respectively, in a metric
for weak convergence. The integer approximation limits permit these
simultaneous choices. The weak law in
Proposition~\ref{prop:finite-location} identifies the diagonal
derivative limit as $\mu_{t_\pm,\eta}$, proving the assertion.

On a smaller closed disk, these weak limits and the common support
gaps give uniform convergence of the Cauchy transforms
$C_{\mu_{t_\pm,\eta_j}}$.
The limiting transforms do not vanish there. Formula
\eqref{eq:weighted-inverse} consequently gives uniform convergence
of the selected inverses at the two ratios $t_\pm$. They remain in fixed, disjoint
neighborhoods of the two saddles at $(z_0,t_0)$, uniformly for $z\in D$.

Choose a small disk about one of these two saddles such that a larger concentric
closed disk contains neither another limiting saddle nor a point of
$B$. The intervening closed annulus contains no solution of
$(w-z_0)C_\eta(w)=t_0$. Uniform Cauchy convergence on this compact
annulus, which avoids $B$, excludes every solution of
\begin{equation}\label{eq:perturbed-saddle-confinement}
 (w-z)C_{\eta_j}(w)=t
\end{equation}
there, uniformly for nearby $z,t$ and sufficiently large $j$.
Take $t_\pm$ close enough and shrink $D$ so that this holds throughout
$D\times[t_-,t_+]$, while the selected endpoints lie on opposite
sides of the annulus. This allows arbitrarily many new saddles near
$B$.

If some $z\in D$ were outside $S_{\eta_j}$, it would be outside every
$S_{t,\eta_j}$ on this interval. Its selected inverse would
then vary continuously with $t$, by
Proposition~\ref{prop:weighted-facts}, and solve
\eqref{eq:perturbed-saddle-confinement}. A continuous path between
the selected endpoints must cross the excluded annulus, by the
intermediate value theorem applied to its distance from the center.
This contradiction proves $D\subset S_{\eta_j}$ for all sufficiently
large $j$.
\end{proof}

\begin{theorem}[Qualitative area universality for finite root locations]
\label{thm:area-universality}
Under the assumptions of Lemma~\ref{lem:regular-interior-persistence},
\begin{equation}\label{eq:area-universality}
 \|\ind_{S_{\eta_j}}-\ind_{S_\eta}\|_{L^1(\C)}\longrightarrow0.
\end{equation}
In particular $\Area(S_{\eta_j})\to\Area(S_\eta)$, including when
positive input masses tend to zero or the number of distinct input
roots tends to infinity.

For polynomial arrays $G_j$ satisfying \eqref{eq:finite-location-input},
\begin{equation}\label{eq:literal-area-universality}
 \|\ind_{\Sh_{G_j}}-\ind_{S_\eta}\|_{L^1(\C)}
 =\|\ind_{\AllSh_{G_j}}-\ind_{S_\eta}\|_{L^1(\C)}
 \longrightarrow0.
\end{equation}
Their Hausdorff limit can be the larger set $K_{\eta,B}$; its added
segments do not change this area limit.
\end{theorem}
\begin{proof}
Since $K_{\eta_j,\supp\eta_j}=S_{\eta_j}$,
Proposition~\ref{prop:finite-location} gives
$S_{\eta_j}\to K_{\eta,B}$. Every point outside this last compact
set is eventually excluded. The difference
$K_{\eta,B}\setminus S_\eta$ has area zero by
Proposition~\ref{prop:finite-location}, as does $\partial S_\eta$
by Proposition~\ref{prop:weighted-tame}.
Lemma~\ref{lem:regular-interior-persistence} gives eventual inclusion
at almost every point of the interior. Therefore
$\ind_{S_{\eta_j}}\to\ind_{S_\eta}$ almost everywhere.
All these sets lie in a common bounded disk, so dominated convergence
proves \eqref{eq:area-universality}.
For polynomials, use $\eta_j=N_j^{-1}\divz G_j$ and
$S_{\eta_j}=\Sh_{G_j}$. The full shadow and the core differ by
at most countably many points, by Proposition~\ref{prop:global-facts}; their
indicators agree almost everywhere for each $j$. This proves the
identity and the convergence in \eqref{eq:literal-area-universality}.
\end{proof}

\begin{corollary}[Area continuity with disappearing atoms]
\label{cor:weak-area-continuity}
Let $\eta_j\Longrightarrow\eta$ be probability measures, each supported
on at most $s$ points in one fixed compact set. Then
$\ind_{S_{\eta_j}}\to\ind_{S_\eta}$ in $L^1(\C)$.
\end{corollary}
\begin{proof}
Every subsequence of the supports has a further Hausdorff-convergent
subsequence with limit $B$ of cardinality at most $s$. Weak
convergence implies $\supp\eta\subset B$. Apply
Theorem~\ref{thm:area-universality} to this further subsequence.
Every subsequence thus has a further subsequence with the stated
$L^1$ limit, which proves the assertion for the whole sequence.
\end{proof}

The qualitative conclusion permits unbounded numbers of distinct input
roots; the constants in Theorem~\ref{thm:area-stability} depend on a fixed
bound for that number.

\appendix
\section{Algebraic selection at rational ratios}\label{app:rational-selection}
For rational differentiation ratios, exponentiation makes the separation
problem itself semialgebraic. Real quantifier elimination---elimination
of quantified real variables from polynomial equalities and
inequalities---then gives an exact terminating procedure on real
algebraic data \cite{BCR}.

\begin{proposition}[Algebraic selection at rational ratios]\label{prop:rational-alpha}
Fix a rational ratio $\alpha=m_0/n_0\in(0,d)$, where $m_0,n_0$ are positive integers. The
support $\FixSh_{\alpha,P}$ is semialgebraic in the observation plane, uniformly in the
coefficients of $P$. For algebraic coefficients, a semialgebraic formula for the entire
support is computable, including its singular points and root-atom thresholds. At an
algebraic non-root observation point, a terminating exact real-algebraic connectivity procedure
also determines the barrier height after exponentiation.
\end{proposition}
\begin{proof}
Multiplying the landscape by $2n_0$ and exponentiating gives the semialgebraic family of
superlevel sets
\[
 |P(w)|^{2n_0}\ge\lambda|w-z|^{2m_0},
 \qquad \lambda>0,
\]
with the marked ends adjoined. For $z\notin Z(P)$, define
\[
 \mathfrak B_{P,\alpha}(z):=\exp\bigl(2n_0\Blev_P(z,\alpha)\bigr).
\]
Finite preconnected covers compatible with the parameters, tests for separation
by relative closures, and unboundedness tests give a first-order semialgebraic predicate
$J(z,\lambda)$ for connectivity of the observation point and infinity in this
compactified superlevel set. For $z\notin Z(P)$, the definition of the barrier gives
$J(z,\lambda)\Longleftrightarrow\lambda\le\mathfrak B_{P,\alpha}(z)$ for every
$\lambda>0$. Thus the graph $b=\mathfrak B_{P,\alpha}(z)$ is defined uniformly in the
coefficients by
\[
 b>0,\qquad
 \forall\lambda>0\quad\bigl(J(z,\lambda)\Longleftrightarrow\lambda\le b\bigr).
\]
For algebraic coefficients and an algebraic observation point, quantifier elimination gives
a polynomial
formula for this graph. Complete real-algebraic sign sampling guarantees an algebraic
witness, so exact algebraic search selects its unique positive value $b$.
This is a terminating exact connectivity procedure; coincident critical
heights require no separate genericity assumption.

Extend $\mathfrak B_{P,\alpha}$ by one at the roots to obtain a total positive
semialgebraic graph. This leaves every non-root germ unchanged. On a root-free disk,
$\Blev_P=(2n_0)^{-1}\log\mathfrak B_{P,\alpha}$ is harmonic exactly when
the positive function is $C^2$ there and
\[
 \mathfrak B_{P,\alpha}\Delta\mathfrak B_{P,\alpha}
 -|\nabla_z\mathfrak B_{P,\alpha}|^2=0
\]
throughout that disk. Differentiability, derivatives and existence of such
a disk have first-order real-algebraic limit formulas. By
\eqref{eq:barrier-potential}, their complement in the non-root plane is
precisely the non-root support, including its singular points. Take its
closure and adjoin exactly the roots with positive residual atom mass,
namely $\alpha<\operatorname{ord}_a(P)$, by
Proposition~\ref{prop:fixed-facts}. A zero-mass root outside that closure
would be an isolated support point of zero mass. These operations are
uniformly semialgebraic and effective over the real algebraic numbers.
\end{proof}

\section{Transient removability and capacity}
\label{subsec:transient-capacity}

For a probability measure $\nu$ carried by a compact set $E\subset\C$, put
\[
 \mathcal E_{\log}(\nu)=\iint\log\frac1{|z-w|}\dd\nu(z)\dd\nu(w),\qquad
 \operatorname{cap}(E)=\exp\bigl(-\inf_{\nu}\mathcal E_{\log}(\nu)\bigr),
\]
where the infimum is taken over all such measures.
A minimizer of this energy is an \emph{equilibrium measure} of $E$.
A set is
\emph{polar} if it is contained in the set where a nonconstant
subharmonic function equals $-\infty$. For these standard
potential-theoretic notions, see \cite[Chapters~3 and~5]{Ransford}.

The \emph{analytic capacity} of a compact set $E$ is
\[
 \operatorname{cap}_{\mathrm{an}}(E)=\sup\{|f'(\infty)|:
 f\text{ holomorphic on }\widehat\C\setminus E,\quad
 |f|\le1,\ f(\infty)=0\},
\]
where $f'(\infty)=\lim_{z\to\infty}zf(z)$. This extremal definition
will be used only through bounded removability.

\begin{theorem}[Transient removability and unchanged capacity]
\label{thm:transient-removability}
For every polynomial with at least two distinct roots:
\begin{enumerate}[label=\textup{(\roman*)},leftmargin=2.2em]
\item Every bounded holomorphic function on
$\C\setminus\AllSh_P$ extends uniquely to a bounded holomorphic
function on $\C\setminus\Sh_P$, with the same supremum norm.
The same assertion holds for bounded harmonic functions.
\item The transient set is polar, and
\[
 \operatorname{cap}(\AllSh_P)=\operatorname{cap}(\Sh_P)>0.
\]
The two compact sets have the same equilibrium measure. They also
have the same analytic capacity:
\[
 \operatorname{cap}_{\mathrm{an}}(\AllSh_P)=\operatorname{cap}_{\mathrm{an}}(\Sh_P).
\]
\end{enumerate}
\end{theorem}
\begin{proof}
By Proposition~\ref{prop:global-facts}, $T_P$ is locally finite in
$\C\setminus\Sh_P$. Each transient is an isolated puncture, across
which a bounded holomorphic or bounded harmonic function extends with
the same bound. Local finiteness makes the extensions compatible, and
continuity gives uniqueness. This proves \textup{(i)}.

If $T_P$ is nonempty, enumerate its points as $t_j$ and choose $b_j>0$
with $\sum_jb_j=1$. Since $T_P$ is bounded, the potential
$\sum_jb_j\log|z-t_j|$ is subharmonic, finite outside a sufficiently
large disk, and $-\infty$ at every $t_j$. Hence $T_P$ is polar.
A finite-energy probability measure on a compact set has no atoms:
the logarithmic kernel is bounded below, while a positive atom gives
infinite energy on the diagonal. It therefore assigns zero mass to the
countable set $T_P$. Thus the finite-energy probability measures carried
by $\AllSh_P$ and by $\Sh_P$ coincide, as do their energy infima and
minimizers. The core is a nondegenerate continuum and consequently has
positive logarithmic capacity \cite[Chapter~5]{Ransford}.
Finally, the extensions in \textup{(i)} identify the admissible functions
in the two analytic-capacity problems, preserving their bounds and
coefficients at infinity. Their analytic capacities are therefore equal.
\end{proof}

The equilibrium measure is different from $\rho_P$: the root atoms
of $\rho_P$, of masses $p_j^2$, give it infinite logarithmic energy,
whereas the equilibrium measure has finite energy. Thus transients
affect neither energy minimization nor bounded analytic continuation.

\end{document}